\documentclass{amsart}
\usepackage{graphicx}
\usepackage{amssymb}
\usepackage{graphicx}
\usepackage{amsfonts}
\usepackage{amsmath}
\usepackage{amsthm}
\usepackage{esint}
\usepackage{mathabx}
\usepackage{mathtools}
\usepackage{subcaption}
\usepackage{nccmath}
\usepackage{enumitem}
\usepackage[colorlinks]{hyperref}

\newtheorem{thm}{Theorem}[section]
\newtheorem{cor}[thm]{Corollary}
\newtheorem{pbm}[thm]{Problem}
\newtheorem{lem}[thm]{Lemma}
\newtheorem{prop}[thm]{Proposition}
\theoremstyle{definition}
\newtheorem{defn}[thm]{Definition}
\theoremstyle{remark}
\newtheorem{rem}[thm]{Remark}
\theoremstyle{remark}
\newtheorem{ex}[thm]{Example}

\renewcommand{\chi}{{ 1}}
\renewcommand{\H}{{ \mathcal H}}

\newcommand{\R}{\mathbb R}

\newcommand{\C}{\mathbb C}
\newcommand{\Hol}{\operatorname{Hol}}
\newcommand{\Id}{\operatorname{Id}}
\newcommand{\Rea}{\operatorname{Re}}
\newcommand{\Ima}{\operatorname{Im}}

\newcommand{\D}{\mathbb{D}}
\newcommand{\B}{\mathcal{B}}
\newcommand{\A}{\mathcal{A}}
\newcommand{\Lmin}{\mathcal L}
\newcommand{\Var}{\operatorname{Var}}
\newcommand{\T}{\mathbb{T}}
\newcommand{\sgn}{\operatorname{sgn}}
\newcommand{\Aw}{\mathcal{A}_w}
\newcommand{\Dm}{D_{\text{mid}}}
\newcommand{\var}{\text{Var}}
\numberwithin{equation}{section}

\begin{document}

\title[Planar Optimal Transport]{On the Planar Optimal Transport Problem}

\author{Arghya Rakshit}
\address{Department of Mathematics, University of Toronto}
\email{\tt arghya.rakshit@utoronto.ca}
\author{Aranya Sen}
\address{Department of Mathematics, UC Irvine}
\email{\tt aranyas1@uci.edu}
\subjclass[2020]{
35J96, 
49Q22. 
}

\keywords{Optimal transportation, Monge-Amp\`{e}re equations}
\date{\today}


\begin{abstract}
We develop a complex analytic framework for the study of quadratic cost planar optimal transport with uniform densities. Our approach is based on rotation of gradient graphs. Using this framework, we carry out the construction, characterization, and approximation of optimal transport maps that are non-degenerate up to the boundary. We also identify obstructions for continuous optimal transport and construct discontinuous optimal transport maps with prescribed analytic singular sets. Using conformal harmonic maps, we also introduce a boundary variational principle à la Douglas for the planar optimal transport problem. 
\end{abstract}
\maketitle
\section{Introduction}
 
Optimal transport asks how to move one distribution of mass to another while minimizing a total cost function. This problem dates back to Monge's work in the eighteenth century \cite{Monge1781}. About a century and a half later, Kantorovich introduced a relaxation in terms of transport plans \cite{K1,K2}. The resulting problem is an
infinite-dimensional linear program.  The theory is simpler for quadratic cost (cost of  moving a unit mass from one point to another is the square of the distance between them), mirroring the relative simplicity of the Hilbert space $L^2$ among the Banach spaces $L^p$, $p>1$.
Under this quadratic cost assumption, for probability measures with finite second moments, if the source measure is absolutely continuous with respect to Lebesgue measure, Brenier in \cite{Brenier} proved there exists a unique optimal transport map, up to sets of source measure zero. Moreover, this map is the gradient of a convex function. See McCann \cite{McCann1995} for a more general existence theory. Unless otherwise specified we have quadratic cost optimal transport in this paper.\\

The interior regularity of optimal transport maps follows from the fundamental work of Caffarelli \cite{Caffarelli1990, Caffarelli1991, Caf1}. Ma, Trudinger, and Wang established interior regularity for more general costs under a structural condition now known as the Ma-Trudinger-Wang condition, together with suitable assumptions
on the densities and domains \cite{MTW}.
For boundary regularity of optimal transport maps when the source and target are convex, see Caffarelli \cite{Caf1, Caf2}, Savin-Yu \cite{SavinYu2020}, Chen-Liu-Wang \cite{CLW2}, Collins-Tong \cite{CT1}, Collins-Firester-Tong \cite{CFT1}, and Chen-Li-Wang \cite{CLW2026}. The boundary regularity theory for Optimal Transport between general domains is far from being well understood. For boundary regularity of nonconvex domain see Mooney-Rakshit \cite{MR2}, Hu-Li \cite{HuLi2026},  Andriyanova-Chen \cite{AC1}. The discontinuity sets of optimal transport maps have been studied by Figalli~\cite{Figalli2010}, Chen, Li, and Liu~\cite{ChenLiLiu2023}, Chen, Li, Li, Liu, and Wu~\cite{ChenEtAlPolyhedral2026,ChenEtAlDimension2026}, Kitagawa and McCann \cite{KitagawaMcCann2019}, and McCann and Rakshit~\cite{McCannRakshit}; see also
Rakshit and Sen~\cite{RS1} for related examples of singular Monge--Amp\`ere potentials. Stability results beyond the convex setting were studied by Chen-Figalli \cite{CF}, Chen-Liu-Wang  \cite{CLW1}, Jhaveri \cite{Jhaveri}, Miura-Otto \cite{MO1}.\\

We work with uniform densities on bounded planar domains $\Omega=\{(x,y)\}$ and $\Lambda=\{(p,q)\}$ of equal area. We call an optimal transport map, $T$, smooth non degenerate if $DT>0$ on the closure of the source domain. We are interested in how boundary geometry governs transport between domains that need not be convex. Nonconvex targets can admit smooth transport (see \cite{McCannRakshit}), but they can also give rise to discontinuities (see \cite{Caf1}). We seek to understand how these different behaviors are reflected in the boundary data. This leads to several questions:
\begin{enumerate}[label=\Roman*.]
    \item Which nonconvex targets admit smooth transport from a given source? 
    \item Can discontinuity be detected from boundary of source and target domains?
    \item Under what class of perturbations is the regularity of optimal transport preserved?
    \item Can one (explicitly) construct an optimal transport map whose discontinuity set can be prescribed?
\end{enumerate}
These questions have not been addressed in their full generality. In this paper we provide answers to the above questions. In doing so, we develop a new viewpoint for the study of planar optimal transport to non-convex target domains. This approach allows us to obtain several sharp results. \\

Our starting point is partial Legendre transform. Partial Legendre transform, and other transforms like Lewy transform, can be seen as rotations of gradient graph. See Figalli \cite{Figalli_book} for an overview of partial Legendre transform and see Mooney-Shankar \cite{MS} for the rotation of gradient graph viewpoint. These transforms have been extremely successful in establishing local regularity theory (see Daskalopoulos-Savin \cite{DS}, Mooney-Rakshit \cite{MR2}, Mooney-Sen \cite{MSen}) as well as global rigidity (see J\"orgens \cite{J}, Jin-Xiong \cite{JX}). However, the implicit nature of these transforms has limited their applicability to boundary value problems on bounded domains. These transforms allow us to study the mass preservation equation written in terms of potential $u=u(x,y)$ as $\det D^2u=1$ in $\Omega$ as a harmonic function $h$ on a mixed source target domain $D=\{(x,q)\}$ locally. The derivative of the harmonic function gives a holomorphic function $f$ defined by $f(z)=p+iy$ where $z=x+iq$. We call $f$ the \emph{holomorphic generator} of the optimal transport map and $D$  the \emph{mixed domain} arising from the optimal transport.

\begin{figure}[ht]
    \centering
    \includegraphics[width=0.85\textwidth]{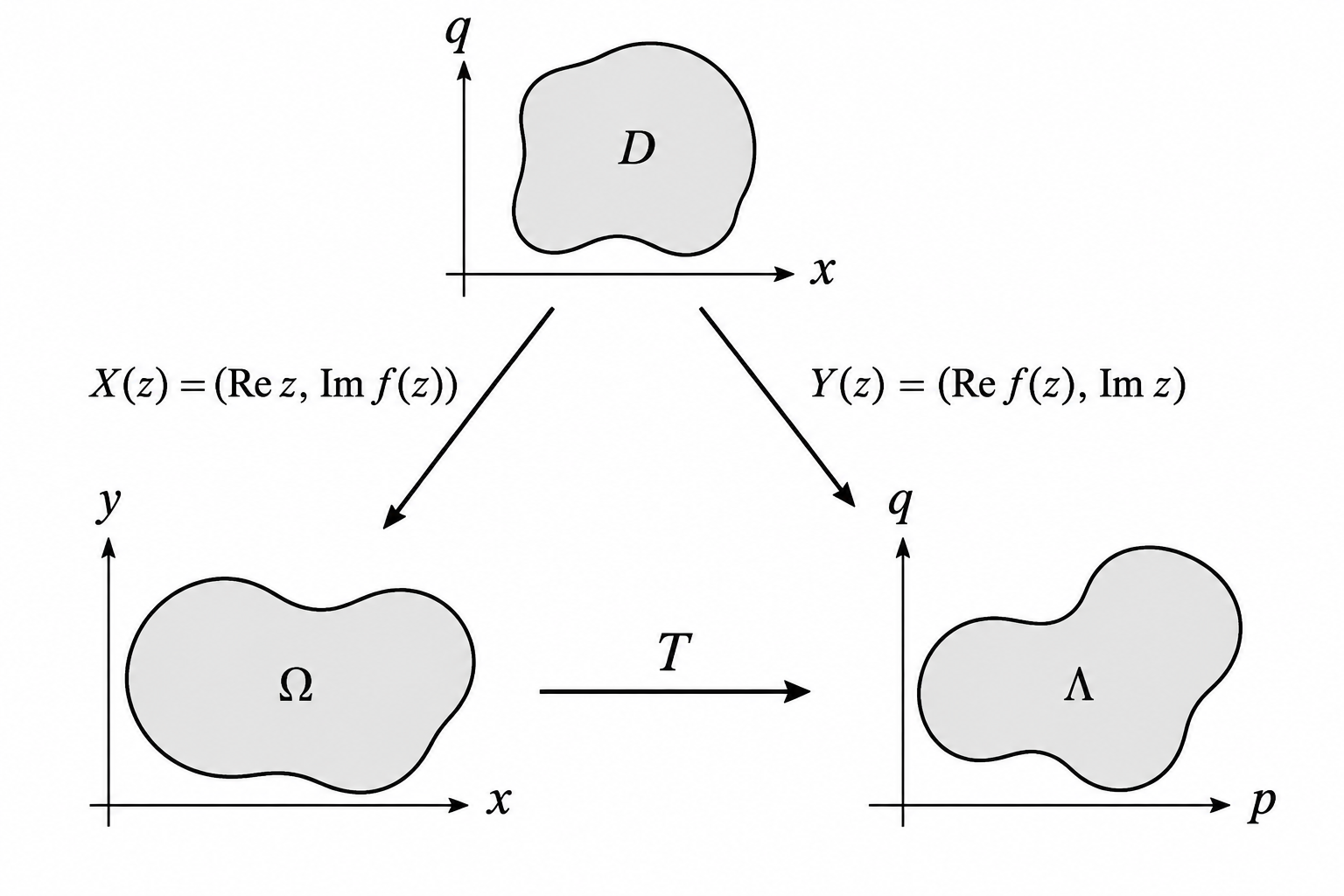}
    \caption{The holomorphic generator $f$ and the mixed domain $D$ and the associated source and target domains.}
    \label{fig:generator}
\end{figure}

  We develop the converse construction. Starting from a holomorphic function $f$ defined in a mixed domain $D$, we recover the source, target, and the transport map. We write down explicit conditions for the map to be globally defined and optimal. \\
Given a holomorphic function $f$ on a planar domain $D$, define
$$
X(z)=(\Rea z,\Ima f(z)),
\qquad
Y(z)=(\Rea f(z),\Ima z).
$$
See Figure \ref{fig:generator} for an overview of the construction described. These assignments determine the source $\Omega=X(D)$, the target $\Lambda=Y(D)$, and the transport through
$$
T(X(z))=Y(z).
$$
\newpage
Building upon this construction, the following theorem gives sufficient conditions on $f$ and $D$ for this construction to produce an optimal transport map.

\begin{thm}\label{thm:global-generator-intro}
Let $D\subset\mathbb C$ be connected and open, and let $f$ be holomorphic on $D$ with $\Rea f'>0$. Suppose the coordinate maps $X$ and $Y$ are injective and have bounded images $\Omega$ and $\Lambda$. Assume that $\Rea\oint_\gamma f(z)\,dz=0$ for every closed piecewise smooth curve $\gamma\subset D$, and that the supporting-plane condition \eqref{eq:support} holds.

Then $T=Y\circ X^{-1}$ is the unique smooth optimal transport from $\Omega$ onto $\Lambda$. Its potential extends to a convex function $u$ on $\mathbb R^2$ and satisfies $D^2u>0$ and $\det D^2u=1$ in $\Omega$.
\end{thm}
The real period condition ensures $u$ is globally defined up to an additive constant.  The supporting-plane condition ensures global convexity, which is an additional requirement when the source is nonconvex.

Explicit optimal transport maps are difficult to construct. Our approach above gives several explicit examples. Here we mention a few examples when one of the domains is nonconvex. One map transports the region between two concentric ellipses to an ellipse. The two source ellipses are not dilates of one another. The map extends continuously to the inner ellipse and collapses it onto an interior line segment. Another map takes a triangle onto a nonconvex parabolic strip. It is real analytic in the interior and up to every boundary point except one vertex. At that vertex, the limiting values fill an entire parabolic arc, so the map has no continuous extension there. Higher-degree polynomial generators also give smooth transports from the unit disk to nonconvex targets. The target boundaries can be arbitrarily close to the unit circle in $C^1$ topology, while their signed curvature becomes arbitrarily negative at some points.

Building upon the holomorphic generator viewpoint for optimal transport, we characterize domains that admit smooth nondegenerate optimal transport. For curves $\alpha(t) = (x(t),y(t))$ and $\beta(t)=(p(t),q(t))$ define
\begin{equation}\label{eq:general-data_intro}
 \zeta_{\alpha,\beta}(t)=x(t)+iq(t),\qquad \eta_{\alpha,\beta}(t)=p(t)+iy(t).
\end{equation}

\begin{thm}\label{thm:intro_thm1}
Let $\Omega$ be a smooth bounded convex domain with a positively oriented boundary parametrization $\alpha(t)$ and let $\Lambda$ be a bounded smooth Jordan domain. There is a smooth optimal transport map of the closures from $\Omega$ to $\Lambda$ if and only if 
$\partial\Lambda$ admits a positively oriented smooth regular parametrization, $\beta(t)$, for which $\zeta$ in \eqref{eq:general-data_intro} is a positively oriented smooth regular Jordan curve and
\begin{equation}\label{eq:intro_moment_conditions}
\int_0^{2\pi}\eta_{\alpha,\beta}(t)\zeta^j_{\alpha,\beta}(t)\zeta_{\alpha,\beta}'(t)\,dt=0 \quad(j=0,1,2,\ldots).
\end{equation}
\end{thm}

See Remark \ref{rem:charac_non_convex} for the characterization when $\Omega$ may be nonconvex and Remark \ref{rem:charac_low_reg} for a lower regularity version of the above theorem.

Theorem \ref{thm:intro_thm1} reduces the problem of finding non degenerate smooth optimal transport to finding the correct boundary parameterization. By reducing the problem to finding harmonic diffeomorphisms, we find all possible boundary parameterizations the give rise to smooth non degenerate optimal transport from the unit ball.

\begin{thm}\label{thm:intro_thm5}
    Let $\Lambda$ be a smooth Jordan domain. There exists a smooth non degenerate optimal transport map $T:B_1(0)\to \Lambda$ if and only if $\partial\Lambda$ admits a parameterization from the set 
    $$\{J\H (\alpha)+c\;|\; \alpha \text{ is a smooth circle diffeomorphism }, \alpha'>0,\;c\in\R^2\}.$$
    Here $J$ is counterclockwise rotation by $\pi/2$ and $\H$ is the Hilbert transform on the circle.
\end{thm}

In his seminal work on the Plateau problem, Douglas \cite{Douglas} reduced the problem of finding minimal surfaces to finding the minimizers of 
\[
    \mathcal{D}(g)=\frac{1}{4\pi}
    \int_0^{2\pi}\int_0^{2\pi}
    \frac{|g(t)-g(s)|^2}{|e^{it}-e^{is}|^2}\,dt\,ds
    =\frac12\int_{\mathbb T}g\cdot Kg\,dt
\]
over parametrizations of a prescribed Jordan curve, where
$K=(-\partial_{tt})^{1/2}$. We also present a boundary variational approach based on the following functional
\[
    \mathcal{O}_\Lambda(g)=\int_{\mathbb T}|Kg|\,dt.
\]
In our framework, the isoperimetric inequality plays the role
of the classical energy-area inequality, and its equality case
characterizes continuous planar optimal transport from the unit disk (see Theorem \ref{thm:infi_charac}). We have the following equivalence for smooth non degenerate optimal transports based on the functional $\mathcal{O}_{\Lambda}$.
\begin{thm}\label{thm:intro_thm6}
    Let $\Lambda$ be a smooth Jordan domain. Then $\Lambda$ admits smooth non degenerate optimal transport from the unit ball $B_1(0)$ if and only if there exists a positively oriented smooth parametrization of $\partial\Lambda$, $g$, such that 
$$
 Kg \neq 0\quad \text{ and }\quad g'\cdot K\left(\frac{Kg}{|Kg|}\right)=0. 
$$
\end{thm}

The nonlocal equation in the above theorem is the Euler-Lagrange equation of the functional $\mathcal{O}_\Lambda$. To prove Theorem \ref{thm:intro_thm6} we actually show that for smooth enough parameterizations (atleast $C^{2,\alpha}$ for some $0<\alpha<1$), stationarity implies minimality. Here stationarity is with respect to inner reparameterizations. Given parameterizations from Theorem \ref{thm:intro_thm5} and Theorem \ref{thm:intro_thm6}, the corresponding optimal transport maps can also be constructed immediately by considering certain harmonic extensions arising from boundary parameterizations.

 By fixing boundary parameterizations of both the domains, one can reduce the problem of finding continuous optimal transport between domains to finding a monotone mapping on the boundary. Using the moment conditions \eqref{eq:intro_moment_conditions} and Lewy transform we are able to show the following family of tests to detect discontinuous optimal transport maps.
\begin{thm}\label{thm:intro_thm2}
     Let $\Lambda$ be a bounded smooth Jordan domain of area $\pi$. Let $\Gamma(s)$ be an arc length parametrization of $\partial\Lambda$ with arclength $l = \mathcal{H}^1(\partial\Lambda)$. Fix a continuous tangent angle lift $\beta:\R\to\R$ such that $\Gamma'(s) = e^{i\beta(s)}$. Define $\Aw$ as the set of continuous nondecreasing functions $\theta:\R\to\R$ satisfying
\begin{equation}
 \theta(s+\ell)=\theta(s)+2\pi,\quad
 \max_{r\in[s-\ell,s]}\big(\beta(r)-\pi\big)\le\theta(s)\le \min_{r\in[s,s+\ell]}\beta(r).
\end{equation}
If the class $\Aw$ is empty or for some integer $k\ge0$,
\begin{equation}
  M_k[\theta] = \int_0^\ell \eta_{\theta,\Gamma}\,
 \zeta_{\theta,\Gamma}^k\, d\zeta_{\theta,\Gamma}\ne0\qquad\text{for every }\theta\in\Aw,
\end{equation}
then the optimal transport map from $B_1(0)$ to $\Lambda$ is discontinuous. 
\end{thm}

For a general convex source see Remark \ref{rem:perturb_gen_convex}. A sufficient condition for non degenerate optimal transport was established by Wolfson in \cite{Wolfson} and later refined by Urbas in \cite{Urbas}. It was shown by Chodosh-Jain-Lindsey-Panchev-Rubinstein in \cite{Chodosh} that these conditions actually lead to discontinuous optimal transport maps. Their conditions are equivalent to the non-existence of an orientation preserving boundary correspondence satisfying the classical obliqueness inequality. The holomorphic moment obstructions given above impose additional compatibility conditions beyond this geometric requirement. Example \ref{ex:example_obstruction}  gives smooth, strictly star-shaped targets for which the radial boundary correspondence is strictly oblique, yet the moment obstruction forces the optimal transport from the disk to be discontinuous.

Next we establish a perturbation result. By using Riemann mapping theorem to fix the mixed domain and a Schauder fixed point argument, we are able to show the stability of optimal transport maps under $C^{1,\alpha}$ perturbations of the boundary. The construction and boundary perturbation arguments do not use estimates on the geometry of Monge--Amp\`ere sections. Instead, we work in mixed coordinates, where the determinant-one equation becomes the Laplace equation. This gives a direct approach to regularity within the configurations considered here.
\begin{thm}\label{thm:intro_thm3}
Let $\Omega_0,\Lambda_0$ be bounded $C^{2,\alpha}$ domains of same area and assume $\Omega_0$ convex. Assume there is a $C^{2,\alpha}(\bar\Omega_0)$ optimal transport map $T_0:\Omega_0\to\Lambda_0$. Consider $C^{1,\alpha}$ domains $\Omega,\Lambda$ whose boundaries are normal graphs, denoted by $h_1$ and $h_2$, over $\partial\Omega_0,\partial\Lambda_0$ respectively. Further assume $\Omega$ is convex. If for $i=1,2$, $\|h_i\|_{C^{1,\alpha}}$ are small while ensuring $|\Omega| = |\Lambda|$, then there exists an optimal transport map $T:\Omega\to\Lambda$ such that for any $\beta<\alpha$, $T\in C^{1,\beta}(\bar\Omega)$.
\end{thm}

It was known by Jhaveri \cite{Jhaveri} that mere $C^1$ perturbations of the boundaries can lead to discontinuous optimal transport. In fact it was shown that given an $\epsilon$, one can actually find small $\alpha'(\epsilon)$ small such that the said $C^{1}$ perturbation can be made small in the $C^{1,\alpha'(\epsilon)}$ norm. Near the identity map, an $\epsilon$ regularity result for $C^{1,\alpha}$ perturbations was obtained by Miura and Otto in \cite{MO1}. Therefore as discussed in \cite{MO1}, the above theorem is sharp. Using the perturbation theorem above we also have the following purely algebraic characterization of smooth nondegenerate optimal transport maps from $B_1(0)$.

\begin{cor}\label{cor:intro_polynomial}
    Nash functions that are optimal transport maps from $B_1(0)$ onto smooth semialgebraic domains whose boundaries lie on real algebraic curves are dense (in the $C^\infty$ topology) among smooth nondegenerate optimal transport maps from $B_1(0)$.
\end{cor}
Note that in the above corollary the topology can be made weaker by considering the $C^{1,\alpha}$ topology for any $0<\alpha<1$. It might be of independent interest to find algebro-geometric conditions guaranteeing smooth optimal transport onto semialgebraic targets.

To show robustness of our approach we also explicitly construct arbitrary real analytic singular sets arising from optimal transport between convex domains and annular domains.

\begin{thm}\label{thm:intro_thm4}
    Let $\Psi$ be a real analytic function defined on an open neighbourhood of the interval $[-1,1]$. Then there exists open uniformly convex sets $\Omega,K_{\text{in}}$ and $K_{\text{out}}$, such that for the optimal transport map
    $$T:\Omega\to K_{\text{out}}\setminus \overline{K_{\text{in}}},$$
    the singular set of the map is precisely the graph of $\Psi$ on $(-1,1)$.
\end{thm}

In the theorem above, the uniformly convex domains $\Omega,K_{\text{in}}$ and $K_{\text{out}}$ are explicitly constructed as well as the holomorphic generator of the map $T$. See Theorem \ref{thm:general_fbp} for the complete description.

We now discuss an important application to the study of minimal Lagrangian diffeomorphisms.

\subsection{Minimal Lagrangian Diffeomorphisms} In \cite{Wolfson} Wolfson studied the question of the existence of smooth (upto the boundary) minimal Lagrangian diffeomorphisms in two dimensions. The question was originally proposed by Schoen \cite{Schoen} in much more general settings. The exact problem statement is:
\begin{pbm}\label{pbm:min_lagrange}
    Let $D_1$ and $D_2$ be simply connected, closed, bounded domains in $\R^2$ with smooth boundaries and equal areas. Find an area preserving diffeomorphism $\psi: D_1\to D_2$ smooth up to the boundary such that the graph of $\psi$, $\text{graph}(\psi)$, is a minimal surface in $\R^4\simeq \R^2\times \R^2$.
\end{pbm}

Wolfson obtained a sufficient condition for non-existence of minimal Lagrangian diffeomorphisms smooth upto the boundary using ideas of Lagrangian angles and Maslov forms. In two dimensions given two domains $D_1$ and $D_2$, one can show that Problem \ref{pbm:min_lagrange} is essentially equivalent to the following second boundary value problem for the Monge-Ampere equation,
\begin{equation}\label{eq:second_BVP}
\det(D^2w) = 1 \text{ in } D_1,\qquad \nabla w(D_1) = D_2,\qquad w\in C^{\infty}(\bar D_1).
\end{equation}
In two dimensions, existence of solutions to the second boundary value problem for general Monge-Ampere equations was obtained by Delanoë in \cite{Delanoe} for uniformly convex domains $D_1,D_2$ and in higher dimensions by Caffarelli in \cite{Caf2} and Urbas in \cite{U2}. Urbas in \cite{Urbas} noted that if a certain obliqueness condition can't be obtained for given domains $D_1,D_2$, then a solution to \eqref{eq:second_BVP} cannot exist. Chodosh-Jain-Lindsey-Panchev-Rubinstein in \cite{Chodosh} noted that the said obliqueness condition actually leads to discontinuous optimal transport.

We note that Theorem \ref{thm:intro_thm1} analytically characterizes all smooth domain pairs $D_1,D_2$ (with the source domain $D_1$ being convex) that admit a solution to \eqref{eq:second_BVP}. Moreover, Theorem \ref{thm:intro_thm2} gives conditions for non-existence of solutions to \eqref{eq:second_BVP} purely based on the geometry of the boundary $\partial D_2$ that are more restrictive than obliqueness.

The paper is organized as follows. In Section \ref{sec:prelims} we recall the partial Legendre transform and basic facts about harmonic and holomorphic functions. In Section \ref{sec:Holo_generator} we develop the holomorphic generator construction and give sufficient conditions for the resulting map to be a globally defined optimal transport. In Section \ref{sec:examples} we give several explicit examples, including transports to nonconvex targets and examples with boundary singularities.  In Section \ref{sec:Characterization} we characterize smooth optimal transport in terms of the boundary correspondence and the Cauchy moments. In Section \ref{sec:discon} we use these moment conditions to obtain obstructions to continuity and construct examples where the optimal transport is discontinuous. In Section \ref{sec:perturb} we prove the stability of smooth optimal transport under small $C^{1,\alpha}$ perturbations of the domains. In Section \ref{sec:harmonic_maps} we prove general versions of Theorem \ref{thm:intro_thm5} and Theorem \ref{thm:intro_thm6}. In Section \ref{sec:poly}, we prove a general version of Corollary \ref{cor:intro_polynomial} and up to perturbations, give the denseness of a family of relatively easy to construct optimal transport maps. Finally, in Section \ref{sec:free_boundary} we study the free-boundary structure arising from the Laurent example and construct transports whose singular set is a prescribed real-analytic graph.

\section*{Acknowledgment}
The authors would like to thank Professor Connor Mooney for introducing the technique of rotation of gradient graphs. A. Rakshit would like to thank Professor Robert J. McCann, Professor Tristan C. Collins and Professor Freid Tong for several helpful discussions related to optimal transport. A. Sen was partially supported by NSF grant DMS-2143668 of C. Mooney.

\section*{AI disclosure}
ChatGPT Plus was used for generating images, choosing mixed domains in certain examples, and polishing the manuscript. The connection with Banchoff-Pohl inequality in Lemma \ref{lem:iso_charac} was discovered with the help of of ChatGPT Plus. All other mathematical ideas and results are the authors’ own.

\section{Preliminaries}\label{sec:prelims}
\subsection{Partial Legendre transform} For a general background on partial Legendre transform and its relation to Monge-Ampere equations see Figalli \cite{Figalli_book}. Here we recall only the definition and its connection to Laplace equation. The partial Legendre transform of a convex function $u:\Omega\to\R$ is
defined by
$$
 u^*(x,q)=\sup_{y:\,(x,y)\in\Omega}\{yq-u(x,y)\}.
$$
For a smooth convex function, the possible mixed variables are given by
$$
 \Omega^*=\{(x,u_y(x,y)):(x,y)\in\Omega\}.
$$
At a point of $\Omega^*$ where the supremum is attained at a unique
interior point $(x,y)$, we have $q=u_y(x,y)$ and
\begin{equation}\label{eq:PLT}
 u^*(x,q)=yq-u(x,y).
\end{equation}
The condition $u_{yy}>0$ allows us to solve locally for $y=y(x,q)$. Let $u$ be a smooth solution of
\begin{equation}\label{eq:MA}
 \det D^2u=1\quad\text{in }\Omega,\qquad D^2u>0.
\end{equation}
Put $p=u_x$ and $q=u_y$. Differentiating \eqref{eq:PLT} gives
\begin{equation}\label{eq:PLT-derivatives}
 u_x^*=-p,\qquad u_q^*=y,\qquad
 u_{qq}^*=\frac1{u_{yy}},\qquad
 u_{xx}^*=-\frac{\det D^2u}{u_{yy}}.
\end{equation}
Therefore the partial Legendre transform gives
$$
 \Delta u^*=0\quad\text{in }\Omega^*.
$$
 Write $h=u^*$ and identify $\Omega^*$ with the
mixed domain $D\subset\C$ by $z=x+iq$. Since $h$ is harmonic, the function
\begin{equation}\label{eq:PLT-generator}
 f=-h_x+ih_q=p+iy
\end{equation}
is holomorphic with $\Rea f' = -h_{xx} = h_{qq}>0$. In Section \ref{sec:Holo_generator} we give an inverse of this construction based on $f$.

\subsection{Harmonic functions in $2$D and Complex analysis} We record several results from complex analysis that we will use throughout the paper. To fix notation, let $\D$ be the unit disk centred at $0$ on the complex plane. We also denote the set of all Holomorphic functions on an open set $\mathcal{U}$ as $\text{Hol}(\mathcal{U})$.\\ 
The first result concerns the boundary behaviour of conformal mappings from the unit disk. The original results are due to Kellogg in \cite{Kellog} and Warschawski in \cite{Warsch}. For a proof see \cite[Theorem 3.6]{Pomm}.
\begin{thm}[Kellogg-Warschawski Theorem]
    Let $D\subset\C$ be a bounded Jordan domain with $C^{m,\alpha}$ boundary for $m\geq 1$ and $0<\alpha<1$. If $\Phi:\D\to D$ is conformal, then $\Phi$ extends to a $C^{m,\alpha}$ map on $\overline\D$. 
\end{thm}

The following weaker version of the above result was shown by Carath\'eodory in \cite{Cara}.
\begin{thm}[Carath\'eodory's theorem]\label{thm:caratheodory}
Let $D\subset\mathbb{C}$ be a bounded Jordan domain and let
$\Phi:\mathbb{D}\to D$ be a conformal bijection. Then $\Phi$
extends uniquely to a homeomorphism
\[
    \Phi:\overline{\mathbb{D}}\longrightarrow\overline{D}.
\]
In particular, $\Phi|_{\partial\mathbb{D}}$ is an
orientation-preserving homeomorphism onto $\partial D$.
\end{thm}

The following theorem deals with polynomial approximation of holomorphic functions on compact sets. It was first shown by Mergelyan in \cite{Mer}. See \cite[Theorem 20.5]{Rudin} for a modern reference.

\begin{thm}[Mergelyan's Theorem]
    Let $K\subset\C$ be compact and suppose $\C\setminus K$ is connected. For every $f\in C(K)\cap \text{Hol}(K^{\circ})$ there are complex polynomials $p_j$ such that
$$
\|p_j-f\|_{C(K)}\longrightarrow0.
$$

\end{thm}

The next theorems deal with the solvability of various boundary value problems concerning holomorphic and harmonic functions. The following theorem was an important result established in the first solutions of the Riemann-Hilbert problem. It was first discovered by Sokhotski and then independently rediscovered by Plemelj in \cite{Plem}. See \cite[Chapter 7]{Kress} for a modern reference.

\begin{thm}[Sokhotski-Plemelj Theorem]
Let $\Gamma=\partial D$ be a positively oriented $C^{1,\alpha}$ Jordan curve, and let $g\in C^{0,\beta}(\Gamma)$ with $0<\beta<1$. Define the Cauchy integral for $z\notin\Gamma$ by
$$
\mathcal Cg(z)=\frac1{2\pi i}\oint_\Gamma\frac{g(\zeta)}{\zeta-z}\,d\zeta.
$$
At every $z_0\in\Gamma$ the limits from the interior ($+$) and exterior ($-$) exist along non-tangential directions, with
\begin{equation}\label{eq:plemelj}
(\mathcal Cg)^\pm(z_0)
=\pm\frac12g(z_0)+\frac1{2\pi i}\text{p.v.}\oint_\Gamma
\frac{g(\zeta)}{\zeta-z_0}\,d\zeta.
\end{equation}
In particular,
\begin{equation}\label{eq:jump}
(\mathcal Cg)^+(z_0)-(\mathcal Cg)^-(z_0)=g(z_0).
\end{equation}
\end{thm}

The Schwarz Problem is another version of the Riemann-Hilbert problem asking for a prescribed real part of a holomorphic function on the unit disk. For the solvability of the Schwarz problem see \cite{Ahl} for a standard reference.  

\begin{thm}[Schwarz integral]
 Let $\phi\in C(\partial\D;\R)$ and $c\in\R$. There exists a unique $f\in\Hol(\D)$ such that $\Rea f$ extends continuously to $\overline\D$, $\Rea f|_{\partial\D}=\phi$, and $\Ima f(0)=c$. It is given by
\begin{equation}\label{eq:schwarz}
 f(z)=\frac{1}{2\pi}\int_0^{2\pi}
\frac{e^{it}+z}{e^{it}-z}\phi(e^{it})\,d t+ic,\qquad |z|<1.
\end{equation}
 More generally, if $\phi\in C^{m,\alpha}$ for $m\geq 0$ and $0<\alpha<1$, the holomorphic extension belongs to the same class up to the circle.
\end{thm}

The next result deals with the solvability of Neumann boundary value problem for harmonic functions (see \cite{Evans} for a standard reference). Let $\Omega\subset\R^d$ be bounded and connected, with outward normal $\nu$. The Neumann problem is
\begin{equation}\label{eq:neumann_n}
\Delta u=0\quad\text{in }\Omega,
\qquad \partial_\nu u=g\quad\text{on }\partial\Omega.
\end{equation}
\begin{thm}
Suppose $\Omega$ is a bounded simply connected $C^{1,\alpha}$ domain, $0<\alpha<1$, and $g\in C^{0,\alpha}(\partial\Omega)$. The problem \eqref{eq:neumann_n} has a solution $u\in C^{1,\alpha}(\overline\Omega)$ if and only if
\begin{equation}
\int_{\partial\Omega}g\,d\mathcal H^{d-1}=0.
\end{equation}
The solution is unique up to an additive constant. With $\int_\Omega u=0$,
$$
\|u\|_{C^{1,\alpha}(\overline\Omega)}
\le C(\Omega,\alpha)\|g\|_{C^{0,\alpha}(\partial\Omega)}.
$$
\end{thm}

For Harmonic maps we have the following important result first conjectured by Rad\'o in \cite{Rado} and proven by Kneser in \cite{Kneser} and independently rediscovered by Choquet in \cite{Cho}.
\begin{thm}[Rad\'o--Kneser--Choquet theorem]
Let $\Omega\subset\mathbb{R}^2$ be a bounded convex domain.
Suppose that $f:\partial\mathbb{D}\to\partial\Omega$
is a continuous nondecreasing positively oriented degree one parametrization. Then the harmonic extension $F$ belongs to $C^\infty(\mathbb{D})\cap C^0(\overline{\mathbb{D}})$ and is an orientation-preserving diffeomorphism
\[
    F:\mathbb{D}\longrightarrow\Omega.
\]
\end{thm}

By combining it with regularity for harmonic functions and an immediate application of Hopf lemma gives the following stronger version.

\begin{thm}[$C^{k,\alpha}$ Rad\'o--Kneser--Choquet theorem]
Let $k\geq 1$, $0<\alpha<1$, and let
$\Omega\subset\mathbb{R}^2$ be a bounded convex domain with
$C^{k,\alpha}$ boundary. Suppose that
\[
    f:\partial\mathbb{D}\longrightarrow\partial\Omega
\]
is an orientation-preserving $C^{k,\alpha}$ diffeomorphism.
Then its harmonic extension $F$ is an orientation preserving
$C^{k,\alpha}$ diffeomorphism
\[
    F:\overline{\mathbb{D}}\longrightarrow\overline{\Omega},\quad \det(DF)>0 \text{ on }\bar\D.
\]
\end{thm}

\section{Holomorphic Generator for Optimal Transport}\label{sec:Holo_generator}
Throughout, $(x,y)\in\Omega$ are source coordinates and $(p,q)\in\Lambda$ are target coordinates. We study the optimal transport map from $\Omega$ to $\Lambda$. We use the convention
\begin{equation}
 z=x+iq\in D,\qquad f(x+iq)=p+iy.
\end{equation}
We construct and characterize the optimal transport map using a holomorphic function $f$ on a mixed domain $D\subset \C$. 

\subsection{Local Calculation}
Let $D\subset\C$ be connected and open, and let $f$ be holomorphic on $D$, that is, $f\in\Hol(D)$. For $z=x+iq$, define
\begin{equation}\label{eq:assignments}
 X(z)=(x,\Ima f(z)),\;\; Y(z)=(\Rea f(z),q),\;\;\;
 \Omega=X(D),\;\; \Lambda=Y(D).
\end{equation}
When $X$ is injective, we define the map $T$ by 
\begin{equation}
 T\bigl(x,\Ima f(x+iq)\bigr)=\bigl(\Rea f(x+iq),q\bigr).
\end{equation}
Below we give the conditions under which this map is optimal.
Writing $y=\Ima f$ and $p=\Rea f$, the Cauchy--Riemann equations give
\begin{equation}\label{eq:CR}
 y_x=\Ima f',\qquad y_q=\Rea f',\qquad p_x=\Rea f',\qquad p_q=-\Ima f'.
\end{equation}
The derivatives here are taken in the mixed coordinates $(x,q)$. The two maps from $D$ to the source and target have Jacobian $\Rea f'$. If
\begin{equation}\label{eq:positive}
 \Rea f'>0\quad\text{in }D,
\end{equation}
then the source assignment is locally invertible:  Since $$ \partial_q\bigl[\Ima f(x+iq)\bigr]=\Rea f'(x+iq)>0, $$ the implicit function theorem allows us to solve $$ y=\Ima f\bigl(x+iq(x,y)\bigr) $$ locally for $q=q(x,y)$. 
Differentiating this identity with respect to $x$ and $y$ gives $$ 0=\Ima f'(z)+\Rea f'(z)\,q_x, \qquad 1=\Rea f'(z)\,q_y, $$ where $z=x+iq(x,y)$. Hence $$ q_x=-\frac{\Ima f'(z)}{\Rea f'(z)}, \qquad q_y=\frac1{\Rea f'(z)}. $$ Writing $p(x,y)=\Rea f\bigl(x+iq(x,y)\bigr)$, the chain rule gives $$
 \begin{aligned}
 p_x&=\Rea f'(z)-\Ima f'(z)\,q_x =\frac{|f'(z)|^2}{\Rea f'(z)},\\
 p_y&=-\Ima f'(z)\,q_y =-\frac{\Ima f'(z)}{\Rea f'(z)}.
 \end{aligned}
$$ Therefore
\begin{equation}\label{eq:Hessian}
 DT(X(z))=
 \begin{pmatrix}
 p_x&p_y\\
 q_x&q_y
 \end{pmatrix}
 =\frac1{\Rea f'(z)}
 \begin{pmatrix}
 |f'(z)|^2&-\Ima f'(z)\\
 -\Ima f'(z)&1
 \end{pmatrix}>0.
\end{equation}
In particular, $$ \det DT(X(z)) =\frac{|f'(z)|^2-(\Ima f'(z))^2}{(\Rea f'(z))^2}=1. $$ The derivative is symmetric. Locally, $T$ is the gradient of a strictly convex potential and preserves area.
\begin{rem}
Condition \eqref{eq:positive} guarantees local invertibility of both
coordinate assignments. If $D$ is convex, both coordinate assignments (to source and target) are also globally injective. Indeed, for each fixed $x$, the vertical slice $$ \{q\in\R:x+iq\in D\} $$ is an interval, whenever it is nonempty. Since $\partial_q[\Ima f(x+iq)]=\Rea f'(x+iq)>0$, the function $q\mapsto\Ima f(x+iq)$ is strictly increasing on this interval. So this proves its injectivity. Similar argument works for target injectivity.

 If a vertical line meets $D$ in separate intervals, points from different intervals may still give the same source point $(x,y)$. We must check separately that this does not happen. 
\end{rem}

Two global issues are left. Since $\nabla u=(p,q)$ locally, constructing a potential on all of $\Omega$ requires that the integral of $T$ along a path from a fixed base point to $(x,y)$ be independent of the path within $\Omega$. Even if such a potential exists, $D^2u>0$ does not by itself give global convexity when $\Omega$ is nonconvex.

\subsection{Global Potential}
Here we address the issue of existence of a global potential:
\begin{prop}\label{prop:periods}
Assume that $X$ is injective and \eqref{eq:positive} holds. There is a function $u$ on $\Omega$ with $\nabla u=T$ if and only if
\begin{equation}\label{eq:periods}
 \Rea\oint_\gamma f(z) dz=0,
\end{equation}
for every closed piecewise smooth curve $\gamma\subset D$.
When this holds, choose $z_0\in D$ and set
\begin{equation}\label{eq:h-integral}
 h(z)=-\Rea\int_{z_0}^z f(\zeta)d\zeta,\;\;\; u(X(z))=\Ima z\Ima f(z)-h(z).
\end{equation}
These functions are well defined, and $u$ is unique up to a constant.
\end{prop}
\begin{proof}
Condition \eqref{eq:periods} gives that the real part of the integral in \eqref{eq:h-integral} is independent of the path. Differentiation gives 
$$ h_x=-\Rea f=-p,\qquad h_q=\Ima f=y.$$ 
Since $y_q=\Rea f'>0$, we may locally write $q=q(x,y)$. Using $u(x,y)=yq(x,y)-h(x,q(x,y))$ we have
$$ u_x=yq_x-h_x-h_q q_x=p,\qquad u_y=q+yq_y-h_q q_y=q.$$ T
hus $\nabla u=T$ on $\Omega$.\\

Conversely, suppose $\nabla u=T$. Define $h(x,q)=q\Ima f(x+iq)-u(X(x+iq))$. Therefore $h_x=-\Rea f$ and $h_q=\Ima f$. For any closed piecewise smooth curve $\gamma(t)=x(t)+iq(t)$ in $D$, 
$$ \frac d{dt}h(\gamma(t)) =-\Rea\bigl(f(\gamma(t))\gamma'(t)\bigr). $$ 
This shows \eqref{eq:periods}. Since $\Omega$ is connected, $u$ is defined up to a constant.
\end{proof}

\subsection{Global Supporting Plane Condition}
For a nonconvex source, $D^2u>0$ doesn't show convexity of $u$. The required global comparison can be expressed directly in terms of $f$.
\begin{prop}\label{prop:support}
Assume the hypotheses of Proposition~\ref{prop:periods} and suppose that $\Omega$ and $\Lambda$ are bounded. The potential $u$ in \eqref{eq:h-integral} is the restriction of a convex function on $\R^2$ if and only if
\begin{equation}\label{eq:support}
 \Rea\int_z^w f(\zeta)\,d\zeta \geq \Rea f(z)\,\Rea(w-z)-\Ima f(w)\,\Ima(w-z) 
\end{equation}
where $(z,w)\in D.$ The real integral is path independent. Equivalently,
\begin{equation}\label{eq:support-physical}
 u(X(w))\geq u(X(z))+Y(z)\cdot(X(w)-X(z)).
\end{equation}
\end{prop}

\begin{proof}
Using \eqref{eq:h-integral} gives
\begin{align*}
 u(X(w))-u(X(z))-Y(z)\cdot(X(w)-X(z)) &=\Rea\int_z^w f(\zeta)\,d\zeta -\Rea f(z)\,\Rea(w-z)\\
 &\hspace{12mm}+\Ima f(w)\,\Ima(w-z).
\end{align*}
Thus \eqref{eq:support} and \eqref{eq:support-physical} are equivalent. Every convex extension satisfies the latter inequality since it agrees with the differentiable function $u$ on the open set $\Omega$.

Conversely, assume \eqref{eq:support-physical} and define
\begin{equation}\label{eq:convex-extension}
 \overline u(X)=\sup_{z\in D} \{u(X(z))+Y(z)\cdot(X-X(z))\},\qquad X\in\R^2.
\end{equation}
This is convex and agrees with $u$ on $\Omega$. It is finite: fixing $z_0\in D$ and writing $M=\sup_{z\in D}|Y(z)|<\infty$, every affine function in the supremum is at most $u(X(z_0))+M|X-X(z_0)|$. The affine function associated with $z_0$ gives a finite lower bound.
\end{proof}

Theorem \ref{thm:global-generator-intro} follows immediately from the following theorem.
\begin{thm}\label{thm:global-generator}
Let $D\subset\mathbb C$ be connected and open, and let
$f\in\operatorname{Hol}(D)$. Suppose the coordinate assignments in \eqref{eq:assignments} have bounded images $\Omega$ and $\Lambda$, and assume:
\begin{itemize}
\item \textbf{Global matching.} Both coordinate assignments from $D$ are injective.
\item \textbf{Local positivity.} $\Rea f'>0$ in $D$.
\item \textbf{Existence of a global potential.} The real-period
condition \eqref{eq:periods} holds.
\item \textbf{Global convexity.} The supporting-plane condition
\eqref{eq:support} holds.
\end{itemize}
Then
$$
T(X(z))=Y(z)
$$
defines a smooth diffeomorphism $T:\Omega\to\Lambda$, which is the unique optimal transport between the uniform densities on $\Omega$ and $\Lambda$, up to null sets. Its potential satisfies
$$
T=\nabla u,\qquad D^2u>0,\qquad
\det D^2u=1,\qquad |\Omega|=|\Lambda|.
$$
\end{thm}

\begin{proof}
Local invertibility and the injectivity assumptions give global coordinate diffeomorphisms, and hence a smooth diffeomorphism $T:\Omega\to\Lambda$. Proposition~\ref{prop:periods} supplies $u$.
Let $\overline u$ be its finite convex extension supplied by
Proposition~\ref{prop:support}. Equation~\eqref{eq:Hessian} gives $D^2u>0$ and $\det DT=1$. Change of variables therefore gives
$$
T_\#(1_\Omega\,dX)=1_\Lambda\,dY,\qquad |\Omega|=|\Lambda|.
$$
For convex $\overline u$, $T=\nabla\overline u$ on $\Omega$,  and $T$ pushes $1_\Omega\,dX$ forward to $1_\Lambda\,dY$, Brenier's theorem implies that $T$ is the unique optimal transport map.
\end{proof}

\begin{rem}
In the above construction both $T$ and its inverse are smooth in the interiors of their respective domains. Under the hypotheses of Theorem~\ref{thm:global-generator}, injectivity of $T$ imply that $f$ is injective. Moreover, Its holomorphic inverse $f^{-1}:f(D)\to D$ generates the inverse optimal transport map $T^{-1}:\Lambda\to\Omega$. Extending this diffeomorphism to the boundaries is a separate question, that we study later.
\end{rem}

\begin{cor}\label{cor:convex-source}
Suppose $\Rea f'>0$ in $D$ and the map from $D$ to $\Omega$
is injective for convex $\Omega$. Also assume that $\Omega,\Lambda$ are bounded. Then the map defined by $T(X(z))=Y(z)$ is the optimal transport map
from $\Omega$ to $\Lambda$.
\end{cor}

\begin{proof}
The positivity and source injectivity assumptions make $X$ a smooth diffeomorphism onto $\Omega$, so $T$ is well defined and smooth. Since $\Omega$ is convex, the symmetry of $DT$ in \eqref{eq:Hessian} gives a single-valued potential $u$ with $\nabla u=T$. Proposition~\ref{prop:periods} therefore gives the
real-period condition.

Since $D^2u>0$, integration along line segments in $\Omega$ gives
$$
u(X')\geq u(X)+T(X)\cdot(X'-X),
\qquad X,X'\in\Omega.
$$
By Proposition~\ref{prop:support}, the supporting-plane condition \eqref{eq:support} holds. The same integration gives
$$
(T(X)-T(X'))\cdot(X-X')>0
\qquad\text{for }X\neq X',
$$
so $T$ is injective. Injectivity of $Y$ follows from $Y(z)= T(X(z))$. Thus we can apply Theorem \ref{thm:global-generator}.
\end{proof}

\section{Examples}\label{sec:examples}
Explicit optimal transport maps are generally difficult to construct. As a first application of the preceding construction, we give several concrete examples.\\
First, we focus on some easy examples in certain known cases.
\subsection{Linear Generators} 
For $\lambda>0$, take
$f(z)=\lambda z,$ defined on
 $$D=\{x+iq:x^2+\lambda^2q^2<1\}.
$$
The assignments in \eqref{eq:assignments} are
$$
 X(z)=(x,\lambda q),\quad Y(z)=(\lambda x,q).
$$
$X$ is injective with image $B_1$, and the target image is the ellipse
$$
 \Lambda=\left\{(p,q):
       \frac{p^2}{\lambda^2}+\lambda^2q^2<1\right\}.
$$
Since $\Rea f'=\lambda>0$ and the source is convex,
Corollary~\ref{cor:convex-source} shows that
$$
 T(x,y)=(\lambda x,y/\lambda)
$$
is the unique optimal transport from $B_1$ to $\Lambda$. Using Proposition~\ref{prop:periods}, we obtain
$$
 h(x,q)=\frac{\lambda}{2}(q^2-x^2),\;\;\; u(x,y)=\frac12\left(\lambda x^2+\frac{y^2}{\lambda}\right),
$$
up to additive constants.\\

The same generator also gives the optimal transport from a square to a rectangle of equal area. Taking $D=\{x+iq:|x|<1,\ |q|<1/\lambda\}$ gives $\Omega=(-1,1)^2$ and $\Lambda=(-\lambda,\lambda)\times(-1/\lambda,1/\lambda)$, with $T(x,y)=(\lambda x,y/\lambda)$.

\subsection{The square-root generator}
This example recovers the radial annulus to disk transport studied in \cite{MR2}.\\
Fix $0<a<R$ and set
 $f(z)=\sqrt{z^2-a^2}$ where
$$ D=\left\{x+iq:\frac{x^2}{R^2}
       +\frac{q^2}{R^2-a^2}<1\right\}\setminus[-a,a].$$
Here the square root is the holomorphic branch on
$\C\setminus[-a,a]$ satisfying $f(z)/z\to1$ at infinity. For $a<r<R$,
$$
 \begin{aligned}
 z&=r\cos\theta+i\sqrt{r^2-a^2}\sin\theta,\\
 f(z)&=\sqrt{r^2-a^2}\cos\theta+ir\sin\theta.
 \end{aligned}
$$
These mixed ellipses fill $D$, and \eqref{eq:assignments} gives
$$
 X(z)=r(\cos\theta,\sin\theta),\quad
 Y(z)=\sqrt{r^2-a^2}(\cos\theta,\sin\theta).
$$
Both maps are injective, and $\Rea f'=\Rea(z/f(z))>0$.
Thus the resulting map is
\begin{equation}\label{eq:radial-map}
 T(X)=\frac{\sqrt{|X|^2-a^2}}{|X|}\,X,
 \qquad
 B_R\setminus\overline{B_a}
 \longrightarrow B_{\sqrt{R^2-a^2}}\setminus\{0\}.
\end{equation}
This map is the gradient of the potential
\begin{equation}\label{eq:radial-potential}
 u(X)=\int_a^{|X|}\sqrt{t^2-a^2}\,dt,\qquad |X|>a.
\end{equation}
Extend $u$ by zero on $\overline{B_a}$. We have $u\in C^1(\R^2)$.
By definition it is globally convex as well. Propositions~\ref{prop:periods} and~\ref{prop:support}
therefore give the real-period and supporting-plane conditions.
Theorem~\ref{thm:global-generator} proves optimality. In particular, $a=1$ and $R=\sqrt2$
give
$$
 B_{\sqrt2}\setminus\overline{B_1}\longrightarrow B_1.
$$

One can construct the preceding examples directly using linearity and symmetry. The next example illustrates the strength of our approach beyond these settings. It gives an explicit formula for the optimal transport from the region between two non-homothetic ellipses to the interior of an ellipse.
\subsection{A Laurent generator: an elliptical annulus to a disk}\label{subsec:ellipse}
Let $\lambda,\mu>0$, and set $r_0=\sqrt{\mu/\lambda}$. Choose $R>r_0$.
Take
$f(z)=\lambda z-\frac{\mu}{z}$ where $D=\{r_0<|z|<R\}$.
For $z=re^{i\theta}$, the maps in \eqref{eq:assignments} are
\begin{equation}\label{eq:ellipse-points}
 \begin{aligned}
 X(z)&=\left(r\cos\theta,
        \left(\lambda r+\frac{\mu}{r}\right)\sin\theta\right),\\
 Y(z)&=\left(\left(\lambda r-\frac{\mu}{r}\right)\cos\theta,
        r\sin\theta\right).
 \end{aligned}
\end{equation}
Write
$$
 E(a,b)=\left\{(x,y): \frac{x^2}{a^2}+\frac{y^2}{b^2}<1\right\}, \qquad K=E(r_0,2\lambda r_0).
$$
For each fixed $r\in(r_0,R)$, the circle $|z|=r$ is mapped by $X$ onto the ellipse $\partial E(r,\lambda r+\mu/r)$ and by $Y$ onto the ellipse $\partial E(\lambda r-\mu/r,r)$.\\
Since both semiaxes in each family are strictly increasing for $r>r_0$, therefore ellipses corresponding to different radii are disjoint. Each circle is mapped injectively onto its image ellipse, hence both maps are injective. As $r\downarrow r_0$, the source ellipses converge to $\partial K$, while the target ellipses converge to the segment $\{0\}\times[-r_0,r_0]$. The resulting bounded images are
\begin{equation}\label{eq:ellipse-domains}
\begin{aligned}
\Omega&=E(R,\lambda R+\mu/R)\setminus\overline K,\\
\Lambda&=E(\lambda R-\mu/R,R)\setminus\bigl(\{0\}\times[-r_0,r_0]\bigr).
\end{aligned}
\end{equation}
Moreover, $\Rea f'(z)\geq\lambda-\mu/|z|^2>0$ on $D$. The residue $-\mu$ is real, so the real-period condition holds.
Proposition~\ref{prop:periods} gives, up to an additive constant,
\begin{equation}\label{eq:ellipse-potential}
 u(X(re^{i\theta}))
 =\frac{\lambda}{2}(r^2-r_0^2)
   -\mu\log(r/r_0)+\mu\sin^2\theta,
 \qquad \nabla u(X(z))=Y(z).
\end{equation}
\begin{figure}[ht]
    \centering
    \includegraphics[width=0.7\textwidth]{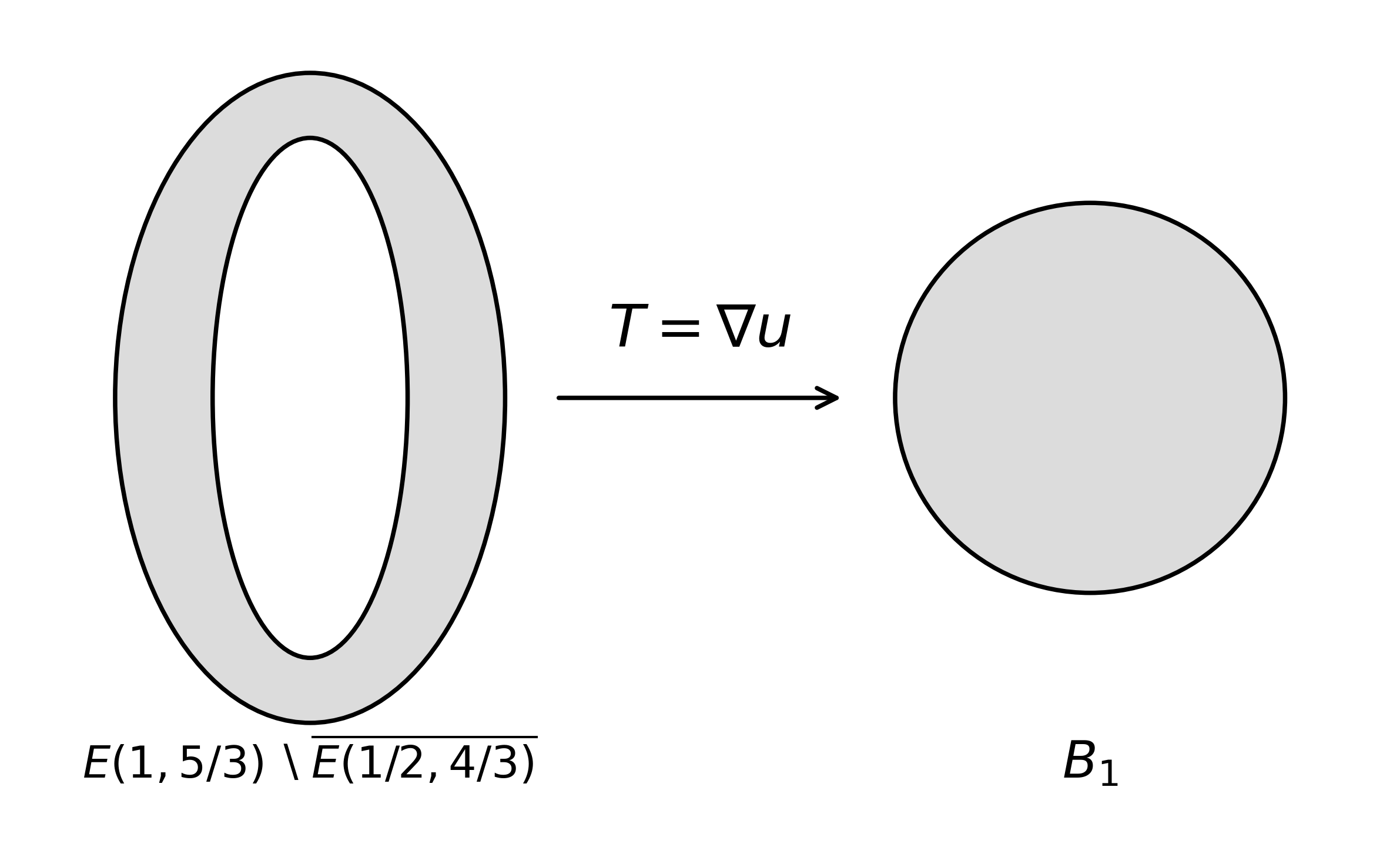}
    \caption{Optimal transport from an elliptical annulus to the unit disk.}
    \label{fig:ellipse-transport}
\end{figure}

Allowing $r\geq R$ in the same formula extends $u$ across and beyond the outer source ellipse, so $u$ is defined on $\R^2\setminus\overline K$. Define $\overline u=u$ there and $\overline u=y^2/(4\lambda)$ on $\overline K$. It is easy to check $\overline u$ is convex.\\
Thus, Theorem~\ref{thm:global-generator} shows that
$$
 T(X(re^{i\theta}))=Y(re^{i\theta}),\qquad r_0<r<R,
$$
is the unique quadratic-cost optimal transport between the uniform
densities on the domains in \eqref{eq:ellipse-domains}, upto measure zero set.

In particular, $\lambda=4/3$, $\mu=1/3$, and $R=1$ give
$$
 E(1,5/3)\setminus\overline{E(1/2,4/3)}
 \longrightarrow B_1.
$$
See Figure~\ref{fig:ellipse-transport}.

\subsection{A triangle transported to a nonconvex parabolic strip}
\label{subsec:triangletoparabolic}

Take
$$
 f(z)=z^2,\qquad
 D=\left\{x+iq:0<x<1,\ |q|<\frac12\right\}.
$$
The assignments in \eqref{eq:assignments} are
$$
 X(z)=(x,2xq),\qquad Y(z)=(x^2-q^2,q).
$$
The source assignment is injective, since $q=y/(2x)$, and its image is the triangle
$$
 \Omega=\{(x,y):0<x<1,\ |y|<x\},
$$
with vertices $(0,0)$ and $(1,\pm1)$. The target image is
$$
 \Lambda=\left\{(p,q):|q|<\frac12,\ -q^2<p<1-q^2\right\}.
$$
This is a nonconvex parabolic strip. For example, $(-1/32,\pm1/4)$ belong to $\Lambda$, whereas their midpoint $(-1/32,0)$ does not.

\begin{figure}[ht]
    \centering
    \includegraphics[width=0.9\textwidth]{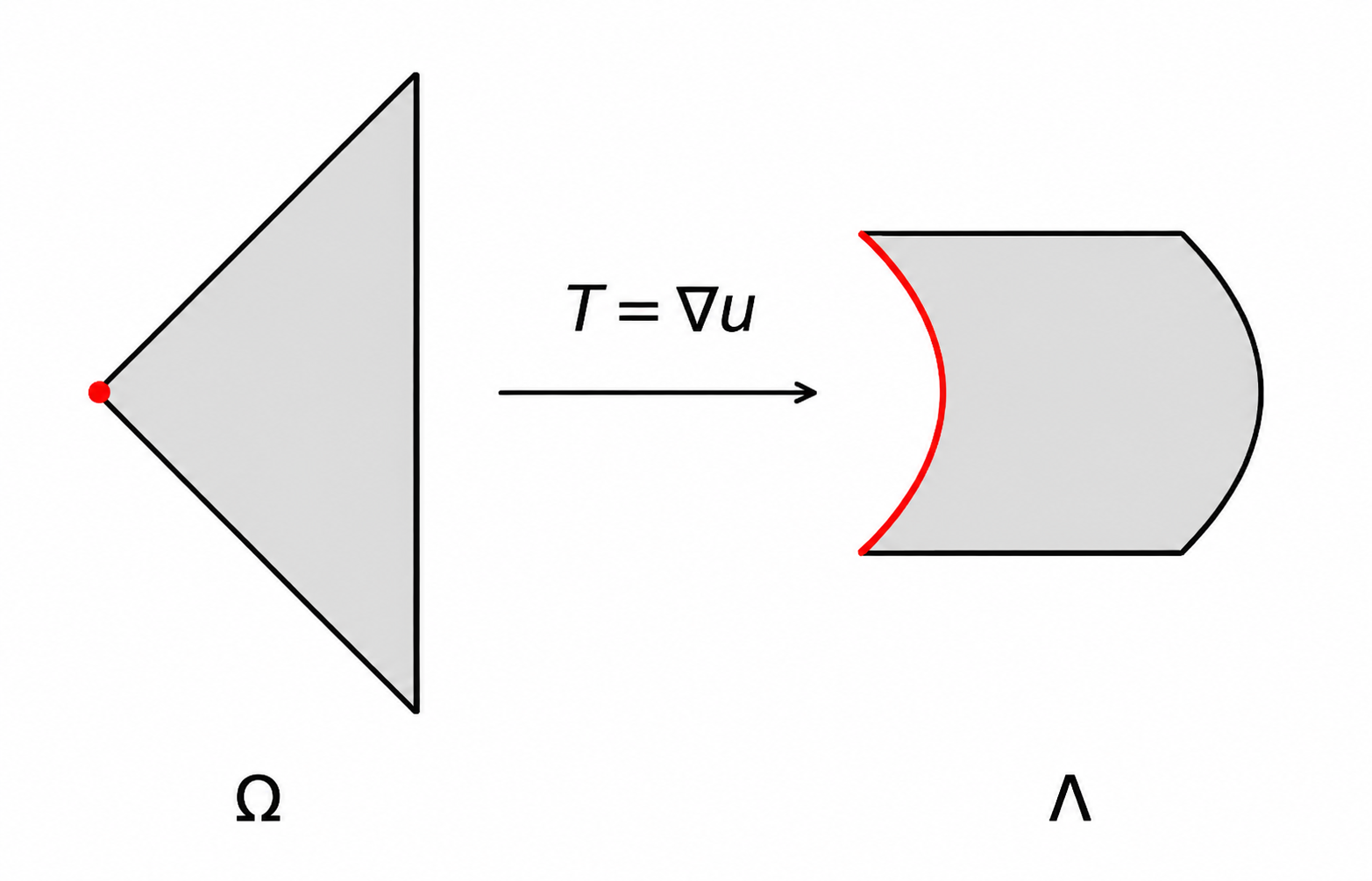}
    \caption{Optimal transport from a triangle to a nonconvex parabolic strip.}
    \label{fig:triangle-to-parabola}
\end{figure}
Since $\Rea f'(z)=2x>0$ in $D$, the source is convex, and the target is bounded. Corollary~\ref{cor:convex-source} shows that
\begin{equation}\label{eq:triangle-parabolic-map}
 T(x,y)=\left(x^2-\frac{y^2}{4x^2},\frac{y}{2x}\right)
\end{equation}
is the unique optimal transport between $\Omega$ and $\Lambda$. Both domains have unit area. The corresponding convex potential is
$$
 u(x,y)=\frac{x^3}{3}+\frac{y^2}{4x},\;\;\; \nabla u=T,\qquad \det D^2u=1.
$$
The formula shows that $T$ is smooth in $\Omega$, but its limit at the vertex $(0,0)$ depends on the direction of approach. For each fixed $|q|<1/2$,
$$
 T(x,2qx)=(x^2-q^2,q)\longrightarrow(-q^2,q)
 \qquad\text{as }x\downarrow0.
$$
The endpoints of this arc are also attained as limits by letting $q\to\pm1/2$ while $x\downarrow0$. Conversely, along any sequence $(x,y)\in\Omega$ approaching $(0,0)$, the quantity $q=y/(2x)$ lies in $(-1/2,1/2)$, and the formula $T(x,y)=(x^2-q^2,q)$ shows that every limiting value belongs to
$$
 \Gamma=\left\{(-q^2,q):|q|\leq\frac12\right\}.
$$
Thus $\Gamma$ is exactly the set of limiting values of $T$ at the vertex, and $T$ has no continuous extension there.

This example exhibits an interesting boundary discontinuity of optimal transport from a convex source to a nonconvex target: the map is smooth throughout the interior, but its limiting gradients at a single boundary point fill an entire parabolic arc of the target boundary.

\subsection{A quadratic generator: a disk to a smooth nonconvex domain}
\label{subsec:quadratic-generator}
Take
$$
 f(z)=\frac12\left(z+\frac32\right)^2,\qquad
 D=\left\{x+iq:x^2+\left(x+\frac32\right)^2q^2<1\right\}.
$$
The maps $X,Y$ in \eqref{eq:assignments} are
$$
 X(z)=\left(x,\left(x+\frac32\right)q\right),\qquad
 Y(z)=\left(\frac{(x+\frac32)^2-q^2}{2},q\right).
$$

\begin{figure}[ht]
\centering
\includegraphics[width=0.8\textwidth]{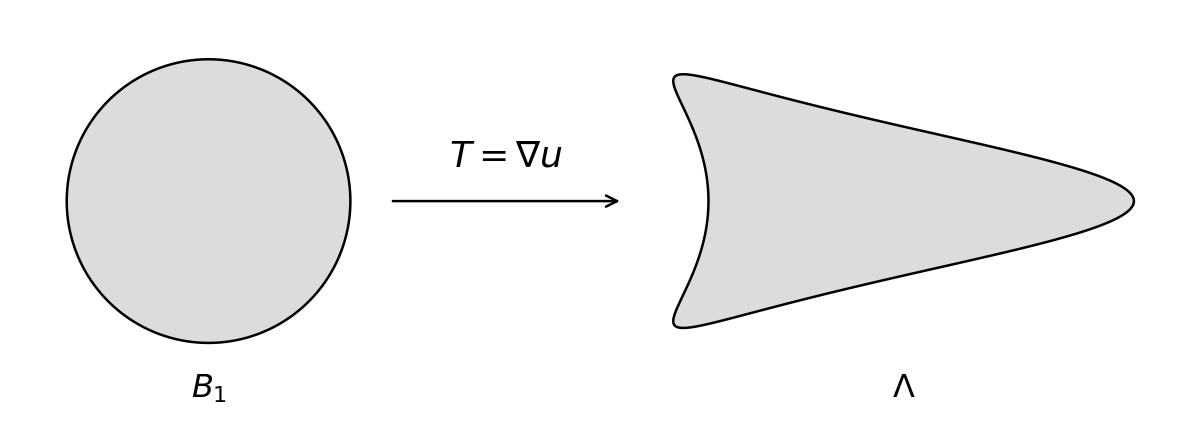}
\caption{Smooth optimal transport to a nonconvex domain.}
\label{fig:smooth-nonconvex}
\end{figure}

\begin{figure}[ht]
\centering
\includegraphics[width=0.5\textwidth]{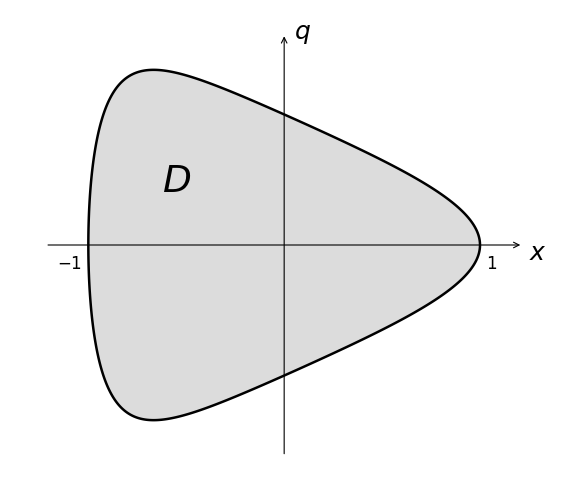}
\caption{The mixed domain for the quadratic generator.}
\end{figure}

The source assignment is a diffeomorphism onto $B_1$. Since $\Rea f'=x+\frac32\geq1/2$ on $\overline D$ and the target image is bounded, Corollary~\ref{cor:convex-source} shows that
\begin{equation}\label{eq:quadratic-map}
 T(x,y)=
 \left(
 \frac{(x+\frac32)^2}{2}
 -\frac{y^2}{2(x+\frac32)^2},
 \frac{y}{x+\frac32}
 \right),
 \qquad \Lambda=T(B_1),
\end{equation}
is the unique quadratic-cost optimal transport between the uniform densities on $B_1$ and $\Lambda$. The potential is
\begin{equation}\label{eq:quadratic-potential}
 u(x,y)=\frac{(x+\frac32)^3}{6}
        +\frac{y^2}{2(x+\frac32)}.
\end{equation}
$T$ is a smooth nondegenerate diffeomorphism of the closures, $\partial\Lambda$ is smooth, and $|\Lambda|=\pi$. This example was studied by McCann and the first author \cite{McCannRakshit}. The target $\Lambda$ is nonconvex, yet the optimal transport map is smooth up to the boundary.

\begin{rem}
The particular choice of mixed domain above is not essential. The same generator on the ellipse
$$
 D=\left\{x+iq:x^2+\frac49q^2<1\right\}
$$
also gives an optimal transport map to a nonconvex target, smooth up to the boundary.
\end{rem}

\subsection{A family of targets using higher order polynomials}\label{ex:higher_order_poly} Consider the polynomial family,
\begin{figure}[ht]
    \centering
    \includegraphics[width=0.8\linewidth]{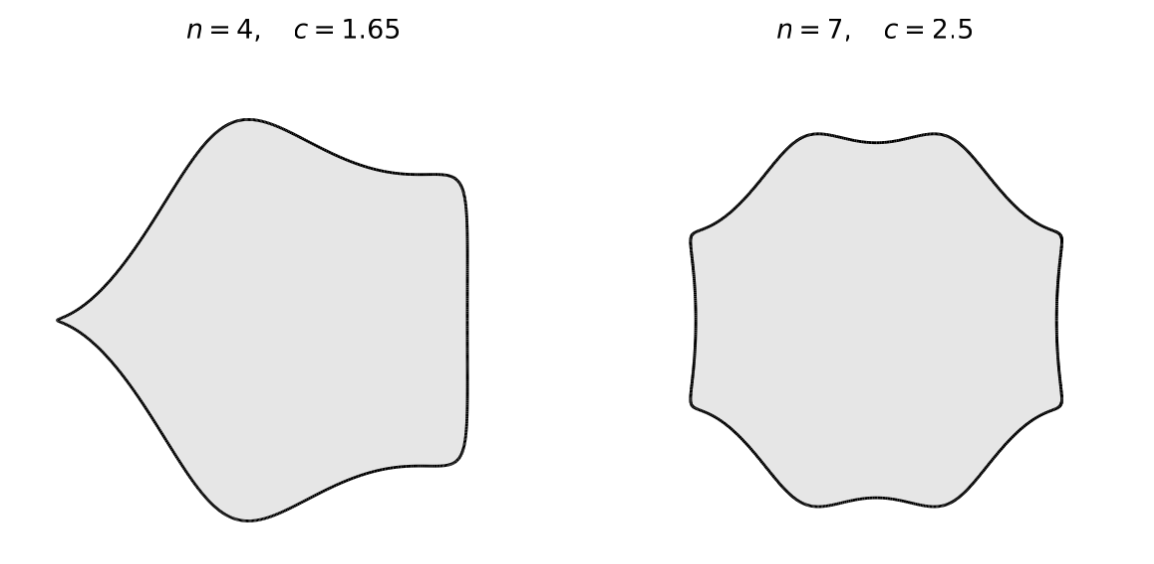}
    \caption{Different nonconvex target domains that have smooth optimal transport maps generated by the polynomials \eqref{eq:poly_family_ex}.}
\end{figure}

For this family one can force the source of these generators to be unit circles by Implicit function theorem. For $(x,y)\in\overline B_1(0)$, solve
\begin{equation}\label{implicit}
\Ima f_{n,c,\vartheta}(x+iq)=y,
\quad \implies\quad
q=y-\frac{c}{n^2}\Ima\big(e^{i\vartheta}(x+iq)^n\big).
\end{equation}
\begin{equation}\label{eq:poly_family_ex}
    f_{n,c,\vartheta} = z+\frac{c}{n^2}e^{i\vartheta}z^n.
\end{equation}
On the interval $|q-y|\le 2c/n^2$, one has
$|x+iq|\le 1+2c/n^2$, so, for large $n$,
$$
\frac{c}{n^2}|x+iq|^n\le \frac{2c}{n^2},\qquad
\left|\frac{c}{n}e^{i\vartheta}(x+iq)^{n-1}\right|\le\frac{2c}{n}<\frac12.
$$

Thus the right side of \eqref{implicit} is a contraction of this interval into itself. Its unique fixed point $q(x,y)$ is smooth by the implicit-function theorem and satisfies

\begin{equation}\label{bounds}
q=y+O_c(n^{-2}),\qquad
f_{n,c,\vartheta}'(x+iq(x,y))=1+O_c(n^{-1}).
\end{equation}
The same strict inequalities yield smooth extension across $\partial B_1(0)$. Note that since $\Rea f_{n,c,\vartheta}'>0$ for $n$ large we get that the map induced by $f_{n,c,\vartheta}$ is indeed an optimal transport map by Corollary \ref{cor:convex-source}.

\section{Characterization}\label{sec:Characterization}
Our next goal is to characterize smooth optimal transport between the unit disk and a possibly nonconvex domain, both equipped with uniform densities. First taking the unit ball as the source, we give a necessary and sufficient analytic condition on the boundary correspondence for the transport to be smooth and nondegenerate up to the boundary. This condition is expressed through the vanishing of certain Cauchy moments.
\subsection{Mixed Boundary and Cauchy Moments}
Let $\Lambda$ be a bounded domain with smooth Jordan boundary. Define
$$ 
e(t)=(\cos t,\sin t),\qquad \gamma(t)=(p(t),q(t)),\qquad t\in[0,2\pi), 
$$ 
where $\gamma$ is a positively oriented regular parametrization of $\partial\Lambda$. All parameterizations are understood to be $2\pi$-periodic. We give conditions under which the prescribed boundary correspondence \(\gamma(t)\) is induced by a smooth optimal transport map from $B_1$ to $\Lambda$. Using our convention $z=x+iq$ and $f(z)=p+iy$, we define
\begin{equation}\label{eq:boundary-data}
\zeta(t)=\cos t+iq(t),\qquad \eta(t)=p(t)+i\sin t.
\end{equation}
Here $\zeta$ parametrizes the curve in the mixed $(x,q)$-plane. Assume that it is a positively oriented regular Jordan curve, and let $D$ be its bounded interior. We seek a holomorphic function $f$ in $D$ with boundary value $f(\zeta(t))=\eta(t)$. To study whether these boundary values admit a holomorphic extension, define the Cauchy moments
\begin{equation}\label{eq:moments}
M_j[\gamma]=\int_0^{2\pi}\eta(t)\zeta(t)^j\zeta'(t)\,dt,\qquad j=0,1,2,\ldots.
\end{equation}
\begin{lem}\label{lem:cauchy}
The function $\eta$ extends to an $f\in\Hol(D)\cap C^\infty(\overline D)$ with $f(\zeta(t))=\eta(t)$ if and only if all the moments \eqref{eq:moments} vanish. The extension is unique and is given by
\begin{equation}\label{eq:Cauchy-reconstruct}
 f(w)=\frac1{2\pi i}\int_0^{2\pi} \frac{\eta(t)\zeta'(t)}{\zeta(t)-w}\,dt, \qquad w\in D.
\end{equation}
\end{lem}
\begin{proof}
Necessity follows from Cauchy's theorem. Conversely, suppose all the moments vanish. Define $f$ by \eqref{eq:Cauchy-reconstruct}, away from $\partial D$. This defines a holomorphic function on each component of $\C\setminus\partial D$. For $|w|>\max_t|\zeta(t)|$, expanding the denominator as a convergent geometric series gives
$$
 f(w)=-\frac1{2\pi i}
       \sum_{j=0}^{\infty}\frac{M_j[\gamma]}{w^{j+1}}=0.
$$
Since $\C\setminus\overline D$ is connected, the identity theorem gives $f=0$ throughout the exterior. Let $f_+$ and $f_-$ denote the exterior and interior boundary values,
respectively as in Sokhotski--Plemelj formulas. Since the exterior function vanishes identically, we have $f_+=0$. Therefore
$
 f_-(\zeta(t_0))=\eta(t_0).
$
Thus the interior function extends continuously to $\overline D$ with the prescribed boundary values.\\
To obtain smoothness up to the boundary, differentiate the Cauchy
integral and integrate by parts to get,
$$
 f'(w)=\frac1{2\pi i}\int_0^{2\pi}
       \frac{\eta'(t)}{\zeta(t)-w}\,dt.
$$
Since $\zeta'$ never vanishes, this is again a Cauchy integral with smooth boundary density $\eta'/\zeta'$. Repeating the argument and applying the same boundary formulas shows that every derivative of $f$ extends continuously to $\overline D$. Hence $f\in\Hol(D)\cap C^\infty(\overline D)$.
\end{proof}

Since \(D\) is simply connected, the holomorphic extension \(f\) automatically satisfies the real-period condition.

\subsection{Characterization} In this section we characterize all target domains that admit a smooth optimal transport map from the unit ball.
\begin{lem}\label{lem:Hopf}
Suppose $\zeta(t)=\cos t+iq(t)$ is a positively oriented smooth regular Jordan curve and $f\in\Hol(D)\cap C^\infty(\overline D)$ satisfies $\Ima f(\zeta(t))=\sin t$. Then $$ \Rea f'\geq c>0\quad\text{on }\overline D. $$
\end{lem}
\begin{proof}
Fix $t_0$ and set $$ W(z)=1-\cos t_0\Rea z-\sin t_0\Ima f(z). $$ On the boundary, $W(\zeta(t))=1-\cos(t-t_0)\geq0$, with a unique zero at $t_0$. Thus, by maximum principle $W>0$ in $D$, and the outward normal derivative at $\zeta(t_0)$ is strictly negative by the Hopf lemma. 
Since $\zeta'(t)=-\sin t+iq'(t)$, differentiating
$\Ima f(\zeta(t))=\sin t$ gives
$$
\begin{aligned}
\cos t_0
&=\Ima\bigl(f'(\zeta(t_0))\zeta'(t_0)\bigr)\\
&=-\Ima f'(\zeta(t_0))\sin t_0
+\Rea f'(\zeta(t_0))q'(t_0).
\end{aligned}
$$
Consequently $$ \nabla W(\zeta(t_0)) =-\Rea f'(\zeta(t_0))(q'(t_0),\sin t_0). $$ 
Thus the outward unit normal is $(q'(t_0),\sin t_0)/|\zeta'(t_0)|$. It follows that
$$
\partial_\nu W(\zeta(t_0)) =\nabla W\cdot \nu=-\Rea f'(\zeta(t_0))|\zeta'(t_0)|<0.
$$
Hence, \(\Rea f'(\zeta(t_0))>0\). Now, this holds for any boundary point.  Moreover, \(\operatorname{Re}f'\) is continuous and \(\partial D\) is compact, so its minimum is attained and is strictly positive. Finally, \(\operatorname{Re}f'\) is harmonic in \(D\). The minimum principle therefore completes the proof.
\end{proof}

We now characterize all smooth, possibly nonconvex, targets that admit a smooth optimal transport from the disk up to the boundary. The criterion is expressed entirely in terms of a boundary parametrization.

\begin{thm}\label{thm:disk}
Let $\Lambda\subset\R^2$ be a bounded domain with smooth Jordan
boundary. The following statements are equivalent.
\begin{enumerate}
\item There is a smooth optimal transport map $T:B_1\to\Lambda$ extending to a smooth diffeomorphism of the closures.
\item The boundary $\partial\Lambda$ admits a positively oriented smooth regular parametrization $\gamma$ such that the associated curve $\zeta$
in \eqref{eq:boundary-data} is a positively oriented smooth regular Jordan curve and
\begin{equation}\label{eq:main-moments}
 \int_0^{2\pi}\eta(t)\zeta(t)^j\zeta'(t)\,dt=0,
 \qquad j=0,1,2,\ldots.
\end{equation}
\end{enumerate}
\end{thm}
\begin{proof}
Assume that $T=\nabla u$ is the transport in the first statement. We have $D^2u>0$ up to the boundary. Then the mixed-coordinate map $(x,y)\mapsto(x,u_y(x,y))$ has positive Jacobian. It is globally injective on $\overline B_1$, since each vertical section is an interval and $u_y$ is strictly increasing there. Its image is a smooth Jordan domain, with boundary $\zeta$.  The partial Legendre transform gives a smooth holomorphic generator with boundary trace $\eta$. Lemma~\ref{lem:cauchy} gives us equation \eqref{eq:main-moments}. 

Conversely, suppose we have equation \eqref{eq:main-moments}. To construct the transport, we first apply Lemma~\ref{lem:cauchy} to obtain a holomorphic function $f$, smooth up to the boundary, with
$f(\zeta(t))=\eta(t)$. Lemma~\ref{lem:Hopf} gives
$\Rea f'\geq c>0$ on $\overline D$.\\
The source assignment $z\mapsto(\Rea z,\Ima f(z))$ sends $\zeta(t)$ to $e(t)$. Moreover
$$
 H(z)=(\Rea z)^2+(\Ima f(z))^2=1
$$
 on $\partial D$ and
$$
 \Delta H=2+2|f'(z)|^2>0.
$$
The maximum principle therefore gives $H<1$ in $D$. So the interior image lies in $B_1$.\\
Since $\Rea\zeta(t)=\cos t$,
for each $x\in(-1,1)$, the vertical line $\Rea z=x$ meets $\partial D$ at exactly two points. Because $\zeta$ is a Jordan curve, the vertical section of $D$ is the interval between these points. On this interval,
$$
 \partial_q\Ima f(x+iq)=\Rea f'(x+iq)>0.
$$
The boundary identity $\Ima f(\zeta(t))=\sin t$ gives the two endpoint values $-\sqrt{1-x^2}$ and $\sqrt{1-x^2}$. Thus $q\mapsto\Ima f(x+iq)$ maps this interval bijectively onto $(-\sqrt{1-x^2},$ $\sqrt{1-x^2})$. Therefore the source assignment is a bijection from $D$ onto $B_1$.

Its boundary restriction is bijective. Its Jacobian is $\Rea f'\geq c>0$ on $\overline D$. Thus, the inverse function theorem gives a smooth inverse up to the boundary. Hence $X$ is a diffeomorphism of the closures.

The reconstructed map has the symmetric positive-definite derivative \eqref{eq:Hessian}, with $\det DT=1$. Its image is bounded. Corollary~\ref{cor:convex-source} gives $T=\nabla u$ for a strictly convex potential. Thus proves optimality. Integration along line segments gives injectivity on $\overline B_1$. Since $T(e(t))=\gamma(t)$, the Jordan curve theorem gives $T(B_1) =\Lambda$. The derivative remains invertible on $\overline B_1$, hence $T$ is a smooth diffeomorphism of the closures.
\end{proof}
\begin{rem}\label{rem:charac_low_reg}
    In the above theorem smoothness assumption on $\partial\Lambda$ and the parametrization $\gamma$ can be weakened to $C^{1,\alpha}$ to ensure an optimal transport map $T\in C^{1,\alpha}(\bar B_1(0))$.
\end{rem}

\subsection{Generalizing the source to arbitrary domains}
Let $\Omega$ be a bounded smooth convex domain. Fix a positively oriented smooth parametrization
$$ \sigma(t)=(x(t),y(t)), $$ 
of $\partial\Omega$. For a positively oriented parametrization $\gamma(t)=(p(t),q(t))$ of a target boundary, use
\begin{equation}\label{eq:general-data}
 \zeta(t)=x(t)+iq(t),\qquad \eta(t)=p(t)+iy(t).
\end{equation}
We prove the following version of Theorem \ref{thm:intro_thm1}.

\begin{thm}\label{thm:convex-boundary}
Let $\Omega$ be as above and let $\Lambda$ be a bounded smooth Jordan domain. There exists a smooth optimal transport map between $\bar\Omega$ and $\bar\Lambda$ if and only if $\partial\Lambda$ admits a positively oriented smooth parametrization for which $\zeta$ in \eqref{eq:general-data} is a positively oriented smooth Jordan curve and
\begin{equation}\label{eq:general-moments}
 \int_0^{2\pi}\eta(t)\zeta(t)^j\zeta'(t)\,dt=0 \quad(j=0,1,2,\ldots).
\end{equation}
\end{thm}
\begin{proof}
Suppose first that there is a smooth optimal transport map from $\bar\Omega$ to $\bar\Lambda$. As in the proof of Theorem~\ref{thm:disk}, vertical sections of $\Omega$ are intervals. Thus the mixed coordinate map is globally injective. Define $\gamma(t)=T(\sigma(t))=(p(t),q(t))$. The partial
Legendre transform gives a holomorphic function $f$ on $D$, smooth up to the boundary, such that $f(x+iq)=p+iy$ whenever $(p,q)=T(x,y)$. At the mixed boundary point $\zeta(t)=x(t)+iq(t)$, therefore
$$
 f(\zeta(t))=p(t)+iy(t)=\eta(t).
$$
Thus $\eta$ is the boundary value of $f$. Cauchy's theorem then gives
$$
 \int_0^{2\pi}\eta(t)\zeta(t)^j\zeta'(t)\,dt
 =\oint_{\partial D}f(z)z^j\,dz=0,
 \qquad j=0,1,2,\ldots.
$$

For the other direction, Lemma~\ref{lem:cauchy} gives the holomorphic extension. Let $\nu=(\nu_1,\nu_2)$ be the outward unit normal to $\partial\Omega$ at $\sigma(t_0)$ and define
$$ W(z)=\nu\cdot\bigl(\sigma(t_0)-(\Rea z,\Ima f(z))\bigr). $$ 
Convexity of the source domain gives $W\geq0$ on $\partial D$. Since $\Omega$ is open, its boundary values are not identically zero.  Thus $W>0$ in $D$ and the Hopf lemma applies at $\zeta(t_0)$.\\
From $\nu=(y'(t_0),-x'(t_0))/|\sigma'(t_0)|$ and differentiating $W$, we get
$$ \nabla W(\zeta(t_0)) =-\frac{\Rea f'(\zeta(t_0))}{|\sigma'(t_0)|} (q'(t_0),-x'(t_0)). $$ 
Therefore 
$$ \partial_{\nu_D}W(\zeta(t_0)) =-\Rea f'(\zeta(t_0)) \frac{|\zeta'(t_0)|}{|\sigma'(t_0)|}<0. $$ 
The minimum principle gives a positive lower bound for $\Rea f'$. Since $t_0$ was arbitrary, the inequality $W>0$ and convexity gives that the source image therefore lies in $\Omega$. Every vertical line meeting $\Omega$ intersects $\partial\Omega$ at exactly two points. Since $\Rea\zeta(t)=x(t)$, it also intersects $\partial D$ at exactly two points.  so the corresponding section of $D$ is an interval. On this interval, $\partial_q\Ima f(x+iq)=\Rea f'(x+iq)>0$, and the boundary values gives the two endpoint heights of the corresponding section of $\Omega$. Thus each vertical section is mapped bijectively onto the corresponding source section. Since the first coordinate is unchanged, $X$ is a bijection from $\bar D$ onto $\bar\Omega$. Therefore since its Jacobian $\Rea f'$ is positive on $\overline D$. The inverse function theorem therefore gives a smooth inverse up to the boundary. Corollary~\ref{cor:convex-source} gives us the convex potential and injectivity, and the prescribed target boundary is $\Lambda$.
\end{proof}

\begin{rem}\label{rem:autometic-obliqueness}
The boundary correspondence in the previous theorem is automatically strictly oblique. The vanishing of all the Cauchy moments gives a holomorphic generator $f$ satisfying
$f(\zeta(t))=\eta(t)$. The supporting-line and Hopf-lemma argument in the proof of Theorem~\ref{thm:convex-boundary} gives
$\Rea f'\geq c>0$ on $\overline D$.
Differentiating the boundary identity yields
$$
 \sigma'(t)\cdot\gamma'(t)
 =\Rea\bigl(\eta'(t)\overline{\zeta'(t)}\bigr)
 =|\zeta'(t)|^2\Rea f'(\zeta(t)).
$$
Thus $\Rea f'>0$ gives strict obliqueness at the boundary.
\end{rem}

\begin{rem}\label{rem:charac_non_convex}
The same boundary construction applies to nonconvex sources, provided we retain the mixed Jordan curve and moment conditions and additionally impose $\Rea f'>0$ on $\overline D$ and the supporting-plane inequality \eqref{eq:support}. The latter ensures that the reconstructed potential has a global convex extension, so Theorem~\ref{thm:global-generator} again gives optimality.
\end{rem}

\subsection{Geometric consequences}
We now study the first few moments to note a few immediate geometric consequences coming out of the moment condition.

\subsubsection{Area preservation}
The imaginary part of the zeroth moment gives the difference between the source and target areas. Using $\zeta=x+iq$ and $\eta=p+iy$, we obtain
\begin{equation}\label{eq:zeroth-general}
 \Ima M_0[\gamma]
 =\int_0^{2\pi}(pq'+yx')\,dt
 =|\Lambda|-|\Omega|.
\end{equation}
Thus $\Ima M_0[\gamma]=0$ if and only if $|\Lambda|=|\Omega|$.
\subsubsection{Four Radial Alignment} The following result exploits the periodicity of the boundary parametrization. 

\begin{thm}\label{thm:four}
Let $T=\nabla u:\overline B_1\to\overline\Lambda$ be a smooth nondegenerate optimal transport map. Assume that $\Lambda$ has barycenter at the origin and define
$$
 b(t)=\det(e(t),T(e(t))),\qquad e(t)=(\cos t,\sin t).
$$
Then
\begin{equation}\label{eq:first-harmonics}
 \int_0^{2\pi}b(t)\,dt
 =\int_0^{2\pi}b(t)\cos t\,dt
 =\int_0^{2\pi}b(t)\sin t\,dt=0.
\end{equation}
There are at least four distinct points $X\in\partial B_1$ such that $0$, $X$, and $T(X)$ lie on a common line. If $0\in\Lambda$, then $T(X)=\lambda X$ with $\lambda>0$ at every such point, so the corresponding radial directions are preserved.
\end{thm}
\begin{proof}
Since $b(t) = \cos(t) q(t)-\sin(t)p(t)$,  $\eta(t)=p(t)+i\sin(t)$ and $\zeta'(t)=-\sin(t)+iq'(t)$ using the vanishing of the $0^{th}$ moment we get,
$$
 0=\int_0^{2\pi}(-p(t)\sin(t)-\sin(t)q'(t))d t
   =\int_0^{2\pi}(-\sin(t)p(t)+\cos(t)q(t))d t.
$$
Therefore
\begin{equation}\label{eq:mean}
\displaystyle \int_0^{2\pi}b(t)d t=0.
\end{equation}

The real part of the first moment gives,
\begin{align*}
 0
 &=\Rea\int_0^{2\pi}(p(t)+i\sin(t))(\cos(t)+iq(t))(-\sin(t)+iq'(t))d t\\
 &=\int_0^{2\pi}
       \bigl(-p(t)\cos(t)\sin(t)-p(t)q(t)q'(t)-\sin(t)\cos(t)q'(t)+\sin(t)^2q(t)\bigr)d t.
\end{align*}
Using periodicity to integrate by parts,
$$
 -\int_0^{2\pi} p(t)q(t)q'(t)d t=\frac12\int_0^{2\pi} p'(t)q^2(t)d t,
$$
$$
 -\int_0^{2\pi} \sin(t)\cos(t)q'(t)d t=\int_0^{2\pi}(\cos(t)^2-\sin(t)^2)q(t)d t.
$$
Hence
\begin{equation}\label{eq:cos-pre}
 0=\int_0^{2\pi} b(t)\cos(t) d t
     +\frac12\int_0^{2\pi}q(t)^2p'(t) d t.
\end{equation}
Since the target boundary $\gamma$ is positively oriented, Green's theorem yields
$$
 0=\int_\Lambda y d A
 =-\frac12\oint_{\partial\Lambda}y^2 d p
 =-\frac12\int_0^{2\pi}q(t)^2p'(t) d t.
$$
Substitution in \eqref{eq:cos-pre} gives
\begin{equation}\label{eq:cos}
 \int_0^{2\pi} b(t)\cos(t)dt=0.
\end{equation}
The imaginary part gives the other harmonic.

If $b\equiv0$, every boundary point is aligned. Otherwise, \eqref{eq:first-harmonics} says that the constant and first Fourier harmonics of $b$ vanish. The  Sturm--Hurwitz theorem therefore implies that $b$ has at least four sign changes, and hence at least four distinct zeros, each giving a
boundary alignment.

Finally, suppose $0\in\Lambda$ and let $X_0=T^{-1}(0)\in B_1$. At an alignment write $T(X)=\lambda X$, where $|X|=1$. Since $D^2u>0$, integration along the segment from $X_0$ to $X$ gives
$$
 0<(T(X)-T(X_0))\cdot(X-X_0)
 =\lambda(1-X\cdot X_0).
$$
Since $1-X\cdot X_0\geq1-|X_0|>0$, we obtain $\lambda>0$.

\end{proof}
Note that the four radial alignments do not imply a decomposition of the transport into invariant sectors.

\section{Discontinuous Optimal Transport and moment obstructions}
\label{sec:discon}
In this section we tackle the question of detecting singularities. Let $\Lambda\subset\R^2$ be a bounded Jordan domain with $C^{2,\alpha}$
boundary for some $0<\alpha<1$, with
$$
  |\Lambda|=|B_1(0)|=\pi.
$$ 
 Choose a positively oriented arclength parametrization
$$
 \Gamma(s)=(p(s),q(s)),\qquad \Gamma(s+\ell)=\Gamma(s),
 \qquad |\Gamma'(s)|=1,
$$
where $\ell=\mathcal H^1(\partial\Lambda)$.  Fix a continuous tangent-angle
lift $\beta:\R\to\R$ so that
$$
  \Gamma'(s)=(\cos\beta(s),\sin\beta(s)),\qquad
  \beta(s+\ell)=\beta(s)+2\pi.
$$
The choice of
lift for $\beta$ fixes the normalization of the admissible lifts below.
Define the closed geometric corridor
\begin{equation}\label{eq:LR}
 L(s)=\max_{r\in[s-\ell,s]}\big(\beta(r)-\pi\big),\qquad
 R(s)=\min_{r\in[s,s+\ell]}\beta(r).
\end{equation}
Define $\Aw$ as the set of continuous nondecreasing functions
$\theta:\R\to\R$ satisfying
\begin{equation}\label{eq:Aw}
 \theta(s+\ell)=\theta(s)+2\pi,\qquad
 L(s)\le\theta(s)\le R(s).
\end{equation}
The class is permitted to be empty. For $\theta\in\Aw$, set
\begin{equation}\label{eq:mixed-traces}
 \zeta_\theta(s)=p(s)+i\sin\theta(s),\qquad
 \eta_\theta(s)=\cos\theta(s)+iq(s).
\end{equation}
Define, for each fixed integer $k\ge0$, the complex moment
\begin{equation}\label{eq:moment}
 M_k[\theta]=\int_0^\ell \eta_\theta(s)\,
 \zeta_\theta(s)^k\, d\zeta_\theta(s).
\end{equation}

\begin{thm}\label{thm:moment_obstruction}
Let $T: B_1(0)\to \Lambda$ be a continuous optimal transport map. Then there exists a correspondence $\theta\in\Aw$
such that
\begin{equation}\label{eq:allzero}
 M_k[\theta]=0\qquad\text{for every integer }k\ge0.
\end{equation}

\end{thm}

\begin{rem}\label{rem:perturb_gen_convex}
    Note that the above theorem the source domain can be replaced by any convex domain where the class $\Aw$ can be appropriately modified by accounting for the source domain with some fixed arclength parameterization $\alpha(s)$.
\end{rem}

For a $BV$ function $f$, we denote its total variation by
$$
\text{Var}_{[a,b]}(f)=\sup_{a=t_0<\cdots<t_m=b}
\sum_{j=1}^m|f(t_j)-f(t_{j-1})|.
$$
We also denote the corresponding Lebesgue–Stieltjes measure by $df$. Given a $\theta\in \Aw$, an arbitrary moment $M_k[\theta]$ it is well defined because of the following lemma.
\begin{lem}\label{lem:BV}
For every continuous nondecreasing degree-one lift $\theta$,
$$
 \var_{[0,\ell]}\big(e^{i\theta}\big)=2\pi,
 \qquad
 \var(\sin\theta)=\var(\cos\theta)=4.
$$
Consequently
\begin{equation}\label{eq:vartraces}
 \var(\zeta_\theta)\le \ell+4,
 \qquad \var(\eta_\theta)\le \ell+4,
\end{equation}
and every moment in \eqref{eq:moment} is well-defined.
\end{lem}

\subsection{A global inverse} In this section we give properties of the inverse mapping $S = T^{-1}$ and provide the choice of $\theta$.
\begin{lem}
\label{lem:external}
Let $S$ denote the optimal transport map from $\Lambda$ to $B_1(0)$.
Then:
\begin{enumerate}
\item $S=\nabla v$ is smooth on $\Lambda$ for a convex $v$ with
\begin{equation}\label{eq:inversePDE}
 DS=D^2v>0,\qquad \det DS=1.
\end{equation}
\item $S$ extends continuously to $\overline\Lambda$, with values in $\overline B_1(0)$.
\item Assume $T:B_1(0)\to \Lambda$ is a continuous optimal transport map. Then
$$
 T:B_1(0)\longrightarrow\Lambda,\qquad S:\Lambda\longrightarrow B_1(0)
$$
are inverse smooth diffeomorphisms.
\end{enumerate}
\end{lem}

Note that statement $(1)$ above follows from \cite{Caf1}, \cite{Caf2} and $(2)$ follows from \cite[Proposition 11]{Chodosh}. Given that we have no convexity information on $\Lambda$, it might still be possible that $T$ is discontinuous purely based on the first 2 points. However if $T$ is continuous then we can upgrade the a.e. inverse statement to the entire ball $B_1(0)$. Based on the above lemma we can now define the choice of the mapping $\theta$ when $T$ is continuous.

\begin{lem}\label{lem:boundary_theta}
Assume $T:B_1(0)\to \Lambda$ is a continuous optimal transport map and let $S$ be it's inverse. Then
\begin{equation}\label{eq:Sboundary}
  S(\partial\Lambda)=\partial B_1(0),
\end{equation}
and there is a continuous nondecreasing lift
$\theta:\R\to\R$ such that
\begin{equation}\label{eq:theta}
  S(\Gamma(s))=(\cos\theta(s),\sin\theta(s)),
 \qquad \theta(s+\ell)=\theta(s)+2\pi.
\end{equation}
\end{lem}
\begin{proof}
    From Lemma \ref{lem:external} we note that $S$ extends continuously to the boundary. We have to show it remains surjective under the additional assumption of $T$ being continuous. Assume for $y\in \partial\Lambda$, $S(y)\in B_1(0)$. Let $y_j\in \Lambda \to y$ then $y_j = T(S(y_j))\to T(S(y))\in \Lambda$. However $y\notin \Lambda$. Now for surjectivity if $x\in \partial B_1(0)$, consider $x_j\in B_1(0)\to x$. Note that there exists a convergent subsequence of $T(x_j)$ therefore $x_{j_k} = S(T(x_{j_k}))\to S(y)$ for some limit $T(x_{j_k})\to y$. Therefore $S(y)=x$. The existence of $\theta$ follows from the fact that an orientation-preserving homeomorphism between Jordan interiors which extends continuously to both boundaries has an order- preserving degree-one boundary trace.
\end{proof}

\subsection{Obliqueness and the closed corridor} We now show that the $\theta$ constructed in Lemma \ref{lem:boundary_theta} belongs to the class $\Aw$. This follows from a weaker version of obliqueness.

\begin{lem}\label{lem:weakflux}
For every $s\in\R$,
\begin{equation}\label{eq:weaknormal}
  S(\Gamma(s))\cdot\nu_\Lambda(\Gamma(s))\ge0,
\end{equation}
where $\nu_\Lambda$ is the outward unit normal.
Equivalently for the tangent angle $\beta$,
\begin{equation}\label{eq:sinpositive}
 \sin\big(\beta(s)-\theta(s)\big)\ge0.
\end{equation}
\end{lem}
\begin{proof}
    For $0<r<1$, let $B_r=\{x:|x|<r\}$,
$\Lambda_r=T(B_r)$, and $\Sigma_r=\partial\Lambda_r$.
By Lemma~\ref{lem:external}, these are smooth Jordan domains exhausting $\Lambda$ as $r\uparrow1$.
At $x=r e(t)$, $e(t)=(\cos t,\sin t)$, let
$\tau_B=Je(t)=(-\sin t,\cos t)$. since $DT>0$ and orientation preserving, the positively oriented tangent at $y=T(x)$ is $\tau_{\Sigma}=DT(x)\tau_B/|DT(x)\tau_B|$.
Writing $\nu_{\Sigma}=-J\tau_{\Sigma}$, we obtain
$$
 e(t)\cdot\nu_{\Sigma}
 =\tau_B\cdot\tau_{\Sigma}
 =\frac{\tau_B\cdot DT(x)\tau_B}{|DT(x)\tau_B|}>0.
$$
  Therefore
$S(y)\cdot\nu_{\Sigma}(y)=r e(t)\cdot\nu_{\Sigma}(y)>0$
on $\Sigma_r$. For any nonnegative $\varphi\in C^1(\overline\Lambda)$,
the divergence theorem on $\Lambda_r$ gives
\begin{equation}\label{eq:fluxr}
 0\le\int_{\Sigma_r}\varphi S\cdot\nu_{\Sigma_r}
 =\int_{\Lambda_r}\big(\varphi\,\text{tr} DS+\nabla\varphi\cdot S\big)\,d y.
\end{equation}
Note that the term $\nabla \varphi\cdot S$ is bounded. Since $DS>0$, the first term can be bounded by noting,
$$
 0\le\int_{\Lambda_\delta}\text{tr} DS\,d y
 =\int_{\partial\Lambda_\delta}S\cdot\nu_\delta\,d\mathcal H^1
 \le \mathcal H^1(\partial\Lambda_\delta).
$$
Here $\Lambda_\delta = \{x\in \Lambda\;|\;d(x,\partial \Lambda)>\delta\}$. Note that since $\Lambda$ was a $C^{2,\alpha}$ domain then $\mathcal H^1(\partial\Lambda_\delta)\to \mathcal H^1(\partial\Lambda)$. This also shows $S\in W^{1,1}(\Lambda,\mathbb{R}^2)$. Therefore since we have uniform bounds for any $r<1$, by taking limits and using the Gauss-Green theorem for $W^{1,1}$ functions we get,
$$0\le\int_{\Sigma}\varphi S\cdot\nu_{\Sigma}.$$
Therefore since $S$ was continuous up to the boundary we get \eqref{eq:weaknormal} and \eqref{eq:sinpositive}.

\end{proof}

We make the choice of $\theta$ clear in the following lemma.
\begin{lem}\label{lem:corridor}
The lift $\theta$ in Lemma~\ref{lem:boundary_theta} can be chosen so that
$$
 \beta(s)-\pi\le\theta(s)\le\beta(s)
 \quad\text{and}\quad L(s)\le\theta(s)\le R(s)
 \qquad(s\in\R).
$$
In particular, $\theta\in\Aw$.
\end{lem}
\begin{proof}
    The difference $d(s)=\beta(s)-\theta(s)$ is continuous and $\ell$-periodic, while $\sin d(s)\ge0$. Since both $\beta$ and $\theta$ have winding number $1$, we can choose a particular lift $\theta$ based on $\beta$ such that $0\leq d(s)\leq \pi$. The conclusion follows from the monotonicity of $\theta$.
\end{proof}

\subsection{Lewy Coordinates}\label{subsec:LewyT} A problem with the mixed domains while performing partial Legendre transform is that there is very little control of the boundary of the mixed domain. By using Lewy transform however we get good enough control of the boundary that allows us to pass to limits. For a general treatment of Lewy transforms see \cite{MS}. For $y\in\overline\Lambda$, set
$$
 x=S(y),\qquad
 m(y)=\frac{x+y}{2}\in\C,
 \qquad
 W(m(y))=\frac{\overline{y-x}}{2}\in\C.
$$
\begin{lem}
\label{lem:midpoint}
The midpoint map $m:\overline\Lambda\to\C$ is continuous and injective.  Its image of the interior, $M=m(\Lambda)$, is a Jordan domain, and $\overline M=m(\overline\Lambda)$.
Moreover, for all $y_1,y_2\in\overline\Lambda$,
\begin{equation}\label{eq:inverseLip}
 |y_1-y_2|\le2|m(y_1)-m(y_2)|,
\end{equation}
and $W$ is globally $1$-Lipschitz on $\overline M$.
The boundary $\partial M$ is rectifiable and admits the positive
parametrization
\begin{equation}\label{eq:muboundary}
 \mu(s)=m(\Gamma(s))=\tfrac12\big(\Gamma(s)+e^{i\theta(s)}\big),
\end{equation}
with
\begin{equation}\label{eq:lenM}
 \var(\mu)\le\frac{\ell+2\pi}{2}.
\end{equation}
\end{lem}
\begin{proof}
    Let $y_j\in\overline\Lambda$, $x_j= S(y_j)$, and $\Delta y=y_1-y_2$, $\Delta x=x_1-x_2$. Cyclic monotonicity gives $\Delta x\cdot \Delta y\ge0$.  Therefore
$$
 2\big(m(y_1)-m(y_2)\big)\cdot\Delta y
 =|\Delta y|^2+\Delta x\cdot\Delta y
 \ge |\Delta y|^2.
$$
Cauchy Schwarz proves \eqref{eq:inverseLip}, and in particular $m$ is injective.  Since the closed Jordan disk
$\overline\Lambda$ is compact, its continuous injective image is homeomorphic to it. On the interior, $Dm=(\Id+DS)/2$ has positive determinant, so $M$ is open. Therefore by \cite{Massey} we also get $m(\partial\Lambda) = \partial M$.\\
By definition $\overline{W(m_1)-W(m_2)}
=(\Delta y-\Delta x)/2$.  By parallelogram identity
$$
 |m_1-m_2|^2-|W(m_1)-W(m_2)|^2
 =\tfrac14\big(|\Delta y+\Delta x|^2
      -|\Delta y-\Delta x|^2\big)
 =\Delta x\cdot\Delta y\ge0.
$$
Thus $|W(m_1)-W(m_2)|\le|m_1-m_2|$ globally,
without assuming $M$ convex.
Finally, \eqref{eq:muboundary} follows from
Lemma~\ref{lem:boundary_theta}. The variation of the arclength curve $\Gamma$ is $\ell$ and that of $e^{i\theta}$ is $2\pi$ by Lemma~\ref{lem:BV}.
\end{proof}

We now study the displacement function $W$.
\begin{lem}\label{lem:holo_W}
The function $W:M\to\C$ is holomorphic, with $|W'(m)|<1$
in $M$. Moreover the functions
\begin{equation}\label{eq:zetaeta_holo_W}
 \zeta(m):=m+W(m),\qquad
 \eta(m):=m-W(m)
\end{equation}
are holomorphic in $M$ and continuous on $\overline M$. Their boundary values are given by
\begin{equation}\label{eq:traceid}
 \zeta(\mu(s))=\zeta_\theta(s),\qquad
 \eta(\mu(s))=\eta_\theta(s),
\end{equation}
where $\theta$ as in Lemma \ref{lem:boundary_theta}.
\end{lem}
\begin{proof}
    Fix $y\in\Lambda$ and write $A=DS(y)$, a symmetric positive definite with $\det A=1$. Let $J_0=\text{diag}(1,-1)$. Therefore,
    $$
 D_y m=\frac{\Id+A}{2},\qquad
 D_y W=\frac{J_0(\Id-A)}{2}.
$$
Since $\Id+A$ is invertible,
\begin{equation}\label{eq:DW}
 D_mW=J_0\,(\Id-A)(\Id+A)^{-1}=J_0B.
\end{equation}
If the eigenvalues of $A$ are
$\lambda$ and $\lambda^{-1}$, the eigenvalues of $B$ are
$$
 \frac{1-\lambda}{1+\lambda}
 \quad\text{and}\quad
 \frac{1-\lambda^{-1}}{1+\lambda^{-1}}
 =-\frac{1-\lambda}{1+\lambda}.
$$
Thus $B$ is symmetric and trace-free. So if in the standard basis
$B=\begin{psmallmatrix}a&b\\b&-a\end{psmallmatrix}$.
Equation \eqref{eq:DW} becomes
$$
 D_mW=\begin{pmatrix}a&b\\-b&a\end{pmatrix},
$$
and therefore $W(m)$ satisfies Cauchy Riemann equations. The other conclusions follow immediately from Lemma \ref{lem:midpoint} and definition of $\theta$.
\end{proof}

\subsection{Moment obstructions}We first note a simple complex analysis fact.
\begin{lem}
Let $M$ be a Jordan domain with rectifiable boundary and
$f\in C(\overline M)\cap\Hol(M)$.  Then
\begin{equation}
  \oint_{\partial M} f(m)\,d m=0.
\end{equation}
\end{lem}
The proof follows immediately from Mergelyan's theorem. We now prove the vanishing of moments.

\begin{proof}[Proof of Theorem \ref{thm:moment_obstruction}]
    Fix $k\ge0$ and parametrize $\partial M$ by
$\mu(s)$ as in \eqref{eq:muboundary}.
The change of boundary parametrization gives
$$
 M_k[\theta]=\oint_{\partial M}\eta\,\zeta^k\,d\zeta.
$$
 $\zeta\circ\mu=\zeta_\theta$ is BV by
Lemma~\ref{lem:BV}. Since $\eta=2m-\zeta$,
\begin{align*}
 M_k[\theta]
 &=2\oint_{\partial M}m\,\zeta^k\,d\zeta
     -\oint_{\partial M}\zeta^{k+1}\,d\zeta\\
 &=\frac{2}{k+1}\oint_{\partial M}m\,d(\zeta^{k+1})
\\
 &=-\frac{2}{k+1}\oint_{\partial M}\zeta(m)^{k+1}\,d m.
\end{align*}
By Lemma \ref{lem:holo_W}, $\zeta^{k+1}$ is homolorphic in $M$ and continuous till the boundary. Therefore $M_k[\theta]=0$.
\end{proof}

Using Theorem \ref{thm:moment_obstruction} and Remark \ref{rem:perturb_gen_convex}, Theorem \ref{thm:intro_thm2} immediately follows. We note a simple corollary of Theorem \ref{thm:intro_thm2}.
\begin{cor}\label{cor:chodosh_test}
    If there exists a point $\Gamma(s)\in \partial\Lambda$ such that $L(s)>R(s)$ then the optimal transport map from $B_1(0)$ to $\Lambda$ is discontinuous. 
\end{cor}
The above test was first studied by Wolfson in \cite{Wolfson} and refined by Urbas in \cite{Urbas}. A convex analysis based proof was also given in \cite{Chodosh} under the condition,
$$\inf_{J\subset \partial\Lambda}\int_J\kappa_{\partial\Lambda}<-\pi.$$
Note that $\int_a^b\kappa_{\partial\Lambda} = \beta(b)-\beta(a)$, where $\beta$ is the tangent angle. Hence,
$$L(b)\geq \beta(a)-\pi > \beta(b) \geq R(b).$$
Therefore Corollary \ref{cor:chodosh_test} is equivalent to the tests proposed in \cite{Chodosh}. Moreover the family of tests in Theorem \ref{thm:intro_thm2} provides a new family of conditions that detect discontinuities beyond obliqueness. A more symmetric version of the tests based on Lewy transform can also be recorded.
\begin{lem}\label{lem:lewy_tests}
    For $\theta\in\Aw$ consider,
    $$m_\theta(s) = \frac{\Gamma(s)+e^{i\theta(s)}}{2},\qquad W_\theta(s) = \frac{\overline{\Gamma(s)-e^{i\theta(s)}}}{2}.$$
    If for some $k\geq 0$,
    $$\mathcal{L}_k[\theta] = \int W_\theta m_\theta^kdm_\theta \neq 0 \qquad\text{for every }\theta\in\Aw$$
    then the optimal transport map from $B_1(0)$ to $\Lambda$ is discontinuous. 
\end{lem}
The above result follows immediately from Lemma \ref{lem:holo_W}. We now show that the infinite family of tests are all in fact necessary by considering the following examples.

\begin{ex}\label{ex:example_obstruction}
For each $k\ge0$, put
\begin{equation}\label{eq:constants}
 N=k+1,\quad M=10^4N2^N,\quad
 r(t)=\frac{4}{\sqrt{17}}\Bigl(1+\frac14\cos Nt+\frac14\cos Mt\Bigr),
\end{equation}
and let $\Lambda_k$ be the Jordan domain bounded by
$\Gamma(t)=r(t)e^{it}$. Then the following statements hold.
\begin{enumerate}
\item $\Lambda_k$ is $C^\infty$, strictly star-shaped, and has area $\pi$. The radial correspondence $\theta_0(t)=t$ has strictly positive twist,
so $\Aw\ne\varnothing$.
\item The radial correspondence satisfies
\begin{equation}\label{eq:low}
 \mathcal{L}_j[\theta_0]=0,\qquad 0\le j\le k=N-1.
\end{equation}
\item For \emph{every} $\theta\in\Aw$, one has
\begin{equation}\label{eq:margin}
\quad
 \Ima\mathcal{L}_N[\theta]
 \ge\left(\frac{\pi}{17}-\frac{18\pi (2\pi+12)}{10^4N}\right)h_-^N>0.
 \quad
\end{equation}
In particular, no single weakly oblique correspondence can annihilate all the
moments $\mathcal{L}_j$.
\end{enumerate}

\begin{figure}[h]
    \centering
    \begin{subfigure}[b]{0.49\textwidth}
    \centering
        \includegraphics[width=0.9\textwidth]{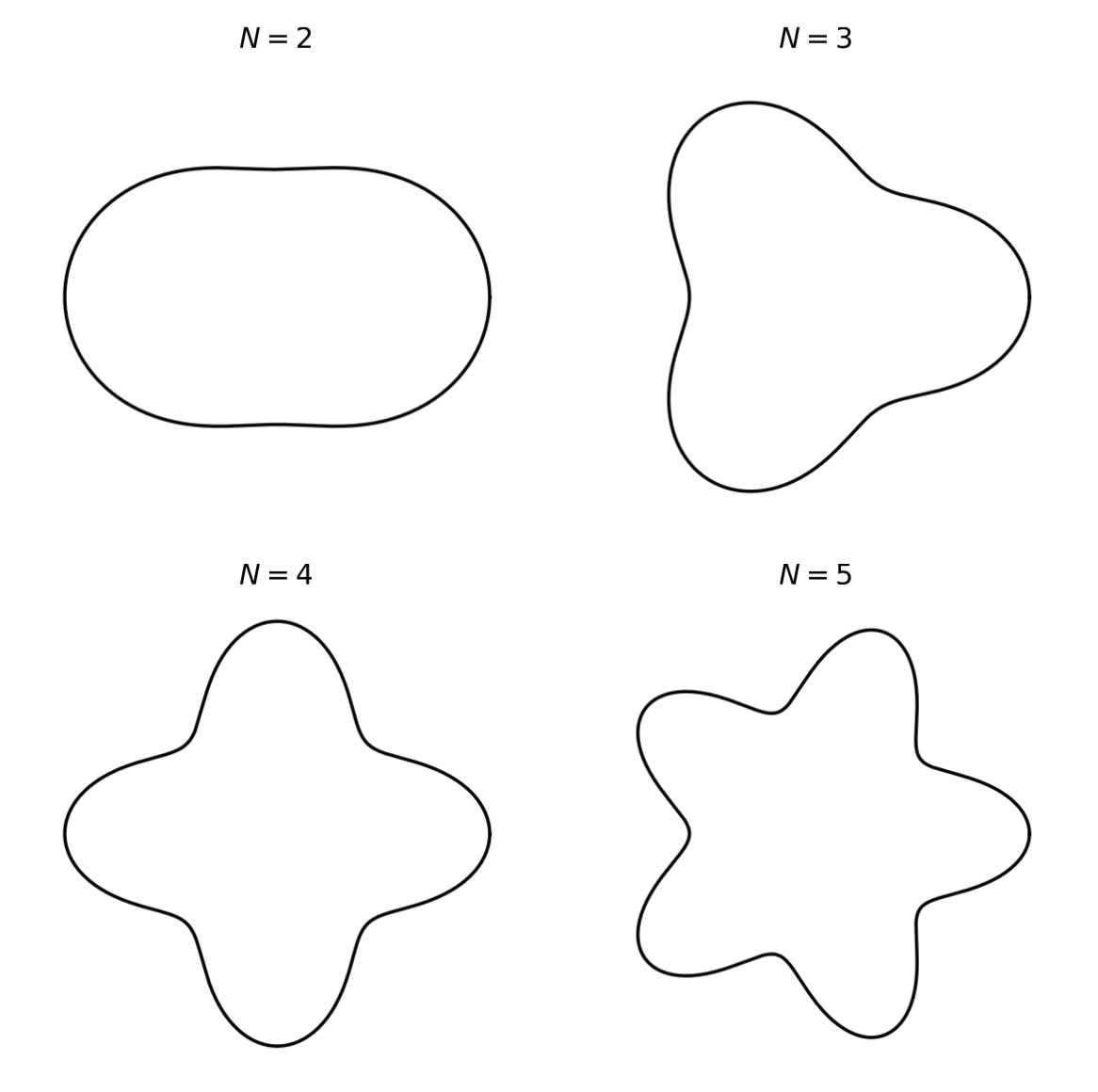}
    \caption{Domains $\Lambda_N$ without the high frequency term $\frac{1}{\sqrt{17}}\cos Mt$.}
    \end{subfigure}
    \hfill
    \begin{subfigure}[b]{0.49\textwidth}
    \centering
        \includegraphics[width=0.9\textwidth]{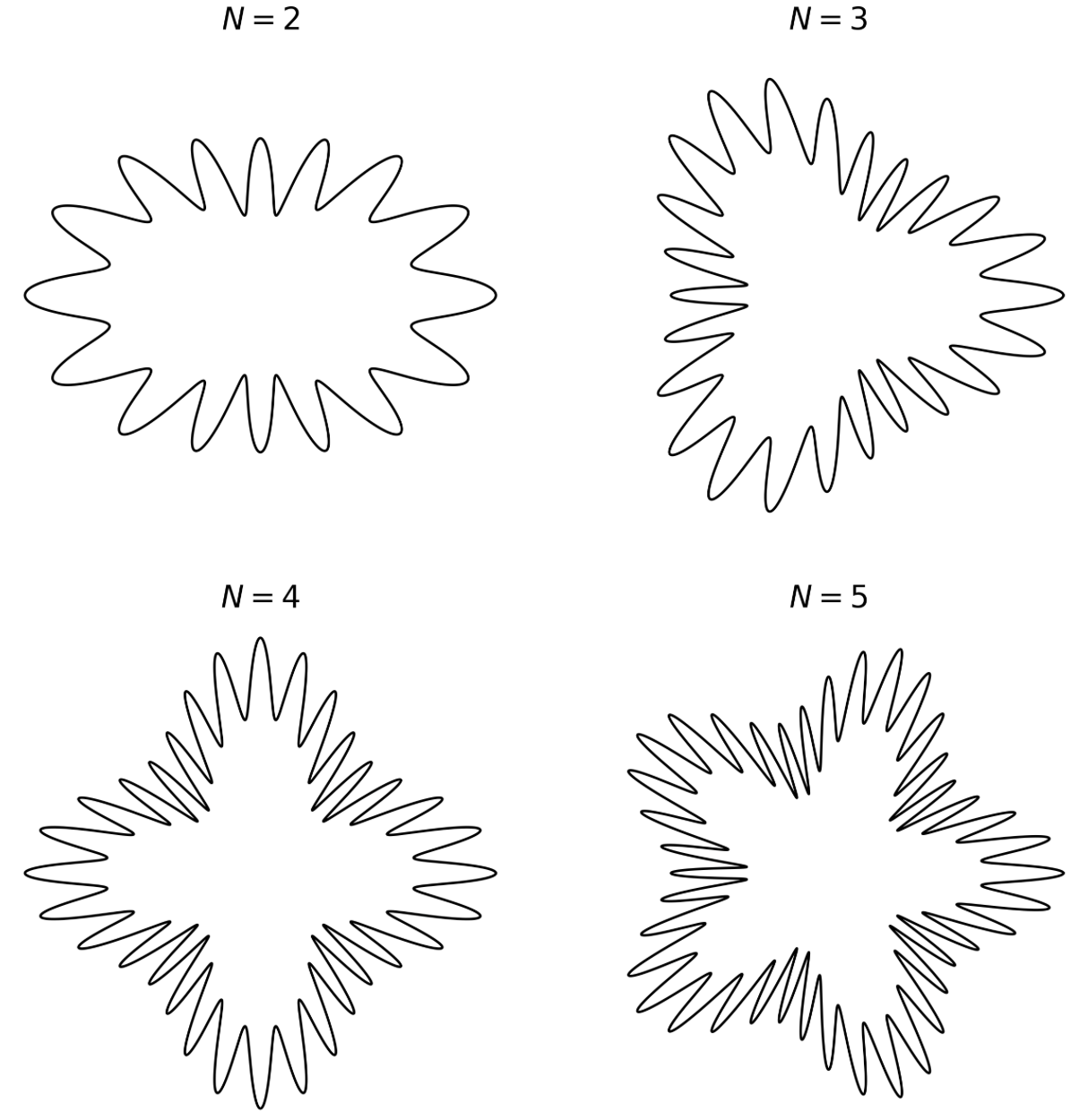}
    \caption{Illustrative figures for $\Lambda_N$ with the high frequency terms $\frac{1}{\sqrt{17}}\cos Mt$ for $M = 8N$.}
    \end{subfigure}
    
\end{figure}

 Note that,
$$
 \frac2{\sqrt{17}}\le r\le\frac6{\sqrt{17}},
 \qquad \frac12\int_0^{2\pi}r^2 \;d t=\pi,
$$
 Since
$$
 \Gamma'(t)=e^{it}(r'(t)+ir(t)),\qquad r(t)>0,
$$
we can write $\beta(t)=t+\alpha(t)$ with $0<\alpha(t)<\pi$, the phase for $r'(t)+ir(t)$.
For $0<t_2-t_1<2\pi$,
$$
 \beta(t_2)-\beta(t_1)
 =t_2-t_1+\alpha(t_2)-\alpha(t_1)>-\pi.
$$
This also shows for $\theta_0(t)=t$ obliqueness is always satisfied, i.e., $\theta_0\in \Aw$. Moreover 
$$
 r'(t)=-\frac{N\sin Nt+M\sin Mt}{\sqrt{17}}.
$$
Therefore every interval of length $2\pi/M$ contains a point with
$\sin Ms=+1$ and a point with $\sin Ms=-1$. At a point with $\sin Ms=+1$, $r'\le-(M-N)/\sqrt{17}$, hence
$\alpha\ge\pi-\delta$ for some small $\delta$. Similarly at a point with $\sin Ms=-1$, $r'\ge(M-N)/\sqrt{17}$, hence $\alpha\le\delta$. Here $\delta$ is given by,
$$
 \delta=\tan^{-1}\frac6{M-N}\le\frac{12}{M}.
$$
Choose a $+1$ point $s_-\in[t-2\pi/M,t]$ and a $-1$ point
$s_+\in[t,t+2\pi/M]$.  Monotonicity and the definitions of $L,R$ yield
$$
 t-\frac{2\pi}M-\delta
 \le s_-+\alpha(s_-)-\pi
 \le\theta(s_-)\le\theta(t)\le\theta(s_+)
 \le s_++\alpha(s_+)
 \le t+\frac{2\pi}M+\delta.
$$
Therefore every such $\theta$ obeys
\begin{equation}\label{eq:strip}
 \|\theta-\theta_0\|_{L^\infty([0,2\pi])}\le \frac{2\pi+12}{M}.
\end{equation}

We now show that for $\theta_0$ the first $k$ moments vanish and the $N = k+1$-th moment is quantitatively away from 0. Note,
$$\mathcal{L}_j[\theta_0] = \int_0^{2\pi}W_0(t)m_0(t)^{j}m_0'(t)dt.$$
Here $m_0,W_0$ are quantities corresponding to $\theta_0$ in Lemma \ref{lem:lewy_tests}. After a change of variables $t\mapsto t+2\pi/N$ we get,
$$\mathcal{L}_j[\theta_0] = e^{\frac{2\pi j}{N}i}\int_0^{2\pi}W_0(t)m_0(t)^{j}m_0'(t)dt.$$

Therefore for $j\not\equiv 0\mod{N}$ the moments must vanish. In particular for $j=1,2,\cdots,N-1$. Note that the $0^{th}$ moment is just the area constraint. For $j=N$ we explicitly calculate the moment. Set $h=(1+r)/2$ and
$$
 F_N(h):=\frac4{N+2}h^{N+2}-\frac2{N+1}h^{N+1}.
$$
For $\theta_0$, one has $m_0(t)=h(t)e^{it}$ and $2W_0(t)=(r(t)-1)e^{-it}$. Therefore integration by parts gives,
\begin{equation}\label{eq:radial}
\begin{aligned}
     2\mathcal{L}_N[\theta_0]=\int_0^{2\pi}&(r-1)h^N(h'+ih)e^{iNt} d t
       =i\int_0^{2\pi}e^{iNt}F_N(h(t)) d t\\
 & F_N'(h)=2(2h-1)h^N=2rh^N>0.
\end{aligned}
\end{equation}
Since $M/N=10^4 2^N$ is even, for $a=\pi/N$
$$
 h(t)-h(t+a)=\frac1{\sqrt{17}}\cos Nt.
$$
Taking the imaginary part of \eqref{eq:radial} and shifting $t$ by
$a$ gives
\begin{align}
 4\Ima \mathcal{L}_N[\theta_0]
 &=\int_0^{2\pi}\cos Nt\,
       \bigl(F_N(h(t))-F_N(h(t+a))\bigr)d t\notag\\
 &=\frac{1}{\sqrt{17}}\int_0^{2\pi}\cos^2 Nt F'_N(\xi_t) d t
 \label{eq:radiallower}
\end{align}
By the definition of $F'_N$, $\xi_t$ is bounded between $h(t),t(t+a)$ and $h \geq h_{-}$ we get that,

\begin{equation}\label{eq:moment_lower}
    \Ima \mathcal{L}_N[\theta_0] \geq \frac{\pi}{17}h_-^N>0. 
\end{equation}
Here $h_-=\frac{1+2/\sqrt{17}}2$. Now using equation \eqref{eq:strip} and a varaition of $\mathcal{L}_N$ we show that $\mathcal{L}_N[\theta]>0$ for all admissible $\theta$. The following lemma calculates the variation of the $N$th moment and may be of independent interest.
\begin{lem}
    For independent variables $t,q$, define
$$
 m(t,q)=\frac{\Gamma(t)+e^{iq}}2,\quad
 \Delta(t,q)=r(t)\cos(q-t)-r'(t)\sin(q-t).
$$
For $\theta\in\Aw$,
\begin{equation}\label{eq:variation}
\mathcal{L}_N[\theta]-\mathcal{L}_N[\theta_0]
 =-\frac{1}{2}\int_0^{2\pi}\int_{\theta_0(t)}^{\theta(t)}
 \Delta(t,q)m(t,q)^Nd qd t.
\end{equation}
\end{lem}
\begin{proof}
    First take a smooth pairing $\theta$, put
$\theta_s(t)=\theta_0(t)+s(\theta(t)-\theta_0(t))$, $\phi(t)=\theta(t)-\theta_0(t)$. Also define,
$$w(t,q)=\overline{y(t)-e^{iq}}.$$
Evaluate $m,w$ at $q=\theta_s(t)$.  Differentiating the integral for $\mathcal{L}_N[\theta_s]$ in $s$ and integrating by parts in $t$,
$$
 \frac d{ds}\mathcal{L}_N[\theta_s]
 =\frac{1}{2}\int_0^{2\pi}m^N\big(w_s m_t-w_t m_s\big)d t.
$$
Using
$$
 m_s=\frac{i}{2}e^{i\theta_s}\phi,
 \qquad w_s=ie^{-i\theta_s}\phi,
 \qquad \Gamma'=e^{it}(r'+ir),
$$
we get,
$$
 w_s m_t-w_t m_s
 =-\phi\,\operatorname{Im}(e^{-i\theta_s}\Gamma')
 =-\phi\,\Delta(t,\theta_s).
$$
Integrating in $s$ and substituting $q=\theta_0(t)+s\phi(t)$ proves \eqref{eq:variation}. A general $\theta\in\Aw$ can be approximated by a smooth nondecreasing and since the total variation remains bounded the as in Lemma \ref{lem:midpoint} the moments $\mathcal{L}_N$ converge.

\end{proof}

The nonvanishing of moments is now immediate from the estimates \eqref{eq:strip}, \eqref{eq:moment_lower} and the formula \eqref{eq:variation}.

\end{ex}

\section{Smoothness preserved under \texorpdfstring{$C^{1,\alpha}$}{C1,alpha} deformations}
\label{sec:perturb}

In this section we prove the following theorem.

\begin{thm}\label{thm:openness}
Let $\Omega_0,\Lambda_0$ be bounded $C^{2,\alpha}$ domains of same area and assume that $\Omega_0$ is convex. Assume there is a $C^{2,\alpha}(\bar\Omega_0)$ optimal transport map $T_0: \Omega_0\to\Lambda_0$. Consider a $C^{1,\alpha}$ domain $\Lambda$ whose boundary is a normal graph, denoted by $h$, over $\partial\Lambda_0$. If $\|h\|_{C^{1,\alpha}}$ was small such that $|\Lambda| = |\Omega_0|$ then there exists an optimal transport map $T:\Omega_0\to\Lambda$ such that for any $\beta<\alpha$, $T\in C^{1,\beta}(\bar\Omega_0)$.
\end{thm}
\begin{rem}
    Note that in the above theorem no convexity assumption is needed either on $\Lambda_0$ or $\Lambda$.
\end{rem}
\begin{rem}
    The Theorem also holds if we allow a $C^{1,\alpha}$ perturbation of the source as long as it remains convex, thus establishing Theorem \ref{thm:intro_thm3}.
\end{rem}

To prove Theorem \ref{thm:openness} we have to deal with moving boundaries. So to fix a domain in $\mathbb{C}$ we use the Riemann mapping theorem to map simply connected domains to the unit disc, $\D$. Since the automorphism group of $\D$ is nontrivial we also normalize it, i.e., have one predefined fixed point and make its derivative there real. We now fix some terminology and setup a fixed point argument.\\

\subsection{Setup}\label{subsec:KW_step} Let $D_0$ be the reference mixed domain and $H_0$ the reference holomorphic function. The hypotheses imply $D_0$ has $C^{2,\alpha}$ boundary and $H_0\in C^{2,\alpha}(\bar D_0)$. Let $Z_0:\D\to D_0$ be a normalized conformal map. By the  Kellogg--Warschawski theorem, $Z_0\in C^{2,\alpha}(\bar\D)$ and $Z_0'\ne0$ on $\bar\D$. Also define $W_0=H_0\circ Z_0\in C^{2,\alpha}(\bar\D)$.

For holomorphic functions $Z,W$ on $\D$ define
\begin{equation}\label{eq:S-R}
 S_{Z,W}=(\Rea Z,\Ima W),\qquad R_{Z,W}=(\Rea W,\Ima Z).
\end{equation}
Here $Z$ will represent the conformal map from $\D$ to an arbitrary mixed domain and $W$ will represent the generator of an arbitrary Optimal transport map with coodinates fixed to $\D$. \\
For $J(Z,W)=\Rea(W'\overline{Z'})$, Cauchy--Riemann equations give
\begin{equation}\label{eq:Jac}
 \det DS_{Z,W}=\det DR_{Z,W}=J(Z,W).
\end{equation}
More explicitly, if $a=\Rea(W'\overline{Z'})$, $b=\Ima(W'\overline{Z'})$, then
\begin{equation}\label{eq:Hessdisk}
 D(R_{Z,W}\circ S_{Z,W}^{-1})
 =\frac1a\begin{pmatrix}|W'|^2&-b\\-b&|Z'|^2\end{pmatrix};
 \quad \det D(R\circ S^{-1})=1.
\end{equation}
Indeed $|Z'|^2|W'|^2=a^2+b^2$, and $a>0$ makes the matrix positive definite. We fix conformal automorphism freedom by
\begin{equation}\label{eq:gauge}
 Z(0)=Z_0(0),\qquad \Ima\left(\frac{Z'(0)}{Z_0'(0)}\right)=0.
\end{equation}
We also fix the source domain to be $\Omega_0$. Let $\tilde{\rho}_0$ be a $C^{2,\alpha}$ function such that the level curve $\tilde{\rho_0}=0$ is $\partial\Omega_0$ and $\nabla \tilde{\rho_0} = \nu_{\Omega_0}$ on the boundary $\partial\Omega_0$. Therefore we ensure, 
\begin{equation}\label{eq:fix_source}
\tilde{\rho_0}(S_{Z,W})=0.
\end{equation}

Similarly choose a $C^{2,\alpha}$ defining function $\rho_0$ of $\partial \Lambda_0$, such that $\rho_0=0$ and $\nabla\rho_0=\nu_{\Lambda_0}$ on its boundary. We define both so that their $C^{2,\alpha}$ norm are bounded on $\R^2$. A nearby $C^{1,\alpha}$ normal graph target can be written
\begin{equation}\label{eq:rhoh}
 \Lambda_h=\{\rho_0+h<0\},\qquad \|h\|_{C^{1,\alpha}}\ll1,
\end{equation}
near the reference boundary. Conversely every small $C^{1,\alpha}$ normal graph admits such an $h$, and the norms are comparable locally. For a scalar $c\in \mathbb{R}$, we now fix the boundary of the target domain by
\begin{equation}\label{eq:fix_target}
\rho_0(R_{Z,W})+h(R_{Z,W})-c=0.   
\end{equation}
Note that we need to include this parameter $c$ to allow us to maintain fixed area constraints. We now setup the fixed point problem we would like to solve.\\

Fix $0<\gamma<\alpha$, and let $\mathcal A^{1,\gamma}=\Hol(\D)\cap C^{1,\gamma}(\bar\D)$ be the Banach space of all Holomorphic functions with $C^{1,\gamma}$ boundary values. Set
$$
 \mathcal X_\gamma=\mathcal A^{1,\gamma}\times\mathcal A^{1,\gamma}\times\R,
 \qquad
 \mathcal Y_\gamma=C^{1,\gamma}(\partial\D;\R)^2\times\C\times \R.
$$
Write $X=(Z,W,c)$, $X_0=(Z_0,W_0,0)$, and define the smooth-reference operator
\begin{equation}\label{eq:B0}
 \mathcal B_0(X)=\left(
 \tilde{\rho}_0(S_X),\ \rho_0(R_X)-c,Z(0)-Z_0(0),\
 \Ima\frac{Z'(0)}{Z_0'(0)}\right).
\end{equation}
We have $\mathcal B_0(X_0)=0$. The target perturbation is defined to be
\begin{equation}\label{eq:Eh}
 E_h(X)=\left(0,h(R_X),0,0\right).
\end{equation}
Thus for a perturbation to have a $C^{1,\gamma}$ optimal transport, we have to solve the following nonlinear problem,
\begin{equation}\label{eq:nonlinear}
    \mathcal B_0(X)+E_h(X)=0.
\end{equation}

\subsection{The Operator \texorpdfstring{$\mathcal{B}_0$}{B0}} We first study the operator $\B_0$.

\begin{lem}\label{lem:B0_diff}
The map $\mathcal B_0:\mathcal X_\gamma\to\mathcal Y_\gamma$ is continuously Fr\'echet differentiable.
\end{lem}
\begin{proof}
First note the last two components of $\B_0$, $\B_0|_{\C\times \R} = ( Z(0)-Z_0(0),\Ima \frac{Z'(0)}{Z_0'(0)})$ are linear therefore continuously differentiable. The first two components $\B_0|_{C^{1,\gamma}(\partial\D;\R)^2}$ can be written as a composition of two maps $X\mapsto (S,R,c)$ and $(S,R,c)\mapsto (\tilde\rho_0(S),\rho_0(R)-c)$. Since the former map is linear, we have to study the differentiability of the latter.\\
We note a simple observation. Let $f:\R^2\to \R$ be a $C^{2,\alpha}$ function. Then the mapping $A \in \mathcal A^{1,\gamma} \mapsto f(A)\in C^{1,\gamma}(\partial\D;\R)$ is continuously differentiable for $\gamma < \alpha$. Consider the error term
$$N(\delta) = f(A+\delta)-f(A)-\nabla f(A)\cdot\delta = \mathcal{O}(\|\delta\|_{\mathcal A^{1,\gamma}}^2). $$
 Further note that on $\partial\D$, 
\begin{align*}
    N'(\delta) &= (\nabla f(A+\delta)-\nabla f(A))\cdot \delta'+(\nabla f(A+\delta)-\nabla f(A))\cdot A'-D^2f(A)A'\cdot\delta\\
   & = \int_0^1 \left<\delta'+A', D^2f(A+t\delta)\delta\right>dt-D^2f(A)A'\cdot\delta
\end{align*}
Since $D^2f \in C^{\alpha}$ and $A,\delta \in \mathcal A^{1,\gamma}$, we get $N'(\delta)\in C^{0,\gamma}(\partial\D,\R)$. Moreover
$$\|N(\delta)\|_{C^{1,\gamma}(\partial\D,\R)}  = o(\|\delta\|_{\mathcal A^{1,\gamma}}).$$
This shows the first two components individually are continuously differentiable as well.
\end{proof}

Let
\begin{equation}\label{eq:L}
 L=D\mathcal B_0(X_0):\mathcal X_\gamma\to\mathcal Y_\gamma.
\end{equation}
Consider,
$$
 N(X)=\mathcal B_0(X)-\mathcal B_0(X_0)-L(X-X_0)
 =\mathcal B_0(X)-L(X-X_0).
$$
By differentiability of $\B_0$ there exists a modulus of continuity $\omega_{\gamma}$ such that,
$$\|N(X)\|_{\mathcal Y_\gamma}\le\omega_\gamma(\|X-X_0\|_{\mathcal X_\gamma})\|X-X_0\|_{\mathcal X_\gamma}.$$
We linearize \eqref{eq:nonlinear} to get,
\begin{equation}
    L(X-X_0)+N(X)+E_h(X)=0
\end{equation}
Therefore if $L$ is invertible at $X_0$, then for a given $h$ we would like to solve the following fixed point problem,
\begin{equation}\label{eq:Th}
 Y=\mathcal T_h(Y),\qquad
 \mathcal T_h(Y)=-L^{-1}\bigl(N(X_0+Y)+E_h(X_0+Y)\bigr).
\end{equation}.

\subsection{Invertibility of \texorpdfstring{$L$}{L}} The invertibility of $L$ is closely related to solving a Neumann boundary problem of the Laplace equation. Let $(z,w,l)$ be a variation of $X$ near $X_0$. Therefore,
\begin{equation}\label{eq:dotSR}
 \dot S=(\Rea z,\Ima w),\qquad \dot R= (\Rea w,\Ima z).
\end{equation}
Along $\partial\D$, write $\nu=\nu_{\Omega_0}(S_0) = (\nu_1,\nu_2)$ and $\nu^*=\nu_{\Lambda_0}(R_0) = (\nu_1^*,\nu_2^*)$. For prescribed real $f_1,f_2\in C^{1,\gamma}(\partial\D)$, the first two linear boundary equations are
\begin{equation}\label{eq:linearBC}
 \nu\cdot\dot S=f_1,\qquad \nu^*\cdot\dot R-\ell=f_2.
\end{equation}
We first reconstruct a zero gauge solution and then treat an arbitrary gauge data case. Consider the following holomorphic variation lemma,
\begin{lem}\label{lem:phi_recover}
    Let $\varphi$ be the holomorphic variation required in $H_0$ for the prescribed data $f_1,f_2\in C^{1,\gamma}(\partial\D)$. Parametrize $\partial \Omega_0$ by positively oriented arc length $s$ and let $\zeta: \partial \Omega_0\to \partial D_0$ be the source domain to mixed domain mapping. Further if $\sigma$ is such that $ds_{\Lambda_0} = \sigma ds$, then there exists a harmonic function $k$ such that
    $$\varphi = k_x-ik_q,$$
    and 
    $$\partial_{\nu_D}k=\frac{\sigma(f_2+\ell)-f_1}{|\zeta_s|}\ \text{on }\partial D_0.$$
    
\end{lem}
\begin{proof}
    Let $G_0(x,q)=(x,\Ima H_0(x+iq))$, $R_{H_0}(x,q)=(\Rea H_0(x+iq),q)$ and $B_0=DT_0$. write $H_0'=a+ib$ with $a>0$. Define the holomorphic function (representing the holomorphic change in $H_0$) on $\D$
$$
 \Phi=w-H_0'(Z_0)z,
$$
and consider its analogue on $D_0$ by $\varphi=\Phi\circ Z_0^{-1}$. The variations split as
$$
 \dot S=DG_0\,z+(0,\Ima\varphi),\qquad
 \dot R=DR_{H_0}\,z+(\Rea\varphi,0),
$$
 Because $B_0=DR_{H_0}(DG_0)^{-1}$, the physical transport variation
\begin{equation}\label{eq:V}
 V:=\dot R-B_0\dot S
 =\left( \Rea\varphi+\frac ba\Ima\varphi,-\frac{\Ima\varphi}{a}\right)
\end{equation}
is independent of $z$. Since $D_0$ is simply connected, there is a real harmonic $k$ (unique up to a constant) with
$$
 \varphi=k_x-ik_q.
$$
Thus
\begin{equation}\label{eq:Vgrad}
 V=\left(\ k_x-\frac ba k_q,\frac{k_q}{a}\right).
\end{equation}

We now note the boundary values of the harmonic function $k$. Let $\tau=(\tau_1,\tau_2)$ be the tangent of $\partial \Omega_0$. Denote the corresponding target tangent and normal by $\tau^*=J\nu^*$ and $\nu^*=(m_1,m_2)$, where $J(x,y)=(-y,x)$. As the boundary map is orientation preserving, there is a strictly positive boundary stretching factor $\sigma$ such that
\begin{equation}\label{eq:Btangent}
 B_0\tau=\sigma\tau^*,\qquad ds_{\Lambda_0}=\sigma\,ds.
\end{equation}
The determinant-one matrix identity $B_0^{-1}=\text{cof} B_0=J^\mathsf TB_0J$ gives
$$
 B_0^{-1}\nu=\sigma\nu^*,\qquad
 B_0\nu^*=\frac1\sigma\nu.
$$
The mixed auxiliary boundary $\zeta(s)=G_0^{-1}(e^{is})=x(s)+iq(s)$ satisfies
$$
 x_s=\tau_1 ,\qquad q_s=\frac{\tau_2-b\tau_1}{a}.
$$
Because $B_0\tau=\sigma(-m_2,m_1)$, comparison of components gives
\begin{equation}\label{eq:mixedtangent}
 q_s=\sigma m_1,\qquad
 ax_s=\sigma(bm_1-m_2).
\end{equation}
The auxiliary outward normal is $\nu_D=(q_s,-x_s)/|\zeta_s|$. Inserting \eqref{eq:mixedtangent} into \eqref{eq:Vgrad} gives
\begin{equation}\label{eq:conormalhol}
 \nu^*\cdot V
 =\frac{m_2-bm_1}{a}k_q+m_1k_x
 =\frac{q_sk_x-x_sk_q}{\sigma}
 =\frac{|\zeta_s|}{\sigma}\partial_{\nu_D}k.
\end{equation}
On the other hand, the linear boundary equations \eqref{eq:linearBC} and $B_0\nu^*=\nu/\sigma$ give
$$
 \nu^*\cdot V
 =\nu^*\cdot\dot R-(B_0\nu^*)\cdot\dot S
 =f_2+\ell-\frac{f_1}{\sigma}.
$$
Combining these identities yields the harmonic Neumann problem
\begin{equation}\label{eq:neumann}
 \quad \Delta k=0\ \text{in }D_0,\qquad
 \partial_{\nu_D}k=\frac{\sigma(f_2+\ell)-f_1}{|\zeta_s|}\ \text{on }\partial D_0.\quad
\end{equation}
The boundary functions $f_i$ are transferred from $\partial\D$ to physical source arclength by the fixed reference boundary diffeomorphism. 
\end{proof}

Using the above lemma we first show that for the zero gauge conditions, the inverse can be constructed.

\begin{lem}\label{lem:uniqueL}
    Let $(f_1,f_2,0,0)\in \mathcal{Y}_{\gamma}$. There exists unique $(z,w,l) \in \mathcal{X}_{\gamma}$ such that $L(z,w,l) = (f_1,f_2,0,0)$.
\end{lem}
\begin{proof}
We know from Lemma \ref{lem:phi_recover} that the holomorphic variation in $H_0$ requires solving a Neumann problem. A Neumann problem is solvable if and only if,
$$
 0=\int_{\partial D_0}\partial_{\nu_D}k\,ds_D
   =\int_{\partial \Omega_0}\bigl(\sigma(f_2+\ell)-f_1\bigr)\,ds.
$$
Since $\int_{\partial \Omega_0}\sigma f_2\,ds=\int_{\partial\Lambda_0}f_2\,ds_{\Lambda_0}$, it uniquely determines
\begin{equation}\label{eq:ell}
 \ell=\frac{\int_{\partial \Omega_0}f_1\,ds-\int_{\partial\Lambda_0}f_2\,ds_{\Lambda_0}}{|\partial\Lambda_0|}.
\end{equation}
With this choice, the classical Neumann theorem on the $C^{2,\alpha}$ domain $D_0$ gives $k\in C^{2,\gamma}(\bar D_0)$ unique modulo constants, with
$$
 \|k-k_{D_0}\|_{C^{2,\gamma}(\bar D_0)}+|\ell|
 \le C\bigl(\|f_1\|_{C^{1,\gamma}}+\|f_2\|_{C^{1,\gamma}}\bigr).
$$
Therefore $\varphi=k_x-ik_q\in C^{1,\gamma}(\bar D_0)$ is uniquely determined and satisfies the corresponding norm estimate. Now let $G_0(x,q)=(x,\Ima H_0(x+iq))$, $R_{H_0}(x,q)=(\Rea H_0(x+iq),q)$ and $B_0=DT_0$. Write $H_0'=a+ib$ with $a>0$. The variations are determined by
$$
 \dot S=DG_0\,z+(0,\Ima\varphi),\qquad
 \dot R=DR_{H_0}\,z+(\Rea\varphi,0).
$$
Note that using $\varphi$ and the equations \eqref{eq:dotSR} $w$ can be recovered if $z$ is known. We write the first boundary equation in \eqref{eq:linearBC},
$$
 \nu\cdot DG_0\,z=f_1-\nu_2\Ima\varphi.
$$
Since the reference $G_0$ is an orientation preserving diffeomorphism we get,
$$
 DG_0^T T\nu=\kappa\nu_D,\qquad \kappa>0.
$$
Consequently
\begin{equation}\label{eq:normz}
 \nu_D\cdot z=\frac{f_1-\nu_2\Ima\varphi}{\kappa}
 \quad\text{on }\partial D_0.
\end{equation}
Write $z(\zeta)=Z_0'(\zeta)q(\zeta)$ for a holomorphic $q$. Along $|\zeta|=1$, conformality gives
$$
 \nu_D(Z_0(\zeta))=\frac{\zeta Z_0'(\zeta)}{|Z_0'(\zeta)|},
 \qquad
 \nu_D\cdot Z_0'q=|Z_0'|\Rea(\bar \zeta q).
$$
Thus \eqref{eq:normz} is the scalar Schwarz boundary equation
\begin{equation}\label{eq:schwarz}
 \Rea(\bar \zeta q(\zeta))=\psi(\zeta),\quad
 \psi=\frac{f_1-\nu_2\Ima\varphi}{\kappa|Z_0'|}\in C^{1,\gamma}(\mathbb{S}^1;\R).
\end{equation}
The Schwarz integral
$$
 E(\zeta)=\frac1{2\pi}\int_0^{2\pi}
 \frac{e^{it}+\zeta}{e^{it}-\zeta}\psi(e^{it})\,dt
$$
defines the holomorphic $C^{1,\gamma}$ function with $\Rea E|_{\mathbb{S}^1}=\psi$ and $E(0)\in\R$. Set
$$
 q(\zeta)=\zeta E(\zeta),\qquad z=Z_0'q,
 \qquad w=H_0'(Z_0)z+\varphi(Z_0).
$$
Then \eqref{eq:schwarz} holds, $z(0)=0$, and $\Ima(z'(0)/Z_0'(0))=\Ima q'(0)=0$, so the gauge variation vanishes.
\end{proof}

We can now easily conclude the invertibility of $L$,

\begin{lem}
    The linear operator $L =D\mathcal B_0(X_0)$ is invertible.
\end{lem}
\begin{proof}
    By linearity of $L$ and the uniqueness result in Lemma \ref{lem:uniqueL} for $(f_1,f_2)=(0,0)$ we know $L$ is injective. Let $(f_1,f_2,\zeta_0,r)\in \mathcal{Y}_{\gamma}$. By Lemma \ref{lem:uniqueL} we first find $(z_1,w_1,l_1)\in \mathcal{X}_{\gamma}$ such that $L(z_1,w_1,l_1) = (f_1,f_2,0,0)$. Therefore we only need to find $(z_2,w_2,l_2)$ such that $L(z_2,w_2,l_2) = (0,0,\zeta_0,r)$. We actually find the preimage $(z_2,w_2,0)$ since $f_1,f_2 = 0$ forces $l=0$. We find a gauge $q$ such that,
    $$z_2 = Z_0'q,\;\text{ and }w_2 = H_0'Z_0'q = W_0'q.$$
    Here we can find $q$ such that,
    \begin{equation}\label{eq:automorphism}
 q(w)=b+itw-\bar b\,w^2,\qquad b\in\C,\quad t\in\R.
\end{equation}
$z(0) = \zeta_0$ forces $b = \frac{\zeta_0}{Z_0'}$ while
$$
 \Ima\frac{z_{2}'(0)}{Z_0'(0)}
 =\Ima\left(\frac{Z_0''(0)}{Z_0'(0)}b+it\right)
$$
uniquely determines $t$.
\end{proof}

\subsection{The perturbation $E_h$ and existence}
Unlike $\mathcal B_0$, the perturbation $E_h$ need not be $C^1$ in $X$ on the same $C^{1,\gamma}$ Banach space. Differentiating $h(R_X)$ in $X$ would introduce $Dh(R_X)$, whose spatial derivative requires $D^2h$. We shall only use its $C^{1,\alpha}$ behaviour.
\begin{lem}\label{lem:composition}
Fix $0<\gamma<\alpha$. There is a small ball $\|X-X_0\|_{\mathcal X_\gamma}\le r_0$ and a constant $C$ such that, for every $h\in C^{1,\alpha}$,
\begin{equation}\label{eq:compositionestimate}
 \|h(R_X)\|_{C^{1,\gamma}(\mathbb{S}^1)}\leq C\|h\|_{C^{1,\alpha}}.
\end{equation}
Moreover, if $0<\beta<\gamma$, $X_j\to X$ in $\mathcal X_\beta$, and $\sup_j\|X_j-X_0\|_{\mathcal X_\gamma}\le r_0$, then
\begin{equation}\label{eq:compositioncontinuity}
 h(R_{X_j})\longrightarrow h(R_X)\quad\text{in }C^{1,\beta}(\mathbb{S}^1).
\end{equation}
\end{lem}
The proof of the above lemma follows from basic calculus. We now focus on the fixed point equation \eqref{eq:Th}. We first find an invariant ball in $\mathcal{X}_{\gamma}$.

\begin{lem}\label{lem:inv_ball}
Let $0<\beta<\gamma<\alpha$. There exist $r>0$ and $\varepsilon>0$ such that, if $\|h\|_{C^{1,\alpha}}\le\varepsilon$, then
$$
 \mathcal T_h(K_r)\subset K_r,\qquad
 K_r=\{Y\in\mathcal X_\gamma:\|Y\|_{\mathcal X_\gamma}\le r\}.
$$
Moreover the set $K_r$ is compact and convex in $\mathcal{X}_{\beta}$ and the mapping $\mathcal{T}_h$ is continuous on it.
\end{lem}
\begin{proof}
Set $M=\|L^{-1}\|_{\mathcal Y_\gamma\to\mathcal X_\gamma}$. By \eqref{eq:Th},
$$
 \|\mathcal T_h(Y)\|_{\mathcal X_\gamma}
 \le M\omega_\gamma(r)r+MC\|h\|_{C^{1,\alpha}}.
$$
For $r$ sufficiently small ensure $M\omega_\gamma(r)\le1/4$ and choose $\varepsilon$ such that $MC\varepsilon\le r/4$. The right side is at most $r/2$, so the entire closed ball maps into itself. The compactness of $K_r\subseteq \mathcal{X}_{\beta}$ follows from compactness of $C^{1,\gamma}\hookrightarrow C^{1,\beta}$. The continuity follows from Lemma \ref{lem:composition}.
\end{proof}

We now conclude with the proof of the openness result,

\begin{proof}[Proof of Theorem \ref{thm:openness}]
By Schauder fixed point we know for every $h$ such that $\|h\|_{C^{1,\alpha}}$ small there exists $X_h = (Z_h,W_h,c_h)$ such that $\B_0(X_h)+E_h = 0$. For a small $h$ it can also be ensured that, 
$$
J_h=\Rea(W_h'\overline{Z_h'})\ge\frac12\min_{\bar\D}\Rea(W_0'\overline{Z_0'})>0.
$$
Thus using \eqref{eq:Hessdisk}, $S_h,R_h$ are orientation-preserving local $C^{1,\gamma}$ diffeomorphisms from $\mathbb{D}$ to $\Omega_0$ and $\tilde{\Lambda}_h = \{\rho_0+h<0\}$ respectively. By \cite{Massey} they are actually $C^{1,\gamma}(\overline{\mathbb{D}})$ diffeomorphisms. By discussion above this gives us a $C^{1,\beta}(\overline{\Omega_0})$ diffeomorphism,
$$T_h: \Omega_0\to \Lambda_h = \{\rho_0+h<c_h\},$$
for any $0<\beta<\alpha$. The map is given by $T_h = R_{Z_h,W_h}\circ S_{Z_h,W_h}^{-1}$. However the perturbation was made in way such that the area of $|\Lambda_h|$ is constant. Therefore $c_h$ has to be zero.
\end{proof}

\section{Optimal transport and Harmonic maps}\label{sec:harmonic_maps}
In this section we present another characterization based on the ideas of Lewy transform (used in Section \ref{subsec:LewyT}) and Kellogg–Warschawski theorem (used in Section \ref{subsec:KW_step}).  We always assume the target domain $\Lambda$ is a Jordan domain of area $\pi$. Recall that given a optimal transport map $T:B_1(0)\to \Lambda$, continuous till the boundary, we denote,
$$
 y=T(x),\qquad
 m(x)=\frac{x+y}{2}\in\C,
 \qquad
 W(m(x))=\frac{\overline{y-x}}{2}\in\C.
$$
From Lemma \ref{lem:midpoint} we have the map $m: \bar B_1(0)\to \C$ is continuous and injective with $W$ being holomorphic as a function of $m$ on the midpoint domain $\Dm = m(B_1(0))$.
 For $k\geq 1$ and $0<\alpha<1$, consider
$$
 \mathcal R^{k,\alpha}=\{\Lambda\;|\;\exists \;T:\overline B_1(0)\to\overline\Lambda\text{ optimal transport map, }T\in C^{k,\alpha}(\overline{B}_1(0)), DT>0\}.
$$
Also define the following continuous analogue
$$
\mathcal{C} = \{\Lambda\;|\;\exists \;T: B_1(0)\to\Lambda\text{ continuous optimal transport map}\}.
$$
Let $J$ be rotation by $\pi/2$, and on $\T=\R/(2\pi\mathbb Z)$ set
\[
 \H e^{ikt}=-i\sgn(k)e^{ikt},\quad K=\H \circ\partial_t = (-\partial_{tt})^{\frac{1}{2}}.
\]
Here $\H$ is the Hilbert transform on the circle. For vector valued inputs it acts componentwise. $K$ is the Dirichlet to Neumann operator on the unit disc for $\Delta$. Also denote $P$ to be the Poisson extension operator taking boundary values on the circle to their harmonic extensions. We prove the following generalization of Theorem \ref{}.

\begin{thm}\label{thm:boundary_charac_unit_circle}
For $k\geq 1$ and $0<\alpha<1$ let $\mathcal U_+^{k,\alpha}=\{\theta\in C^{k,\alpha}_{loc}(\R):\theta(t+2\pi)=\theta(t)+2\pi,\ \theta'>0\}$, modulo additions of $2\pi$. For every $\theta\in\mathcal U_+^{k,\alpha}$, the boundary values
\begin{equation}\label{eq:traces_charac}
 x=e^{i\theta},\qquad y=J\H x+c,\qquad c\in\R^2,
\end{equation}
produce a positively oriented embedded boundary $\partial\Lambda=y(\T)$, $|\Lambda|=\pi$, such that $\Lambda\in\mathcal R^{k,\alpha}$ and 
$$ T=P[y]\circ P[x]^{-1}: B_1(0)\to \Lambda $$
is the $C^{k,\alpha}$ non degenrate optimal transport map. Conversely every $\Lambda\in\mathcal R^{k,\alpha}$ arises this way.
\end{thm}

\begin{proof}
    We first show that given $T$ we can always find a $\theta\in \mathcal{U}_{+}^{k,\alpha}$ so that \eqref{eq:traces_charac} give the boundary correspondence. Consider the simply connected midpoint domain $\Dm$. By Kellogg–Warschawski theorem we can map it to the unit disc, $Z_0: \overline{\D}\to \overline{\Dm}$ such that $Z_0\in C^{k,\alpha}(\bar\D)$. By Lemma \ref{lem:holo_W} and definitions, for
    $$X=  Z_0-\bar V,\quad Y =  Z_0+\bar V,$$
    we have $X:\bar\D\to \bar B_1(0)$ and $Y:\bar\D\to \bar\Lambda$ are $C^{k,\alpha}$ diffeomorphims for holomorphic $V = W\circ Z_0$. Using Taylor expansions for $Z_0,V$ and parametrization $e^{it}$ for $\partial\D$, the Fourier expansions of $x= X(e^{it}),y=Y(e^{it})$ are,
    $$\hat{x}(k) = \sgn(k)\hat{y}(k),\;\;k\neq 0.$$
    Here $\hat{x}(k),\hat{y}(k)$ are the $k$-th Fourier coefficient of $x,y$ respectively. In particular, it is immediate that $y = J\H x+c$ for some $c\in \R^2$ coming from the zeroth Fourier coefficient. Moreover since $X:\bar\D\to \bar B_1(0)$ was a $C^{k,\alpha}$ diffeomorphism till the boundary, there exists $\theta\in \mathcal{U}_{+}^{k,\alpha}$ such that $x=e^{i\theta}$. Similarly since $Y:\bar\D\to\bar\Lambda$ was a $C^{k,\alpha}$ diffeomorphism we get $y$ is a $C^{k,\alpha}$ parametrization of $\partial\Lambda$.\\
    
    Now assume we are given a pair $x,y$ as in \eqref{eq:traces_charac}. Let $X = P[x],Y = P[y]$. Since $x,y$ are related by the relation in \eqref{eq:traces_charac} and $P$ commutes with $\H$, we also get,
    $$X = J\H Y+c.$$
    Now for $z=\xi+i\eta$, Cauchy Reimann equations immediately give,
    $$X_{\xi} = -JY_{\eta},\;\; X_{\eta} = JY_{\xi}.$$
    In particular $DY = \text{cof}(DX)$, $\det(DX) = \det(DY) = d$ and $(DX)^TDY = dI$. Note that since $X(e^{it}) = x$ is an orientation preserving diffeomorphism of the unit circle, by Rado-Kneser-Choquet Theorem, $X$ is a harmonic diffeomorphism. By defining $T = Y\circ X^{-1}$, we get 
    $$DT = (DT)^{T} >0,\;\;\det(DT)=1.$$
    Moreover since $B_1(0)$ is convex, $T$ is injective. Therefore $Y$ is also injective and since $Y\in C^{k,\alpha}(\bar\D)$, we get $Y:\bar\D\to \bar\Lambda$ is a $C^{k,\alpha}$ diffeomorphism (since $DY>0$). Therefore $T$ is a non degenerate optimal transport map between $\bar\B_1$ and $\bar\Lambda$.

\end{proof}

\begin{rem}
    In the above Theorem the source domain can be taken to be any $C^{k,\alpha}$ convex source $\Omega$ (and also ensuring $|\Lambda| = |\Omega|$), where we replace $x = e^{i\theta}$ simply by $x = \gamma(e^{i\theta})$ for an regular parametrization $\gamma$ of $\Omega$. 
\end{rem}

Now given a target domain $\Lambda$, we would like to find the parameterization $\theta\in\mathcal{U}_{+}^{k,\alpha}$ that would guarantee (when possible) smooth non degenerate optimal transport. We connect it to a minimization problem using isoperimetric inequality. Note that given a pair of smooth parametrizations $x,y$ as in \eqref{eq:traces_charac} we note that,
$$x' = JHy' = JKy.$$
Moreover by the isopermetric inequality we get,
$$\left(\int_{\T}|Ky|dt\right)^2 = \left(\int_{\T}|x'|dt\right)^2 \geq 4\pi^2.$$
Therefore $\Lmin (\Lambda) \geq 2\pi$ with equality at the parameterizations given by \eqref{eq:traces_charac} since $x$ parametrizes the unit circle. Motivated by this observation we consider a variational problem. Fix a positively oriented  regular parametrization $\Gamma:\R/\mathbb Z\to\partial\Lambda$  (with regularity specified in each theorem) and let
$$
 \mathcal S=\{\sigma\in C^0(\R):\sigma\text{ nondecreasing},\ \sigma(2\pi \frac j3)=\frac j3\text{ for } (j=0,1,2,),\sigma(t+2\pi) = \sigma(t)+1\}.
$$

For $y=\Gamma\circ\sigma$ and $x=J\H y$, define
\begin{equation}\label{energy}
 \mathcal{O}_{\Lambda}(y)=\Var(x) = \int_{\mathbb{S}^1}|dx|,\quad \Lmin(\Lambda)=\inf_{\sigma\in\mathcal{S}}\mathcal{O}_{\Lambda}(\Gamma\circ\sigma),
\end{equation}
where $\mathcal{O}_{\Lambda}=+\infty$ if $dx$ is not a finite measure. Here $dx$ is the distributional derivative of $x$. Note that when $x,y$ are absolutely continuous, $x' = JKy$ and it reduces to minimizing $\int_\T |Ky|$. Note that the 3 point normalization is the class $\mathcal{S}$ can be arranged for any regular degree one parameterization $y$ by using automorphisms of the unit disc without changing $\Var(x)$. We now show that for finite $\mathcal{O}_{\Lambda}$, the parameterizations $x,y$ are indeed absolutely continuous.

\begin{lem}\label{lem:att_charac}
Let $y = \Gamma\circ \sigma$ for some $\sigma\in \mathcal{S}$. If $\mathcal{O}_{\Lambda}(y)$ is finite, then both $y$ and $x = J\H y$ are absolutely continuous. Moreover, the infimum in \eqref{energy} is attained.
\end{lem}
\begin{proof}
Since $y$ parametrizes $\partial\Lambda$, $\Var(y) = \H^{1}(\partial\Lambda)$. Moreover since $\mathcal{O}_{\Lambda}(y) = \Var(x)<\infty$, we have $x,y$ are BV functions. Since $x = J\H y$, we have $\hat{x}(k) = \sgn(k)\hat{y}(k)$ for the Fourier coefficients. Therefore the negative Fourier coefficients vanish for the measures $d(x+y)$ and $d(\overline{x-y})$. The Theorem of F. and M. Riesz (see \cite{Rudin}) gives $dx,dy\ll dt$. In particular
\begin{equation}\label{length}
 \mathcal{O}_{\Gamma}(y)=\int_\T |x'|\,dt=\int_\T|Ky|\,dt.
\end{equation}

Now for a minimizing sequence,
$$
 \Var(y_j)=\H^1(\partial\Lambda),\qquad \sup_j\Var(x_j)<\infty.
$$
Moreover since $\int_\T x_j = 0$ we have,
$$|x_j(t)| \leq \frac{1}{2\pi}\int_\T|x_j(t)-x_j(s)|ds\leq \Var(x_j) .$$
Therefore the BV norms of $x_j,y_j$ are bounded and hence up to a subsequence $x_j\to x$ and $y_j\to y$ in $L^1(\T)$ with the limits also satisfying $x = J\H y$. Moreover by compactness of BV functions, $x,y$ are also BV. Furthermore, Helly's selection theorem gives up to a subsequence that $\sigma_j\to \sigma$ for some monotone $\sigma$ and therefore $y = \Gamma \circ \sigma$. Since we also have imposed a three point normalization, injectivity of $\Gamma$ forces $\sigma\in \mathcal{S}$ due to continuity of $y$. By lower semicontinuity of total variation we get that $\mathcal{O}_{\Lambda}(y) = \Var(x) = \Lmin(\Lambda)$.
\end{proof}

\begin{lem}\label{lem:iso_charac}
    Let $y = \Gamma\circ \sigma$ for some $\sigma\in \mathcal{S}$ and $x = J\H y$. If $\mathcal{O}_{\Lambda}(y)$ is finite then,
    $$\A(x) = \frac{1}{2}\int_\T \det(x,x') = \frac{1}{2}\int_\T \det(y,y') = \pi.$$
    Moreover $\mathcal{O}_{\Lambda}(y)\geq 2\pi$ with equality if and only if $x=c+e^{i\theta}$ for continuous nondecreasing degree 1 lift $\theta$.
\end{lem}
\begin{proof}
    From the proof of Theorem \ref{thm:boundary_charac_unit_circle} we know if $X=P[X]$ and $Y=P[y]$ then $\det(DX) = \det(DY)$ in $\D$. Therefore by Green's theorem we have,
    \begin{align}
        \int_\T \det(X(re^{it}),X_t(re^{it}))&dt =2\int_{\D_r}\det DX \\
        &=2\int_{\D_r}\det DY\\
        &=\int_\T \det(Y(re^{it}),Y_t(re^{it}))dt
    \end{align}
    taking limits $r\uparrow 1$ we get
    $$\frac{1}{2}\int_\T \det(x,x') = \frac{1}{2}\int_\T \det(y,y') =|\Lambda|= \pi,$$
    since $y = \Gamma\circ\sigma$ is a positively oriented parametrization of $\partial\Lambda$. Now by using Banchoff-Pohl inequality and Green's theorem, see \cite[Corollary 2, Theorem 1]{Vogt},
    $$\pi = \frac{1}{2}\int_\T\det(x,x') = \int_{\R^2}w(p)dp \leq\int_{\R^2}(w(p))^2dp  \leq \frac{(\Var(x))^2}{4\pi}.$$
    Here $w(p)$ is the winding number of the curve $x(\T)$ around the point $p$. Equality happens if and only if the parametrization $x$ traverses the circle only once (since $w(p) \in \{0,1\}$ for equality). Therefore we get $x = c+e^{i\theta}$ for a continuous nondecreasing degree 1 lift $\theta$.
\end{proof}

We now prove the following continuous optimal transport characterization.
\begin{thm}\label{thm:infi_charac}
    For the infimum $\Lmin$,
$$
\Lmin(\Lambda)=2\pi\iff\Lambda\in\mathcal C.
$$
At equality the pair $x,y$ reconstructs the optimal transport map $T = P[y]\circ P[x]^{-1}$ for a centered parameterization $x$.
\end{thm}
\begin{proof}
    Combining Lemma \ref{lem:att_charac} and \ref{lem:iso_charac} we get that the infimum is attained and at the equality we have $x = c+e^{i\theta}$ and $y = \Gamma\circ \sigma$ for some $\sigma\in \mathcal{S}$ and a continuous nondecreasing degree 1 lift $\theta$. The if and only if statement now follows similar to the proof of Theorem \ref{thm:boundary_charac_unit_circle} by using the weak Rado-Kneser-Choquet theorem and the midpoint mixed domain being defined using the inverse optimal transport map ($S = T^{-1}$ that is continuous till the boundary by Lemma \ref{lem:external}) exactly as in Section \ref{subsec:LewyT}. Instead of using Kellog-Warschawski theorem we use Caratheodory's theorem.
\end{proof}

Combining the above theorem with Theorem \ref{thm:boundary_charac_unit_circle} we note that we have to find minimizers of
$$\mathcal{O}_{\Lambda}(\sigma) = \int_\T |K(\Gamma\circ \sigma)(t)|dt,\;\;\sigma\in \mathcal{S}.$$
Writing $y = \Gamma\circ \sigma$ we note that if at the minimizer $\Lmin(\Lambda) = \int_\T |Ky|dt = 2\pi$, then not only does the minimizer determines the optimal transport map but also the regularity of the minimizer determines the regularity of the optimal transport map (as long as their derivatives don't vanish). We therefore have the following version of Theorem \ref{}.

\begin{thm}
    Let $y\in C^{k,\alpha}(\T;\R^2)$, $k\ge2$, $0<\alpha<1$, be a positively oriented regular Jordan parametrization of $\partial\Lambda$, and suppose $Ky\ne0$ everywhere. Then
\begin{equation}\label{eq:equi_charac}
 y\text{ is stationary}\quad\Longleftrightarrow\quad
 y'\cdot K\left(\frac{Ky}{|Ky|}\right)=0
 \quad\Longleftrightarrow\quad \mathcal{O}_{\Lambda}(y)=2\pi.
\end{equation}
Stationarity is with respect to smooth inner reparametrizations. 
\end{thm}
\begin{proof}
    For $\eta$ smooth periodic consider $y_{\epsilon} = y\circ (Id+\epsilon\eta)$ for small $\epsilon$. Self-adjointness of $K$ gives
\begin{equation}\label{eq:var_charac}
 \left.\frac{d}{d\varepsilon}\mathcal{O}_{\Lambda}(y_\varepsilon)\right|_{0}
 =\int_\T \frac{Ky}{|Ky|}\cdot K(\eta y')\,dt
 =\int_\T \eta y'\cdot K\left(\frac{Ky}{|Ky|}\right)dt.
\end{equation}
Thus the first equivalence holds. $E(y)=2\pi$ implies stationarity (in fact minimality) by Theorem \ref{thm:infi_charac} and our assumption $Ky\neq 0$. \\

We now assume stationarity and prove minimality. Let $x = J\H y$ and denote $X = P[x]$. Let $p = \frac{x'}{|x'|} = \frac{JKy}{|Ky|}$ be the unit tangent vector to the curve $x(\T)$ and denote $U = P[p]$. We first show that $p$ is a circle diffeomorphism. Then by rigidity due to the stationarity equation and harmonic conjugacy it is shown that $x$ had to parameterize the unit circle.\\

We first note that for a function $f$ on $\T$, its harmonic extension is given by,
$$P[f] = \sum_{n\in\mathbb{Z}}r^{|n|}\hat{f}(n)e^{int}$$
in radial coordinates $(r,t)$. Therefore the radial derivative $\partial_r P[f] = r^{-1}P[Kf]$ since $Ke^{int} = \H \partial_t e^{int} = |n|e^{int} $. Similarly its tangential derivative on $\partial\D$ can be shown to be $\partial_t P[f] = f'$. All derivatives taken are at $\partial\D$ unless specified. Note that $(J\H)^2 = Id$ (for mean 0 functions) and therefore,
$$X_r=-Jy',\quad X_t=x'=|x'|p,\quad
 U_r=Kp = JK\left(\frac{Ky}{|Ky|}\right),\quad U_t=p'.$$
 Using $x=J\H y$ and stationarity we also get,
 $$-y'\cdot K\left(\frac{Ky}{|Ky|}\right)= X_r\cdot Kp = X_r\cdot U_r = 0.$$
Therefore $U_r$ is parallel to $y'$. Also note that $X_t\cdot U_t = |x'|p\cdot p'=0$ since $|p|=1$. For a relation for the mixed derivatives we first note that if $u$ is a harmonic function that $ru_r$ and $-u_t$ are harmonic conjugates of each other. Therefore $rX_r,rU_r$ are (componentwise) harmonic conjugates to $-X_t,-U_t$. Therefore, by multiplying the corresponding holomorphic functions componentwise and summing, we note that,
 \begin{align}
     a &= r^2X_r\cdot U_r-X_t\cdot U_t\\
     b&= -r(X_r\cdot U_t+X_t\cdot U_r)
 \end{align}
 are harmonic conjugates of one another. Since we know that $a=0$ on $\partial\D$ we get that $b$ has to be a constant. Since $b(0)=0$ we get the following mixed derivative equation,
 $$X_r\cdot U_t+X_t\cdot U_r = 0.$$
 Note that since $|p|=1$, the harmonic function $1-p(t_0)\cdot U$ is nonnegative with a minimum at $z = e^{it_0}$. Therefore by Hopf lemma, $p\cdot U_r>0$. Since $U_r$ was parallel to $y'$ we get that $p\cdot y'\neq 0$. Using the mixed derivative equation we get,
 $$|x'|p\cdot U_r-Jy'\cdot p'=0.$$
 Since $y$ was a regular parametrization, $Jy'\cdot p'>0$ gives $p'\neq 0$. Therefore $p(\T)$ is a parametrization of the unit circle since $|p|=1$ and $|p'|>0$. Since $y$ was a regular parametrization and we know $p\cdot y' \neq 0$, the curves $p$ and $y'/|y'|$ have to be homotopic to one another and in particular the degree of the curves should be the same. Since $y'/|y'|$ represents the unit tangent to $\partial\Lambda$, its winding number is 1 and therefore the winding number of $p$ is also 1. Therefore $p$ is a circle diffeomorphism and by Rado-Kneser-Choquet theorem $U$ is a harmonic diffeomorphism. Now note that $X_r\cdot U_r = X_t\cdot U_t = 0$ shows $X_r = aJU_r,X_t = bJU_t$ for some real valued functions $a,b$. Substituting in the mixed derivative identity we get,
 $$(a-b)\det(DU) = 0.$$
 Since $\det DU>0$ we get that $a=b$ and $DX = a JDU$. If we write $V=-JX$ we get that $DV=aDU$ for two harmonic maps $V,U$. Note that if we only look at the first components $V_1,U_1$ and define the holomorphic functions,
$$f  =V_{1x}-iV_{1y}, \quad g =U_{1x}-iU_{1y}, $$
we get that $f = ag$ and $g\neq 0$. Therefore $f/g = a$ is a holomorphic function with real boundary values. Hence it has to be constant and therefore $DX = aJDU$ for a constant $a$. Considering tangential components we get $x' = aJp'$ and therefore $x = aJp+c$. Since $p$ is a circle diffeomorphism, then $x$ also has to be one and therefore $\mathcal{O}_{\Lambda}(y) = \Var(x) = 2\pi$.
\end{proof}

\section{A dense subset of Target domains}\label{sec:poly}

We first define an admissible class of polynomial generators,

\begin{defn}
Let $P\in\C[z]$ and write $P(x+iq)=A(x,q)+iB(x,q)$. A polynomial is called admissible over a source domain $\Omega$ if there exists a bounded $C^{2,\alpha}$ Jordan domain $D_P$ such that
$$
 X_P(x,q)=(x,B(x,q)):\bar D_P\longrightarrow\bar \Omega
$$
is a diffeomorphism and $ A'>0$ on $\bar D_P$. 
\end{defn}
Associated to an admissible polynomial $P$, we denote its target domain as
$$\Lambda_P = Y_P(D_P),$$ 
with $Y_P(x,q) = (A(x,q),q)$ and the optimal transport map 
$$T_P = Y_P\circ X_p^{-1}: \Omega\to \Lambda_P.$$
We are concerned with the following theorem in this section,
\begin{thm}\label{thm:density}
Let $0<\beta<\alpha<1$, and let $\Omega$ be a $C^{2,\alpha}$ bounded convex domain. Suppose there exists a optimal transport map
$$T: \Omega\to \Lambda,$$
such that $T\in C^{1,\alpha}(\bar\Omega)$. Then for every $\epsilon>0$ there exists an admissible polynomial $P$ over $\Omega$ such that
$$\|T_P-T\|_{C^{1,\beta}(\bar\Omega)} \leq \epsilon.$$
Consequently the boundaries of the polynomial targets $\Lambda_P$ are dense in the class of $C^{1,\alpha}$ targets admitting nondegenerate optimal transport, in every $C^{1,\beta}$ topology of boundary embeddings.  
\end{thm}

\begin{rem}
    For domains that are $C^{k,\alpha}$ smooth, a similar density statement can be made under stronger norms. In particular for $C^{\infty}$ target domains $\Lambda$ and source domain $B_1(0)$, Corollary \ref{cor:intro_polynomial} follows immediately.
\end{rem}

The proof for the above theorem also gives us the following interior version of the statement,
\begin{thm}\label{thm:density_interior}
    Let $\Omega$ be a $C^{2,\alpha}$ bounded convex domain and $\Lambda$ be a bounded domain. Suppose there exists a optimal transport map
$$T: \Omega\to \Lambda,$$
such that $T$ is continuous on $\Omega$. Then there exists a sequence of admissible polynomials over $\Omega$, $P_j$, such that their corresponding optimal transport maps
$$T_{P_j}:\Omega\to \Lambda_{P_j}$$
converge to $T$ in $C^{\infty}_{loc}$. Moreover $\Lambda_{P_j}$ converge to $\Lambda$ in Hausdorff distance and measure.
\end{thm}
We fix some notation. Choose an interior point $x_*\in\Omega$ and, for $r>0$, set
$$
 \Omega_r=x_*+r(\Omega-x_*).
$$
Convexity and $x_*\in\Omega$ imply that
\begin{equation}\label{eq:nested}
 1\le r<s\quad\Longrightarrow\quad
 \bar\Omega_r\Subset\Omega_s.
\end{equation}
Here the inclusion is strict whenever $r<s$. Consider the rescaling,
\begin{equation}\label{eq:Tt}
 T_t(x):=t^{-1}T\bigl(x_*+t(x-x_*)\bigr)
\end{equation}
for $x\in\Omega_{1/t}$. In particular $DT_t>0$ and $\det DT_t=1$ on the enlarged domain. We have the following convergence of maps $T_t$.
\begin{lem}\label{lem:scaling}
For every $0\le\beta<\alpha$,
\begin{equation}\label{eq:scalingconv}
 \|T_t-T\|_{C^{1,\beta}(\bar\Omega)}\longrightarrow0
 \qquad\text{as }t\uparrow1.
\end{equation}
\end{lem}
\begin{proof}
The formulas \eqref{eq:Tt} and the uniform continuity of $T,DT$ on $\bar\Omega$ show convergence in $C^1(\bar\Omega)$. Further,
$$
 [DT_t]_{C^\alpha(\bar\Omega)}
 =t^\alpha[DT]_{C^\alpha(x_*+t(\bar\Omega-x_*))}
 \le [DT]_{C^\alpha(\bar\Omega)}.
$$
Thus $DT_t-DT$ is uniformly bounded in $C^\alpha$ and tends to zero uniformly. The interpolation inequality
$$
 [f]_{C^\beta}\le C\|f\|_{C^0}^{1-\beta/\alpha}[f]_{C^\alpha}^{\beta/\alpha},
 \qquad 0<\beta<\alpha,
$$
shows convergence in $C^\beta$ of the derivatives.
\end{proof}

For each $T_t$ we write its holomorphic generator as $H_t$ with mixed domain written as $D_t$. We also denote the map from the mixed domain to the source domain as $X_t$. Fix $t<1$ and choose
\begin{equation}\label{eq:radii}
 1<r_0<r_1<r_2<1/t.
\end{equation}

For $j=0,1,2$ define
$$
 K_j=X_t^{-1}(\bar\Omega_{r_j}).
$$
The mixed map is a homeomorphism on the corresponding closures and $\bar\Omega_{r_2}\Subset\Omega_{1/t}$, so
\begin{equation}\label{eq:collar}
K_0\Subset\operatorname{int}K_1\Subset\operatorname{int}K_2\Subset D_t.
\end{equation}
Each $K_j$ is a closed Jordan disk, because it is the image of the closed convex Jordan domain $\bar\Omega_{r_j}$ under a homeomorphism. 
By  Mergelyan's theorem applied to $K_2$, there exist polynomials $P_j$ with
\begin{equation}\label{eq:merg}
 \|P_j-H_t\|_{C^0(K_2)}\longrightarrow0.
\end{equation}
Since $K_1\Subset\operatorname{int}K_2$, Cauchy's estimates on small complex disks uniformly contained in $\operatorname{int}K_2$ imply for every finite integer $m\ge0$ that
\begin{equation}\label{eq:cauchy}
 \|P_j-H_t\|_{C^m(K_1)}\longrightarrow0.
\end{equation}
We show explicitly why the source remains exactly $\Omega$.

\begin{lem}\label{lem:branch}
For all sufficiently large $j$, the map
$$
 X_{P_j}(x,q)=(x,\Ima P_j(x+iq))
$$
is injective on $K_1$, $\operatorname{Re}P_j'>0$ on $K_1$, and
\begin{equation}\label{eq:covered}
 \bar\Omega_{r_0}\subset\operatorname{int}X_{P_j}(K_1).
\end{equation}
Consequently, setting
\begin{equation}\label{eq:Dpj}
 D_{P_j}:=X_{P_j}^{-1}(\Omega)\cap\operatorname{int}K_1,
\end{equation}
we obtain an admissible branch over $\Omega$.
\end{lem}
\begin{proof}
The proof follows from the fact that since $X_t$ was a diffeomorphism on $K_2$, a small $C^{1}$ perturbation remains a diffeomorphism. Moreover for a large $j$ one can ensure the boundaries $X_t(K_1)$ and $X_{P_j}(K_1)$ remain close so that \eqref{eq:covered} follows. The admissibility follows from the fact that $K_1\Subset D_t $ therefore by compactness $\Rea H_t'\geq c_t>0$ for small $c_t$ and \eqref{eq:cauchy}.
\end{proof}

As an immediate corollary of the inverse function theorem and \eqref{eq:cauchy} we get,
\begin{cor}\label{lem:invconv}
For fixed $t<1$ and $r_0,r_1,r_2$ as in \eqref{eq:radii}, the polynomial models constructed above satisfy
\begin{equation}\label{eq:Tjconv}
 T_{P_j}\longrightarrow T_t
 \quad\text{in }C^m(\bar\Omega)
 \quad\text{for every finite }m\ge0.
\end{equation}
\end{cor}

The proof of Theorems \ref{thm:density} and \ref{thm:density_interior} follows from Lemma \ref{lem:scaling} and Corollary \ref{lem:invconv}.\\

Now we consider the examples \ref{ex:higher_order_poly}. In particular we show that these examples could never be obtained as $C^{1,1}$ or $C^2$ perturbations of nearby convex objects despite being close in Hausdorff sense.
\begin{prop}
    Fix $c>0$ and $\vartheta\in\R$. For all sufficiently large integers $n$,
$$
H_n(z)=z+\frac{c}{n^2}e^{i\vartheta}z^n
$$
induces an optimal transport map between $B_1(0)$ and a smooth target $\Lambda_n$ whose signed boundary curvature satisfies, uniformly in $t$,
$$
\kappa_n(t)=1+c\cos((n+1)t+\vartheta)+O_c(n^{-1}).
$$

\end{prop}
\begin{rem}
    In particular, for $c>1$ and large $n$, the target is nonconvex. Taking $c$ arbitrarily large produces targets with arbitrarily negative curvature while preserving smooth transport. Therefore its impossible to obtain the targets by small $C^{1,1}$ perturbations of convex domains.
\end{rem}
\begin{proof}
    On $\partial B_1(0)$ set
$z_n(t)=\cos t+i q(\cos t,\sin t)$ and $\gamma_n(t)=T_n(\cos t,\sin t)$ where $T_n:B_1(0)\to \Lambda_n$ is the optimal transport map. Here $q$ is defined as in Example \ref{ex:higher_order_poly}. From the discussion in example \ref{ex:higher_order_poly},
$$
z_n(t)=e^{it}+O_c(n^{-2}).
$$
Note that by the implicit relation
\begin{equation}\label{eq:ex_implicit}
    (\Rea z_n(t),\Ima H_n(z_n(t))) = e^{it},
\end{equation}
we also obtain the following information on the derivative $z_n$,
\begin{equation}\label{eq:zprime}
    z_n'(t) = ie^{it}+O_c(n^{-1}).
\end{equation}
Therefore 
$$\gamma_n(t) = (\Rea H_n(z_n(t)),\Ima z_n(t)),\qquad \gamma_n'(t)=ie^{it}+O_c(n^{-1}).$$
Note that by smoothness of $H_n$ we know that $z_n$ is smooth. Calculating the signed curvature of the curve $\gamma_n$ gives,
$$\det(\gamma_n',\gamma_n'') =1-\Ima\big(H_n''(z_n)(z_n')^3\big).$$
The information on $q''$ follows from differentiating the implicit relation \eqref{eq:ex_implicit}.
\end{proof}
\begin{figure}[h]
    \centering
    \includegraphics[width=0.75\linewidth]{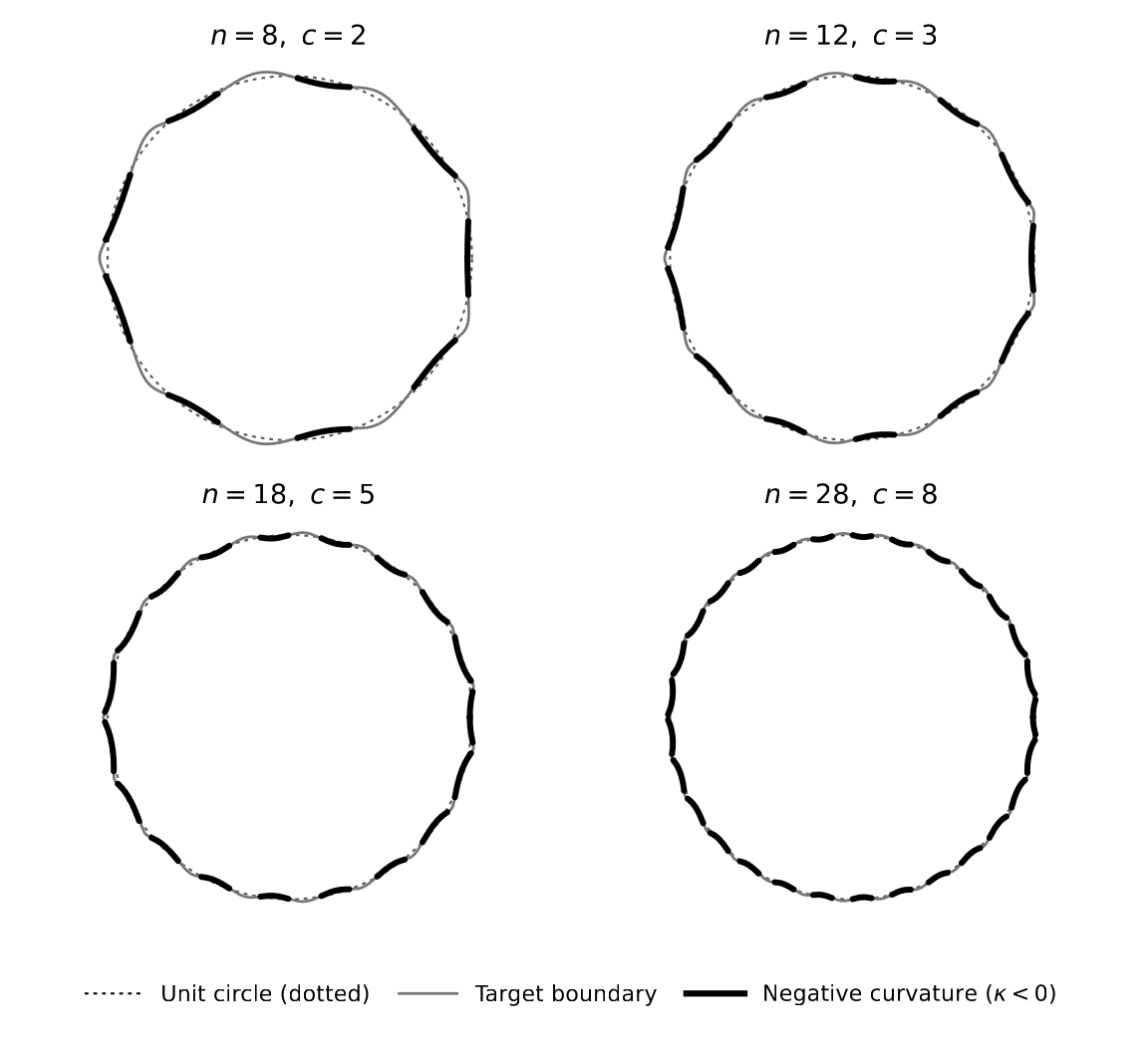}
    \caption{Different target domains obtained by large $n$ and $c$.}
\end{figure}

\section{Real analytic singular sets}\label{sec:free_boundary}
In this section we deal with singular sets that can be represented as real analytic graphs. We first return to the Laurent example from Section~\ref{subsec:ellipse} to study the singular set of the inverse transport and the behavior at its endpoints. Related connections with the thin obstacle problem were discussed in \cite[Remark~3.6]{RS1}; see also \cite{McCannRakshit} for the phantom transport viewpoint developed by McCann and one of us. We write explicit formulas to determine the jump of the inverse map, its singular Monge--Amp\`ere measure, and the sharp boundary regularity of the forward map. We then extend the construction from a straight segment to a prescribed analytic graph.

\subsection{Further analysis of the elliptical example}
Let the parameters $\lambda,\mu,r_0,R$, the generator $f$, mixed domain $D$, and the coordinate assignments $X(z),Y(z)$ be as in the Laurent generator example. We use the following notations
\begin{equation}\label{eq:ellipsedomains}
 \begin{aligned}
 K_{\rm hole}&=K=E(r_0,2\lambda r_0),\\
 K_{\rm out}&=E(R,\lambda R+\mu/R),\\
 L_{\rm out}&=E(\lambda R-\mu/R,R),\\
 S_0&=\{0\}\times[-r_0,r_0].
 \end{aligned}
\end{equation}
Thus $\Omega=K_{\rm out}\setminus\overline{K_{\rm hole}}$ and $\Lambda=L_{\rm out}\setminus S_0$. The optimal transport map transports the uniform density on $\Omega$ to that on $L_{\rm out}$. We retain the globally defined convex extension $\overline u$ from the example:
$$
 \overline u=u\quad\text{outside }\overline{K_{\rm hole}},
 \qquad
 \overline u(x,y)=\frac{y^2}{4\lambda}
 \quad\text{on }\overline{K_{\rm hole}}.
$$
\subsubsection{The collapsing boundary and the inverse transport}
The boundary values in \eqref{eq:ellipse-points} give
$$
 X(r_0e^{i\theta})
 =(r_0\cos\theta,2\lambda r_0\sin\theta),
 \qquad
 Y(r_0e^{i\theta})=(0,r_0\sin\theta).
$$
By \eqref{eq:ellipse-potential} and the formulas above, the values
and gradients match continuously across $\partial K_{\rm hole}$.
Hence $\overline u\in C^1(\R^2)$. In particular,
\begin{equation}\label{eq:ellipsechords}
 \nabla\overline u(x,2\lambda q)=(0,q),
 \qquad
 |q|\leq r_0,\quad |x|\leq\sqrt{r_0^2-q^2}.
\end{equation}
Thus every horizontal chord of the filled cavity maps to one point of $S_0$. The two cavity tips $(0,\pm2\lambda r_0)$ map to $(0,\pm r_0)$.\\
Denote the full Legendre conjugate of $\overline u$ by
$$
 \bar{v}(p,q)
 =\sup_{(x,y)\in\R^2}
 \{xp+yq-\overline u(x,y)\}.
$$
The quadratic growth of $\overline u$ gives $\bar{v}$ finite everywhere. Away from $S_0$, it is  the potential of the is the optimal transport. Using \eqref{eq:ellipse-potential}, we obtain
\[
 \nabla \bar{v}(Y(z))=X(z),
 \]
 \[
 \bar{v}(Y(re^{i\theta}))
 =\frac{\lambda}{2}\bigl(r^2-r_0^2\cos2\theta\bigr)
   +\mu\log(r/r_0),
   \]
 for $r>r_0$. Consequently, $\nabla \bar{v}$ gives the reverse optimal transport from $L_{\rm out}$ to $\Omega$.\\
Using convex conjugacy and \eqref{eq:ellipsechords}, we have
\begin{equation}\label{eq:ellipsesubgradient}
 \begin{aligned}
 \bar{v}(0,q)&=\lambda q^2,\\
 \partial \bar{v}(0,q)
 &=[-\sqrt{r_0^2-q^2},\sqrt{r_0^2-q^2}]
   \times\{2\lambda q\},
 \end{aligned}
\end{equation}
for \(|q|\leq r_0.\)
Thus each subgradient along the segment is exactly one horizontal chord of the cavity. For $|q|<r_0$, the two boundary values of the inverse transport are
\begin{equation}\label{eq:ellipsetraces}
 \lim_{p\to0^\pm}\nabla \bar{v}(p,q)
 =\left(\pm\sqrt{r_0^2-q^2},2\lambda q\right).
\end{equation}
The jump is perpendicular to $S_0$ and has magnitude $2\sqrt{r_0^2-q^2}$. The nondifferentiability set of $\bar{v}$ is $\{0\}\times(-r_0,r_0)$, whose closure is $S_0$. At the endpoints, the subdifferentials are singletons.

\subsubsection{The singular mass}
For a convex function $v$, the Monge–Ampère measure $\operatorname{MA}_v(B)=|\partial v(B)|$. Formula \eqref{eq:ellipsesubgradient} gives us
\begin{equation}\label{eq:ellipse-fb-mass}
 \operatorname{MA}_{\bar{v}}
 =\mathcal L^2
 +4\lambda\sqrt{r_0^2-q^2}\,
   \mathcal H^1\!\restriction_{S_0}.
\end{equation}
Indeed, $\det D^2\bar{v}=1$ away from $S_0$. For a Borel set $A\subset(-r_0,r_0)$, the subgradients over $\{0\}\times A$ sweep out horizontal chords of length $2\sqrt{r_0^2-q^2}$ at height $y=2\lambda q$. Fubini's theorem therefore gives
$$
 \left|\partial \bar{v}(\{0\}\times A)\right|
 =4\lambda\int_A\sqrt{r_0^2-q^2}\,dq.
$$
There are no endpoint atomic measures because the corresponding subdifferentials are singletons. In particular,
\begin{equation}\label{eq:ellipse-fb-total-mass}
 \operatorname{MA}_{\bar{v}}(S_0)
 =4\lambda\int_{-r_0}^{r_0}\sqrt{r_0^2-q^2}\,dq
 =2\pi\mu
 =|K_{\rm hole}|.
\end{equation}
Thus the singular mass of the inverse potential is exactly the area of the cavity.
\begin{rem}\label{rem:ellipsephantom}
The computation above gives an explicit formula for the phantom target measure of \cite{McCannRakshit}, associated with the reverse transport $\nabla\bar v:L_{\rm out}\to\Omega$:
$$
 \mu_{\rm ph}
 =(\nabla\overline u)_\#\bigl(\mathbf 1_{K_{\rm hole}}\,dx\,dy\bigr)
 =4\lambda\sqrt{r_0^2-q^2}\,
 \mathcal H^1\!\restriction_{S_0}.
$$
The corresponding phantom map is $\nabla\overline u(x,y)=(0,y/(2\lambda))$ on $K_{\rm hole}$. By \eqref{eq:ellipse-fb-mass}, $\mu_{\rm ph}$ is precisely the excess Monge--Amp\`ere measure of $\bar v$ inside $L_{\rm out}$. Its total mass is \( |K_{\rm hole}| \), the support is the closed segment \(S_0\). Note that the endpoints carry no atoms.
\end{rem}

\subsubsection{The cavity as a contact set: an obstacle-type interpretation}
The quadratic filling also gives a contact interpretation in the source plane. Set
\(
 \psi(x,y)=\frac{y^2}{4\lambda}.
\)
For $r>r_0$, equation \eqref{eq:ellipse-potential} gives
$$
 \begin{aligned}
 \overline u(X(re^{i\theta}))-\psi(X(re^{i\theta}))
 &=\frac{\lambda}{2}(r^2-r_0^2)-\mu\log(r/r_0)
   -\frac{\lambda(r^2-r_0^2)^2}{4r^2}\sin^2\theta\\
 &\geq
 \frac{\lambda}{4}\left(r^2-\frac{r_0^4}{r^2}\right)
 -\mu\log(r/r_0)>0.
 \end{aligned}
$$
The last expression vanishes when $r=r_0$ and has derivative $\lambda(r^2-r_0^2)^2/(2r^3)>0$. Hence,
\begin{equation}\label{eq:ellipse-fb-contact}
 \overline u\geq\psi,
 \qquad
 \{\overline u=\psi\}=\overline{K_{\rm hole}},
 \qquad
 \operatorname{MA}_{\overline u}
 =\mathbf 1_{\R^2\setminus\overline{K_{\rm hole}}}\,dx\,dy.
\end{equation}
For the measure identity, the gradient maps the cavity and its boundary into the null segment $S_0$, while its Jacobian is one outside the cavity. Thus $\partial K_{\rm hole}$ is the boundary of the contact set for the forward potential, whereas $S_0$ is the closure of the discontinuity set of the inverse transport.

\subsubsection{Sharp regularity at the cavity tips}
Recall that the generator is $f(z)=\lambda z-\mu/z$ on the mixed domain $D=\{z\in\C:r_0<|z|<R\}$, where $r_0=\sqrt{\mu/\lambda}$. Although $f$ extends holomorphically across the inner circle
$|z|=r_0$, the real part of its derivative vanishes at two points:
\begin{equation}\label{eq:ellipse-fb-jacobian}
 a(re^{i\theta})
 :=\Rea f'(re^{i\theta})
 =\lambda+\frac{\mu}{r^2}\cos2\theta,
 \qquad
 a(r_0e^{i\theta})=2\lambda\cos^2\theta.
\end{equation}
The exceptional mixed points are $z=\pm ir_0$, corresponding to the cavity tips $(0,\pm2\lambda r_0)$. Away from these tips, the forward map is real analytic up to the source boundary from the annular side.

\begin{prop}\label{prop:ellipseregularity}
The potential $\overline u$ satisfies
$$
 \overline u\in C^{1,1/2}(\overline{K_{\rm out}}),
 \qquad
 \overline u\notin C^{1,\beta}(\overline{K_{\rm out}})
 \quad\text{for }\frac12<\beta\leq1.
$$
Moreover, for $1\leq\ell<\frac{5}{2}$, we have \(
 u\in W^{2,\ell}(\Omega).\) Moreover the Sobolev exponent is sharp.
\end{prop}

\begin{proof}
In this proof, $\simeq$ denotes two-sided bounds by positive constants depending only on the fixed parameters.
We first prove the Sobolev regularity. Since \(f\) is holomorphic on a neighborhood of \(\overline D\), we have
\(
\det DX=\Rea f'>0
\)
on \(\overline D\setminus\{\pm ir_0\}\). Thus using inverse function theorem we have \(T=Y\circ X^{-1}\) is real analytic in \(\Omega\). It is also real analytic up to the boundary away from the points \((0,\pm2\lambda r_0)\). Thus $D^2u=DT$ is bounded outside any fixed neighborhoods of these two tips. It remains to analyze the behavior near them.

Near the two cavity tips, where \(z=\pm i r_0\), write $r=r_0+\delta$ and $\theta=\theta_0+\vartheta$, where $\theta_0=\pi/2$ or $3\pi/2$. Equation \eqref{eq:ellipse-fb-jacobian} above gives
$$
 a(re^{i\theta})
 =\lambda\left(1-\frac{r_0^2}{r^2}\right)
  +\frac{2\lambda r_0^2}{r^2}\sin^2\vartheta
 \simeq\delta+\vartheta^2.
$$
Since $f'$ is bounded near $\overline D$, equation \eqref{eq:Hessian} gives
$$
 |D^2u(X(re^{i\theta}))|
 \simeq\frac1{a(re^{i\theta})}
 \simeq(\delta+\vartheta^2)^{-1}.
$$
Moreover, the Jacobian of the source assignment is $a$, so
$$
 dx\,dy=a(re^{i\theta})\,r\,dr\,d\theta
 \simeq(\delta+\vartheta^2)\,d\delta\,d\vartheta.
$$
Hence, to analyze the local integrability of $|D^2u|^\ell$ at either tip is,  we need  finiteness of
$$
 \int_0^{\varepsilon^2}\int_{-\varepsilon}^{\varepsilon}
 (\delta+\vartheta^2)^{1-\ell}\,d\vartheta\,d\delta.
$$
Splitting the integration into the regions $\delta\geq\vartheta^2$ and $\delta<\vartheta^2$, the two contributions are comparable to
\(
 \int_0^{\varepsilon^2}\delta^{3/2-\ell}\,d\delta
\) and 
\(
 \int_0^\varepsilon\vartheta^{4-2\ell}\,d\vartheta.
\)
Both are finite exactly when $\ell<5/2$, and both diverge logarithmically at $\ell=5/2$. This gives us the desired Sobolev regularity. 

We next prove the H\"older continuity of the gradient. The $C^{1,1/2}$ estimate follows from \cite[Theorem~1.1]{MR2}. Inside the cavity, $\nabla\overline u=(0,y/(2\lambda))$ is Lipschitz and agrees with the exterior trace. The estimate therefore extends across the cavity boundary, giving $\overline u\in C^{1,1/2}(\overline{K_{\rm out}})$.

It remains to prove that the exponent is sharp. Along the upper vertical axis,
$$
 T\left(0,\lambda r+\frac{\mu}{r}\right)=(0,r),
 \qquad
 \nabla\overline u(0,2\lambda r_0)=(0,r_0).
$$
Since $\mu=\lambda r_0^2$, the source displacement satisfies
$$
 \lambda r+\frac{\mu}{r}-2\lambda r_0
 =\frac{\lambda(r-r_0)^2}{r}.
$$
Thus a H\"older estimate for $\nabla\overline u$ with exponent $\beta>1/2$ would imply
$$
 r-r_0
 \leq C\left(\frac{\lambda(r-r_0)^2}{r}\right)^\beta
 \leq C'(r-r_0)^{2\beta},
$$
which is impossible as $r\downarrow r_0$. Hence $\overline u\notin C^{1,\beta}(\overline{K_{\rm out}})$ for every $\frac12<\beta\leq 1$.
\end{proof}
\begin{rem}
The estimate $\overline u\in C^{1,1/2}(\overline{K_{\rm out}})$ can also be obtained directly from the explicit coordinate formulas \eqref{eq:ellipse-points}, without appealing to \cite{MR2}.
\end{rem}

\subsection{Singular sets as analytic graph}
 We first introduce some notations. For some $R>1$, let
\begin{equation}\label{eq:annulus}
 \A_R=\{z=x+iq\in\C:1<|z|<R\},
\end{equation}
and let $f=P+iQ$ be holomorphic on a neighborhood of $\overline{\A_R}$. We denote a real analytic function on a neighbourhood of $[-1,1]$ by $\Psi$. The desired singular set is denoted by
\begin{equation}\label{eq:slit}
 S_\Psi=\{(U,V)=(\Psi(s),s):-1\le s\le1\}.
\end{equation}

Since $\Psi$ is real analytic on a real neighborhood of $[-1,1]$, it admits a holomorphic extension to a complex neighborhood of this interval that is real-valued on the real axis. Set
\begin{equation}\label{eq:sw}
 s(z)=\frac{z-z^{-1}}{2i},\qquad w(z)=\frac{z+z^{-1}}2,
 \qquad s(z)^2+w(z)^2=1.
\end{equation}
For $z=e^{it}$ these functions are real, with
\begin{equation}\label{eq:swcircle}
 s(e^{it})=\sin t,\qquad w(e^{it})=\cos t.
\end{equation}
Define
\begin{equation}\label{eq:H}
 f_{\kappa,\Psi}(z)=i\kappa s(z)+\Psi(s(z))
             -i w(z)\Psi'(s(z)),\qquad \kappa>0.
\end{equation}
For $R>1$ sufficiently close to one, $s(\overline{\A_R})$ lies in the complex domain of analyticity of $\Psi$; hence $f_{\kappa,\Psi}$ is holomorphic on a neighborhood of the closed annulus. The parameter $\kappa$ will then be chosen sufficiently large.

\begin{thm}\label{thm:general_fbp}
For every real-analytic $\Psi$ on a neighborhood of $[-1,1]$, one can
choose $R>1$ and $\kappa>0$ so that $f_{\kappa,\Psi}$ is a holomorphic generator for some map $T:\Omega\to \Lambda$ that has the following properties:
\begin{enumerate}
\item The physical source is
\begin{equation}\label{eq:sourcegeometry}
 \Omega=K_{\rm out}\setminus\overline{K_{\rm hole}},
\end{equation}
where both $K_{\rm out}$ and $K_{\rm hole}$ are bounded uniformly convex analytic Jordan domains, $\overline{K_{\rm hole}}\subset K_{\rm out}$.
\item The physical target is
\begin{equation}\label{eq:targetgeometry}
 \Lambda=L_{\rm out}\setminus S_\Psi,
\end{equation}
where $L_{\rm out}$ is an uniformly convex analytic Jordan domain and
$S_\Psi\subset L_{\rm out}$.
The inner source $\partial K_{\text{hole}}$ maps two-to-one onto the target slit, except at the two slit endpoints.
\item The map $T: \Omega\to\Lambda$ is the optimal transport map between the two sets, which have equal area.
The exact areas are
\begin{equation}\label{eq:areas}
 \begin{aligned}
 |K_{\rm hole}|&=\pi\kappa,\\
 |K_{\rm out}|&=\frac{\pi\kappa}{2}(R^2+1),\\
 |\Omega|=|L_{\rm out}|=|\Lambda|
 &=\frac{\pi\kappa}{2}(R^2-1).
 \end{aligned}
\end{equation}
These areas are independent of the shape of $\Psi$.
\end{enumerate}
\end{thm}

\begin{rem}
    Theorem \ref{thm:intro_thm4} follows immediately by considering the inverse of the optimal transport map constructed above.
\end{rem}

 We note that the choice of $f_{\kappa,\Psi}$ is made so that the Jacobians of the mapping $X,Y$ associated to $f_{\kappa,\Psi}$ are positive.

 \begin{lem}\label{lem:psi_global_diff}
     Let $\Psi$ be a real analytic function on a neighbourhood of $[-1,1]$. Then there exists $\kappa$ large such that the mappings $X(x,q) = (x,\Ima f_{\kappa,\Psi})$ and $Y(x,q) = (\Rea f_{\kappa,\Psi},q)$ are global orientation preserving diffeomorphisms onto their images.
 \end{lem}
 
\begin{proof} Let $a,b$ such that
\begin{equation}\label{eq:ab}
 a(z)=\Rea f_{\kappa,\Psi}'(z),\qquad b(z)=\Ima f_{\kappa,\Psi}'(z).
\end{equation}
Note  that on $|z| = e^{it}$,
\begin{equation}\label{eq:XY_sec8}
    X(e^{it}) = (\cos(t),\kappa \sin(t)-\cos(t)\Psi'(\sin(t))),\;\; Y(e^{it}) = (\Psi(\sin(t)),\sin(t)). 
\end{equation}
Using $s'(z)=(1+z^{-2})/(2i)$ and $w'(z)=(1-z^{-2})/2$,
one obtains the exact identity
\begin{equation}\label{eq:Hp}
f_{\kappa,\Psi}'(z)=A_0(z)\bigl(\kappa-w(z)\Psi''(s(z))\bigr)
                     -i\Psi'(s(z)).
\end{equation}
In polar coordinates $z=re^{it}$, $r\ge1$,
\begin{equation}\label{eq:A0}
 \Rea A_0(re^{it})=\frac{1-r^{-2}}2+r^{-2}\cos^2t
 \ge\frac{(r-1)+\cos^2t}{R^2}
 \quad(1\le r\le R).
\end{equation}
On $|z|=1$, direct substitution in \eqref{eq:Hp} gives
\begin{equation}\label{eq:ainner}
a(e^{it})=\cos^2t\bigl(\kappa-\cos t\,\Psi''(\sin t)\bigr).
\end{equation}
In particular, if $\kappa>\|\Psi''\|_{L^\infty[-1,1]}$, the Jacobian is positive except at the two forced points $z=\pm i$. These are the points that map to the slit endpoints (their degeneracy
is expected as the target inner trace is double-covered and is not a smooth immersion at those two parameters). Moreover by \eqref{eq:A0} for a large 
\begin{equation}
    \kappa> C(\|\Psi''\|_{L^\infty[-1,1]}+\|\Psi'\|_{L^\infty[-1,1]})
\end{equation} 
we can also ensure positivity of $a$ on $\A_R$.\\
To ensure global injectivity of $X$  we note that it is enough to check injectivity on each line $x=x_0$. For $f_{\kappa,\Psi} = P+iQ$, $a = P_x=Q_y$. Therefore on each connected component of $\{(x_0,y)\;|\;y\in\R\}\cap \A_R$, $Q$ is injective. Note  that on $|z| = e^{it}$,
\begin{equation}\label{eq:XY_sec8}
    X(e^{it}) = (\cos(t),\kappa \sin(t)-\cos(t)\Psi'(\sin(t))),\;\; Y(e^{it}) = (\Psi(\sin(t)),\sin(t)). 
\end{equation}
Therefore for $z_{\pm}(x_0)=x_0\pm i\sqrt{1-x_0^2}$ we have,
\begin{equation}
    \begin{aligned}
        Q(z_+)-Q(z_-) &= 2\kappa \sqrt{1-x_0^2}-x_0(\Psi'(\sqrt{1-x_0^2})-\Psi'(-\sqrt{1-x_0^2}))\\
        &\geq 2\sqrt{1-x_0^2}(\kappa-\|\Psi''\|_{L^\infty[-1,1]})\geq 0.
    \end{aligned}
\end{equation}
Therefore $Q$ is globally injective on $\{(x_0,y)\;|\;y\in\R\}\cap \A_R$. The global injectivity of $Y$ along lines $q=q_0$ follow similarly.
\end{proof}

\subsubsection{The geometry of domains $\Omega$ and $\Lambda$.} Before proving the holomorphic generator indeed produces an optimal trasport map we first note the geometry of the domains $\Omega$ and $\Lambda$. Put $v=\sin t$ and $w=\cos t$. Since
$f_{\kappa,\Psi}(e^{it})=\Psi(v)+i(\kappa v-w\Psi'(v))$, the physical traces are
\begin{equation}\label{eq:innertraces}
 \begin{aligned}
 \gamma_{\rm hole}(t)&=X(e^{it})=
       \bigl(w,\ \kappa v-w\Psi'(v)\bigr),\\
 \gamma_{\rm slit}(t)&=Y(e^{it})=
       \bigl(\Psi(v),\ v\bigr).
 \end{aligned}
\end{equation}
Here $\gamma_{\rm hole}$ represents the boundary of $K_{\rm hole}$ and $\gamma_{\rm slit}$ represents the graph of $\Psi$ as the singular inside the target domain $L_{\rm out}$. We now prove the following structure lemma for the domains $\Omega$ and $\Lambda$.

\begin{lem}\label{lem:dom_geo_fbp}
    Let $\Omega = X(\A_R)$ and $\Lambda = Y(\A_R)$ for the mappings $X,Y$ associated to $H_{\kappa,\Psi}$ for $\kappa$ large. Then $\Omega = K_{\rm out}\setminus \overline{K_{\rm hole}}$ for two uniformly convex domains $K_{\rm out}$ and $K_{\rm hole}$. Also $\Lambda = L_{\rm out}\setminus S_{\Psi}$ for uniformly convex domain $L_{\rm out}$ and $S_{\Psi}$ as in \eqref{eq:slit}. The areas of $L_{\rm out}, K_{\rm out}$ and $K_{\rm hole}$ are as in point (3) of Theorem \ref{thm:general_fbp}.
\end{lem}
\begin{proof}
    Write $g(t)=-\cos t\,\Psi'(\sin t)$, so that $\gamma_{\rm hole}(t)=(\cos t,\kappa\sin t+g(t))$.
Its signed curvature numerator is
\begin{equation}\label{eq:holecurv}
 \begin{aligned}
 K_{\rm hole}(t)
 &:=\det(\gamma_{\rm hole}'(t),\gamma_{\rm hole}''(t))\\
 &=\kappa+\cos t\,g'(t)-\sin t\,g''(t)\\
 &=\kappa-\cos t(1+2\sin^2t)\Psi''(\sin t)
    +\sin t\cos^3t\Psi'''(\sin t).
 \end{aligned}
\end{equation}
In particular, if $\kappa>\|g'\|_\infty+\|g''\|_\infty$ then $K_{\rm hole}$ is a uniformly convex set and the embedding $\gamma_{\rm hole}$ is regular. Denote
\begin{equation}\label{eq:uv}
 u_R(t)=\Rea E_\Psi(Re^{it}),\qquad
 v_R(t)=\Ima E_\Psi(Re^{it}),
\end{equation}
\begin{equation}\label{eq:abR}
 a_R=\frac\kappa2(R+R^{-1}),\qquad
 b_R=\frac\kappa2(R-R^{-1}).
\end{equation}
Here $E_{\Psi} = \Psi(s(z))-iw(z)\Psi'(s(z))$. From $i\kappa s=\frac\kappa2(z-z^{-1})$,
we obtain the exact physical outer boundary parametrizations
\begin{equation}\label{eq:outertraces}
 \begin{aligned}
 \gamma_{\rm out}^{\Omega}(t)=X(Re^{it})
    &=\bigl(R\cos t,\ a_R\sin t+v_R(t)\bigr),\\
 \gamma_{\rm out}^{\Lambda}(t)=Y(Re^{it})
    &=\bigl(b_R\cos t+u_R(t),\ R\sin t\bigr).
 \end{aligned}
\end{equation} 
A direct differentiation gives the curvature numerators
\begin{equation}\label{eq:outercurv}
\begin{aligned}
 K_{\rm out}^{\Omega}(t)
 &=R\bigl[a_R+\cos t\,v_R'(t)-\sin t\,v_R''(t)\bigr],\\
 K_{\rm out}^{\Lambda}(t)
 &=R\bigl[b_R-\sin t\,u_R'(t)-\cos t\,u_R''(t)\bigr].
 \end{aligned}
\end{equation}
Therefore for $\kappa$ large both of them are strictly positive. The area for the domains can be calculated directly using Green's formula. 
\end{proof}

Using the equation \eqref{eq:innertraces} we also give an exact description of $K_{\rm hole}$ and the mapping $T$ restricted to it.

\begin{lem}\label{lem:hole_description}
    The region $K_{\rm hole}$ is foliated by straight lines as
    \begin{equation}\label{eq:holeparam}
 K_{\rm hole}=\left\{(x,y):
 \begin{array}{l}
 y=\kappa q-x\Psi'(q),\\
 -1<q<1,\quad |x|<\sqrt{1-q^2}
 \end{array}\right\}.
\end{equation}
Moreover the mapping $T$ restricted to $\partial K_{\rm hole}$ is given by,
$$T(\pm \sqrt{1-q^2},\kappa q\mp \sqrt{1-q^2}\Psi'(q)) = (\Psi(q),q)$$
for $|q|\leq 1$.
\end{lem}

Note that at the end points of the foliation we have,
$(x,y)_{\pm} = (\pm \sqrt{1-q^2},\kappa q\mp \sqrt{1-q^2}\Psi'(q))$, therefore,
$$(x,y)_{+}-(x,y)_{-} = 2\sqrt{1-q^2}(1,-\Psi'(q)).$$

Therefore the difference of points is also giving us a vector in the normal direction to the graph of $\Psi$ at the point $(\Psi(s),s)$.
\begin{figure}[h]
    \centering
    \includegraphics[width=\linewidth]{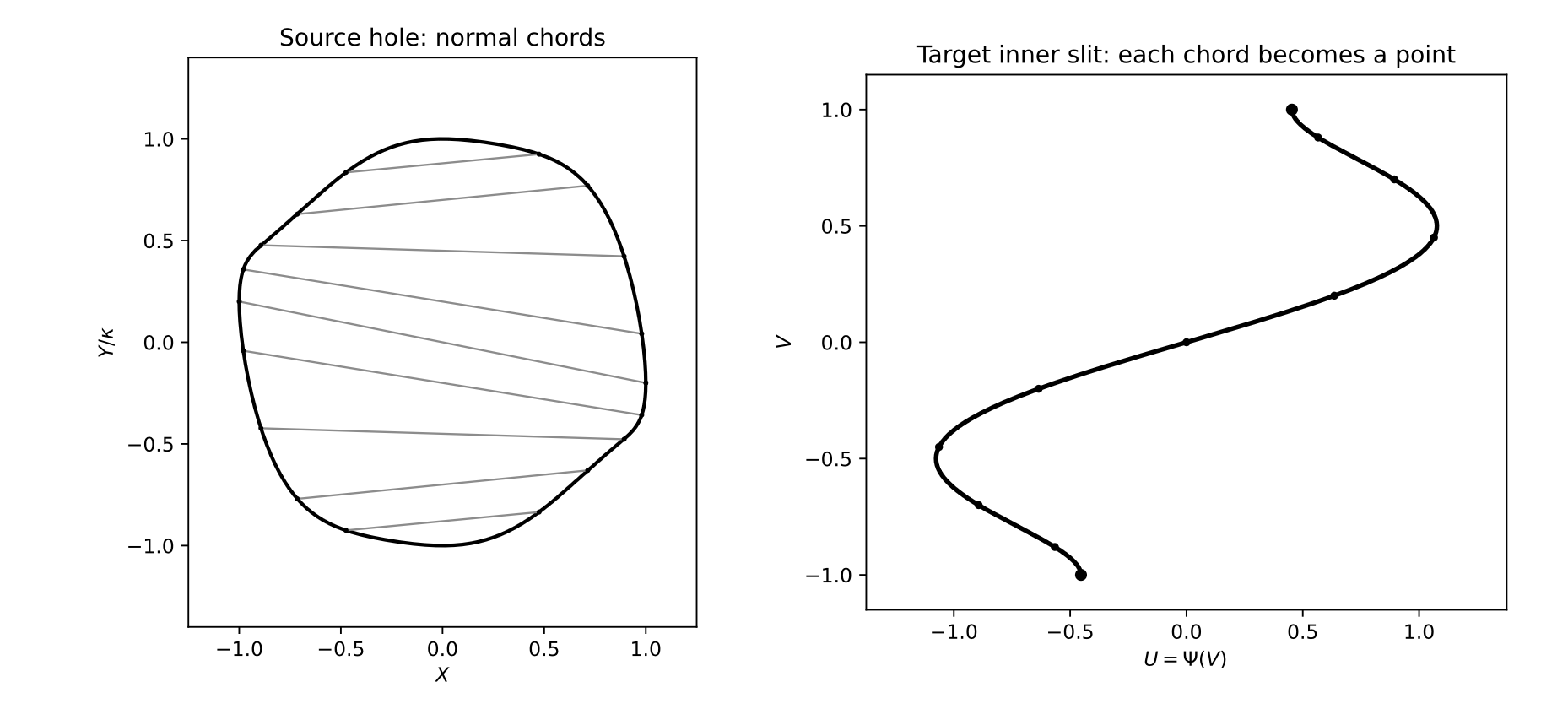}
    \caption{Each line segment in the source (on the left) collapses to a point on the graph (on the right) with the normal to the point given by the direction of the line segment.}
\end{figure}

\subsubsection{Global convexity and Optimal mapping.} We now show that the map we have produced is indeed the gradient of a convex function. 

\begin{lem}
    There exists a function $\phi_{\rm out}$ on $\Omega$ such that $\nabla \phi_{\rm out} = T$.
\end{lem}
\begin{proof}
    By Proposition \ref{prop:periods} we just have to show that the integral,
    $$\Rea \oint_{\partial \mathbb{S}^{1}}f_{\kappa,\Psi}(z)dz = 0.$$
    Note that,
    $$f_{\kappa,\Psi}(z) = \Psi(s(z))
             -i w(z)\Psi'(s(z)) + \frac{1}{2}\kappa (z-z^{-1}) = \frac{d}{dz}[z\Psi(s(z))+\kappa \frac{z^2}{4}]-
             \frac{\kappa}{2} z^{-1}.$$
The first portion is exact therefore its integral is zero and since
$$\oint_{\partial \mathbb{S}^{1}} \frac{1}{z}dz = 2\pi i $$
we get the desired conclusion.
\end{proof}

Inside the domain $K_{\rm hole}$ we construct a convex function $\phi_{\rm in}$ such that it can be glued together with $\phi_{\rm out}$ in a $C^1$ fashion. This shows that $T$ is indeed the optimal transport map.

\begin{lem}\label{lem:T_opt_fbp}
    There exists a convex function $u$ on $K_{\rm out}$ such that $\nabla u = T$. 
\end{lem}
\begin{proof}
    Consider the partial Legendre transform of the function $\kappa \frac{q^2}{2}-x\Psi(q)$,
    $$\phi_{\rm in}(x,y) = \sup_{|q|\leq 1}\left(x\Psi(q)+yq - \kappa \frac{q^2}{2}\right).$$
    note that for $(x,q)\in B_1(0)\subset \C$, one exactly recovers $y = \kappa q-x\Psi'(q)$. Moreover since $\kappa - C\Psi''>0$, we can write $q = q(x,y)$ for $(x,y)\in K_{\rm hole}$. Therefore
    $$\phi_{\rm in}(x,y) = x\Psi(q(x,y))+yq(x,y) - \kappa \frac{q(x,y)^2}{2},$$
    and at $(x,y)\in K_{\rm hole}$ we have,
    $$\nabla \phi_{\rm in}(x,y) = (\Psi(q),q).$$
    Therefore at $(x,y)\in \partial K_{\rm hole}$ we have $\nabla \phi_{\rm in} = T$. 
    Therefore since $\nabla \phi_{\rm in} = \nabla \phi_{\rm out}$ on $\partial K_{\rm hole}$, they only differ by a constant on the boundary and after adjusting for the constant the can be glued together so that the function

    $$u = 
\begin{cases}
    \phi_{\rm in}  & \text{if } x \in \overline{K_{\rm hole}} \\
    \phi_{\rm out} & \text{if } x\in \Omega = K_{\rm out}\setminus \overline{K_{\rm hole}}
\end{cases}$$
    is $C^{1}$ and convex on $K_{\rm out}$. The convexity can be seen by restricting to lines. On $\Omega$, $D^2u = D^2\phi_{\rm out}>0$ and inside $K_{\rm hole}$, $\phi_{\rm in}$ is convex.
    
\end{proof}

The proof of Theorem \ref{thm:general_fbp} follows from Lemmata \ref{lem:psi_global_diff}, \ref{lem:dom_geo_fbp} and \ref{lem:T_opt_fbp}.

\end{document}